\documentclass[12pt]{article}
\usepackage{amssymb}
\usepackage{amsmath}
\usepackage{latexsym}

\newtheorem{theorem}{Theorem}[section]
\newtheorem{proposition}[theorem]{Proposition}

\newtheorem{lemma}[theorem]{Lemma}

\newtheorem{remark}[theorem]{Remark}
\newtheorem{example}[theorem]{Example}
\newtheorem{exercise}[theorem]{Exercise}
\newtheorem{definition}[theorem]{Definition}

\begin{document}

\title{Vertex Operators and Modular Forms}
\author{Geoffrey Mason\thanks{Research supported by NSA and NSF}\\
Michael Tuite}
\maketitle

\tableofcontents

\newpage
\part{Correlation Functions and Eisenstein Series}\label{PartI} 

\section{The Big Picture}
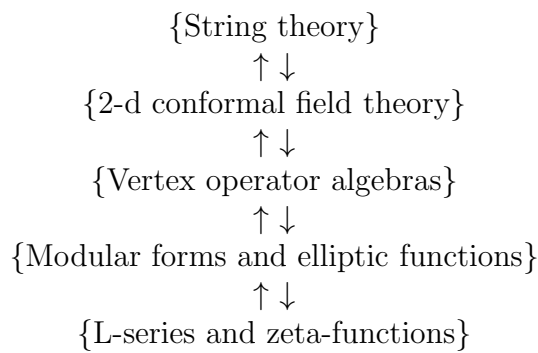
\label{Section_Intro}
 \begin{eqnarray*}
\begin{array}{c}
\{\mbox{String   theory} \} \\
\uparrow  \ \downarrow \\
\{\mbox{2-d  conformal field theory}\} \\
\uparrow \ \downarrow \\
\{\mbox{Vertex operator algebras}\}\\
\uparrow  \ \downarrow \\
\{\mbox{Modular forms and elliptic functions}\}\\
\uparrow  \ \downarrow \\
\{\mbox{L-series and zeta-functions}\}
\end{array}
 \end{eqnarray*}

 \bigskip
The \emph{leitmotif} of these Notes is the idea of a \emph{vertex operator
algebra} (VOA) and the relationship between VOAs  and \emph{elliptic functions 
and modular forms}.  This is to some extent
analogous to the relationship between a finite group and its  irreducible characters; the algebraic structure
determines a set of numerical invariants, and arithmetic properties of the invariants provides feedback in the form of restrictions on the algebraic structure. One of the main points of these Notes is to explain how
this works, and to give some reasonably interesting examples.

\medskip
VOAs may be construed as an axiomatization of
$2$-dimensional conformal field theory, and it is via this connection that vertex operators enter into physical theories. A sketch of the VOA-CFT connection via the Wightman Axioms can be found in the Introduction to \cite{K1}. Although we make occasional comments to relate our development of VOA theory to physics, no technical expertise in physics is necessary to understand these Notes. As mathematical theories go, the one we are discussing here is relatively new.  There are a number of basic questions which are presently unresolved, and we will get far enough in the Notes to explain some of them.

\medskip
To a modular form one may attach (via the Mellin transform) a Dirichlet series, or $L$-function, and Weil's Converse Theorem says that one can go the other way too. So there is a close connection
between modular forms and certain $L$-functions, and this is one way in which our subject matter 
relates to the contents of other parts of this book. Nevertheless, as things stand at present,
it is the Fourier series of a modular form, rather than its Dirichlet series, that is important in VOA theory.
As a result, $L$-functions will not enter into our development of the subject.

\medskip
The Notes are divided into three parts. In Part~I we give some of the foundations of VOA theory,
and explain how modular forms on the full modular group (Eisenstein series in particular) and elliptic functions naturally intervene in the description of $n$-point  correlation functions. This is a general phenomenon and the simplest VOAs, namely the free boson (Heisenberg VOA) and the Virasoro VOA, suffice to illustrate the computations. For this reason we delay the introduction of more complicated VOAs until Part~II, where we describe several families of VOAs and their representations. We also cover some aspects of vector-valued modular forms, which is the appropriate language to describe the modular properties
of $C_2$-cofinite and rational VOAs. We give some applications to holomorphic VOAs to illustrate how
modularity impinges on the algebraic structure of VOAs.
In Part~III we describe two current areas of active research of the authors. The first concerns the development of VOA theory on a genus two Riemann surface and the second is concerned  with the  relationship between exceptional VOAs and Lie algebras and the Virasoro algebra.  

\medskip
There are a number of Exercises at the end of each Subsection.  They provide both practice in the ideas and also a subtext to which we often refer during the course of the Notes. Some of the Exercises are straightforward, others less so... Even if the reader is not intent on working out the Exercises, he or she should read them over before proceeding.

\medskip
These Notes constitute an expansion of the lectures we gave at MSRI in the summer of 2008 during the Workshop \emph{A window into zeta and modular physics}. We would like to thank the organizers of the 
Workshop, in particular Klaus Kirsten and Floyd Williams, for giving us the opportunity to participate in the program.
 
\section{Vertex Operator Algebras}\label{Section_VOAs}
\subsection{Notation and Conventions}\label{Subsect_Notation}

$\mathbb{Z}$ is the set of integers, $\mathbb{R}$ the real
numbers, $\mathbb{C}$ the complex numbers, $\mathfrak{H}$ the complex upper half-plane
\begin{eqnarray*}
\mathfrak{H} = \{ \tau \in \mathbb{C} \ | \ \Im(\tau) > 0 \}.
\end{eqnarray*}

\noindent
All linear spaces $V$ are defined over $\mathbb{C}$; linear transformations are
$\mathbb{C}$-linear; End$(V)$ is the space of \emph{all} endomorphisms of $V$. 
 For an indeterminate $z$, 
\begin{eqnarray*}
V[[z, z^{-1}]] &=& \left\{ \sum_{n \in \mathbb{Z}} v_nz^n \ | \ v_n \in V \right\},  \\
V[[z]][z^{-1}] &=&  \left\{ \sum_{n = -M}^{\infty}  v_nz^n \ | \ v_n \in V \right\}.
\end{eqnarray*}
These are linear spaces with respect to the obvious addition and scalar multiplication. The  formal residue is
\begin{eqnarray*}
\mathrm{Res}_z \sum_{n \in \mathbb{Z}} v_n z^n =  v_{-1}.
\end{eqnarray*}
For integers $m, n$ with $n \geq 0$, 
\begin{eqnarray*}
\binom{m}{n} = \frac{m(m-1)\hdots (m-n+1)}{n!}.
\end{eqnarray*}
For indeterminates $x,y$ we adopt the convention that
\begin{eqnarray}
(x+y)^m=\sum_{n\ge 0}\binom{m}{n} x^{m-n}y^n, \label{xypower}
\end{eqnarray}
i.e. for $m<0$ we formally expand in the second parameter $y$. 
  
  \medskip
We use the following \emph{$q$-convention}:
\begin{eqnarray}\label{qconv}
q_x = e^x; q = q_{2\pi i \tau} = e^{2\pi i \tau} \ (\tau \in \mathfrak{H}),
\end{eqnarray}
where $x$ is anything for which $e^x$ makes sense.

\subsection{Local Fields}\label{Subsect_localfields}
We deal with  formal series 
\begin{eqnarray}\label{psdef}
a(z) = \sum_{n \in \mathbb{Z}} a_nz^{-n-1} \in \mbox{End}(V)[[z, z^{-1}]].
\end{eqnarray}
$a(z)$ defines a linear map
\begin{eqnarray*}
&&a(z): V \rightarrow V[[z, z^{-1}]] \\
&& \ \ \ \ \ \ \ \ \ v \mapsto \sum_{n \in \mathbb{Z}} a_n(v)z^{-n-1}.
\end{eqnarray*}
The endomorphisms $a_n$ are called the \emph{modes} of $a(z)$. We often refer to 
the elements in $V$ as \emph{states}, and call $V$ the \emph{state-space} or \emph{Fock space}.

\begin{remark}\label{Rem_zpowers}
The convention for powers of $z$ in (\ref{psdef}) is standard in mathematics. A different  convention is common in the physics literature. Whenever a mathematician and physicist discuss fields, they should first agree on their conventions. 
\end{remark}
\begin{definition}\label{Def_field}
 $a(z) \in$ End$(V)[[z, z^{-1}]]$  is a \emph{field} if it satisfies the following truncation condition $\forall\ v \in V$:
\begin{eqnarray*}a(z)v \in V[[z]][z^{-1}].
\end{eqnarray*}
 I.e. for $v \in V$
there is an integer $N$ (depending on $v$) such that $a_n(v) = 0$ for all $n > N$.
\end{definition}
Set
\begin{eqnarray*}
\mathfrak{F}(V) = \{ a(z) \in \mbox{End}(V)[[z, z^{-1}]] \ | \ a(z) \ \mbox{is a field}\}.
\end{eqnarray*}
$\mathfrak{F}(V)$ is the field-theoretic analog of End$(V)$. It is a subspace of
End$(V)[[z, z^{-1}]].$

\bigskip

The introduction of a \emph{second} indeterminate facilitates 
the study of products and commutators of fields. Set
\begin{eqnarray*}
\left[\sum_m a_mz_1^{-m-1}, \sum_n b_nz_2^{-n-1}\right] = \sum_{m, n} [a_m, b_n]z_1^{-m-1}z_2^{-n-1},
\end{eqnarray*}
which lies in $\mbox{End}(V)[[z_1, z_1^{-1}, z_2, z_2^{-1}]].$
The idea of \emph{locality} is crucial. 
\begin{definition} $a(z), b(z) \in \mbox{End}(V)[[z, z^{-1}]]$ are called \emph{mutually local}
if there is a nonnegative integer $k$ such that 
\begin{eqnarray}\label{local1def}
(z_1-z_2)^k[a(z_1), b(z_2)] = 0.
\end{eqnarray}
\end{definition}

If (\ref{local1def}) holds, we write $a(z) \sim_k b(z)$ and say that $a(z)$ and $b(z)$ are \emph{mutually local of order $k$}. Write 
$a(z)\sim b(z)$ if $k$ is not specified. $a(z)$ is  a \emph{local field} if $a(z) \sim a(z)$. 
(\ref{local1def}) means that the coefficient of each monomial $z_1^{r-1} z_2^{s-1}$ in the expansion of the left hand side  vanishes. Explicitly, this means that
\begin{eqnarray}\label{localid}
\sum_{j = 0}^k (-1)^j\binom{k}{j}[a_{k-j-r}, b_{j-s}] = 0.
\end{eqnarray}
Locality defines a \emph{symmetric relation} which is generally neither
reflexive nor transitive.

\bigskip
Fix a \emph{nonzero} state $\mathbf{1} \in V$. We say that $a(z) \in \mathfrak{F}(V)$
is \emph{creative} (with respect to $\mathbf{1}$) and \emph{creates the state $u$}  if 
\begin{eqnarray*}
a(z)\mathbf{1} = u+ \hdots \in V[[z]].
\end{eqnarray*}
We sometimes write this in the form $a(z)\mathbf{1} = u+O(z)$. In terms of modes, 
\begin{eqnarray*}
a_n\mathbf{1} = 0, \ n \geq 0, \ a_{-1}\mathbf{1} = u.
\end{eqnarray*}

\bigskip
\noindent
\begin{exercise}\label{Exercise_2.2.1}
 Let $\partial a(z) = \sum_n (-n-1)a_nz^{-n-2}$ be the formal derivative
of $a(z)$. Suppose that $a(z), b(z) \in \mathfrak{F}(V)$ and $a(z) \sim_k b(z)$. Prove that
$\partial a(z) \in \mathfrak{F}(V)$ and $\partial a(z) \sim_{k+1} b(z)$.
\end{exercise}

\begin{exercise}\label{Exercise_local_Trunc}(Locality-Truncation Relation) 
Suppose that $a(z), b(z)$ are creative fields with $a(z) \sim_k b(z)$. By choosing $s=1$ and $r=k-n$ for $n\ge k$ in (\ref{localid}), show that $a_n b=0$ for all $n\ge k$ i.e. the order of truncation $N$ is $k-1$. 
\end{exercise}

\subsection{Axioms for a Vertex Algebra}\label{Subsect_VAaxioms}
For various approaches to the contents of this Subsection, see
 \cite{B}, \cite{FHL}, \cite{FLM}, \cite{Go},  \cite{K1}, \cite{LL}, \cite{MN}. 
 \begin{definition} A \emph{vertex algebra} (VA) is a quadruple $(V, Y, \mathbf{1}, D)$ where 
\begin{eqnarray*} 
&&Y: V \rightarrow \mathfrak{F}(V), \ v \mapsto Y(v, z) = \sum v_n z^{-n-1} \ \mbox{is a linear map}, \\
&& \mathbf{1} \in V, \ \mathbf{1} \not= 0, \\
&&D \in \mbox{End}(V),  \ D\mathbf{1} = 0,
\end{eqnarray*}
and the following hold for all $u, v \in V$:
\begin{eqnarray*}
&&\ \ \ \ \ \ \ \ \ \ \ \ \ \ \ \ \ \ \mbox{locality}: Y(u, z) \sim Y(v, z) \\
&&\ \ \ \ \ \ \ \ \ \ \ \ \ \ \   \mbox{creativity}: Y(u, z)\mathbf{1} = u + O(z)\\
&&\mbox{translation-covariance}: [D, Y(u, z)] = \partial Y(u, z)
\end{eqnarray*}
\end{definition}

We often refer to the Fock space $V$ itself as a vertex algebra rather than $(V, Y, \mathbf{1}, D)$. 
$\mathbf{1}$
is called the \emph{vacuum} state and $Y$  the \emph{state-field correspondence}. 
The physical interpretation of creativity is that $Y(u, z)$ \emph{creates}
the state $u$ from the vacuum. 
This set-up models the creation and annihilation of bosonic states from the vacuum. Most of the subtlety is tied to locality and its consequences.

\medskip
 There are a number of equivalent formulations of these axioms. We discuss some of them. Another approach, via so-called \emph{rationality} (\cite{FHL}) is also discussed in Subsection \ref{Subsect_MatrixElements}. The \emph{Jacobi Identity} of 
 \cite{FLM} is equivalent to the identity
\begin{eqnarray}\label{BI}
\sum_{i \geq 0} \binom{p}{i} (u_{r+i}v)_{p+q-i} = \sum_{i \geq 0} (-1)^i \binom{r}{i}  \left\{ u_{p+r-i}v_{q+i} -(-1)^r v_{q+r-i}u_{p+i} \right\},
\end{eqnarray}
which holds in a VA for all $u, v \in V$ and all $p, q, r \in \mathbb{Z}$. 
Conversely, if $Y(v, z) \in \mathfrak{F}(V), v \in V,$ are creative
fields with respect to $\mathbf{1}$ and satisfy (\ref{BI}) then $(V, Y, \mathbf{1}, D)$ is a vertex algebra
with $Du = u_{-2}\mathbf{1}$. 

\medskip
Specializing (\ref{BI}) in various ways leads to some particularly useful identities first written down in \cite{B}:
\begin{eqnarray*} 
{[}u_m, v_n{]} &=& \sum_{i \geq 0} \binom{m}{i} (u_iv)_{m+n-i} \ (\mbox{commutator}), \\
(u_mv)_n  &=&\sum_{i \geq 0} (-1)^i \binom{m}{i}\{ u_{m-i}v_{n+i} - (-1)^mv_{m+n-i}u_i \} \ 
(\mbox{associator}),  \\
u_mv &=& \sum_{i \geq 0}(-1)^{m+i+1}\frac{1}{i!} D^iv_{m+i}u \ (\mbox{skew-symmetry}).
\end{eqnarray*}
These identities may be stated more compactly using vertex operators, and it is often more efficacious to
use the vertex operator format. We state one more consequence of (\ref{BI}),
the \emph{associativity} formula, in the operator format. For large enough $k$ (and recalling convention (\ref{xypower}) 
\begin{eqnarray}\label{associativity}
(z_1+z_2)^kY(u, z_1+z_2)Y(v, z_2)w = (z_1+z_2)^kY(Y(u, z_1)v, z_2)w.
\end{eqnarray}
We have (\cite{FKRW}, \cite{MP})
\begin{theorem}\label{thmexist} Let $V$ be a linear space with $0 \not= \mathbf{1} \in V$
and $D \in$ End$(V)$. Suppose $S \subseteq \mathfrak{F}(V)$ is a set of mutually local, creative, translation-covariant fields which \emph{generates} $V$ in the sense that
\begin{eqnarray*}
V = \mbox{span}\{a^1_{-n_1} \hdots a^k_{-n_k} \mathbf{1} \ | \ a^i(z) \in S, n_1, \hdots , n_k \geq 1, k \geq 0 \}.
\end{eqnarray*}
Then there is a unique vertex algebra $(V, Y, \mathbf{1}, D)$ such that 
$Y(a^i_{-1}\mathbf{1}, z) = a^i(z)$.
\end{theorem}

\bigskip

\begin{exercise}\label{Exercise_2.3.1}
 Prove that the state-field correspondence is \emph{injective}.
\end{exercise}

\begin{exercise}\label{Exercise_2.3.2} Prove that
\begin{eqnarray*}
Y(u, z)\mathbf{1} = q_z^{D}u \ (= \sum_{n \geq 0} \frac{z^n}{n!} D^nu).
\end{eqnarray*}
\end{exercise}

\begin{exercise}\label{Exercise_2.3.3} Deduce the commutator, associator and skew-symmetry formulas from
(\ref{BI}).
\end{exercise}

\begin{exercise}\label{Exercise_2.3.4} Assume $V$ is a linear space and $\{Y(v, z) \ | \ v \in V \} \subseteq \mathfrak{F}(V)$ are mutually local fields such
that $Y(v, z)$ is creative (with respect to $\mathbf{1} \not= 0$) and creates
$v$. Prove that (\ref{BI}) and the associator formula are \emph{equivalent}.
\end{exercise}

\begin{exercise}\label{Exercise_2.3.5} $A$ is a commutative, associative algebra with identity element
$1$ and derivation $D$.  Show that there is a vertex algebra $(A, Y, 1, D)$ with
$Y(a, z)b = \sum_{n\leq -1} \frac{(D^{-n-1}a)b}{(n+1)!}z^{-n-1}.$
\end{exercise}

\begin{exercise}\label{Exercise_2.3.6} $(V, Y, \mathbf{1}, D)$ is a VA. Assume either (a) $Y(v, z) \in$ End$(V)[[z]]$
for $v \in V$, (b) $D$ is the zero map, or (c) $\dim V$ is finite.  Prove in each case
that $V$ is of the type described in Exercise \ref{Exercise_2.3.5}.
\end{exercise}

\begin{exercise}\label{Exercise_2.3.7} Show that the commutator formula is equivalent to
the identity $[u_m, Y(v, z)] = \sum_{i \geq 0} \binom{m}{i} Y(u_iv, z)z^{m-i}$.
\end{exercise}

\begin{exercise}\label{Exercise_skew} Show that the skew-symmetry formula is equivalent to
the identity $Y(u,z)v=  q_z^D Y(v,-z)u$.
\end{exercise}

\begin{exercise}\label{Exercise_translation} Show that $q_y^{D}Y(u,x)q_y^{-D}=Y(u,x+y)$. 
\end{exercise}%

\subsection{Heisenberg Algebra}\label{Subsect_HeisVOA}
In this and the following Subsection we will use Theorem \ref{thmexist} to construct two fundamental examples
of VAs. We must look for generating sets $S$ of mutually local, creative, translation-covariant fields. In our two examples, $S$ consists of a \emph{single} field. The construction relies on some basic techniques from Lie theory (Verma modules, Poincar\'{e}-Birkhoff-Witt Theorem, etc) which are reviewed in the Appendix.

\medskip
Let $A = \mathbb{C}a$ be a $1$-dimensional linear space. The \emph{affine algebra} 
$\hat{A} = A[t, t^{-1}] \oplus \mathbb{C}K$
is the Lie algebra with central element $K$ and bracket
\begin{eqnarray}\label{aff1}
[a \otimes t^m, a \otimes t^n] = m \delta_{m, -n}K. 
\end{eqnarray}
\begin{remark}\label{Remark_2.7} Set $p_m = \frac{1}{\sqrt{m}} a \otimes t^m$($m > 0$) and $q_{-m} = 
\frac{1}{\sqrt{-m}} a \otimes t^m$($m < 0.$) Then (\ref{aff1}) reads
\begin{eqnarray}\label{CCR}
[p_m, q_n] = \delta_{m, n} K.
\end{eqnarray}
These are essentially the \emph{canonical commutator relations} of Quantum Mechanics. 
\end{remark}

\medskip
 Set $\hat{A}^{\geq} = \langle a\otimes t^n, K \ | \ n \geq 0 \rangle$,
 $\hat{A}^- =  \langle a\otimes t^n \ | \ n<0\}$. These are a Lie ideal and Lie subalgebra
 of $\hat{A}$ respectively. Let $\mathbb{C}v_h$ be the $1$-dimensional 
$\hat{A}^{\geq}$-module defined for a scalar $h$ via
\begin{eqnarray*}
K.v_h &=& v_h, \\
 (a\otimes t^n).v_h &=& h\delta_{n, 0}v_h \ ( n \geq 0).
  \end{eqnarray*}
 The induced (Verma) module is
\begin{eqnarray}\label{MVerma}
M_h = \mbox{Ind}_{\mathcal{U}(\hat{A}^{\geq})}^{\mathcal{U}(\hat{A})} \mathbb{C}v_{h}
=  \mathcal{U}(\hat{A}) \otimes_{\mathcal{U}(\hat{A}^{\geq})} \mathbb{C}v_h = \mathcal{U}(\hat{A}^-) \otimes \mathbb{C}v_h
\end{eqnarray}
where $\mathcal{U}( \ )$ denotes universal enveloping algebra and the third equality in the last display is just a linear isomorphism. 

\medskip
Let $a_n \in \mbox{End}(M_h)$ be the induced action
of $a \otimes t^n$ on $M_h$, with $a(z) = \sum_n a_nz^{-n-1}$. In what follows we identify $v_h$ with $1 \otimes v_h$.
Let $v = a_{-n_1} \hdots a_{-n_k}.v_h$ with $n_1\geq \dots \geq n_k \geq 1.$
For $n > n_1, a_n$ commutes with each $a_{-n_i}$ by (\ref{aff1}). Therefore, 
$a_n.v = a_{-n_1} \hdots a_{-n_k}a_n.v_h = 0$. This shows that
$a(z) \in \mathfrak{F}(M_h)$. As for locality,
\begin{eqnarray}\label{locality1}
&& \sum_{j = 0}^2 (-1)^j\binom{2}{j}[a_{2-j-r}, a_{j-s}] \notag \\
&=& \sum_{j = 0}^2 (-1)^j\binom{2}{j}(2-j-r) 
\delta_{2-j-r, s-j}K   \\
&=& \left\{ (2-r)
-2(1-r) -r  \right \} \delta_{r+s, 2}K = 0. \notag
\end{eqnarray}
By (\ref{localid}) this shows that $a(z) \sim_2 a(z)$. Because
\begin{eqnarray*}
a(z)v_h = hv_hz^{-1} + \sum_{n \leq -1} a_n.v_h z^{-n-1},
\end{eqnarray*}
we see that $a(z)$ is creative with respect to $v_h$ if (and only if) $h=0$.
  In this case, Theorem \ref{thmexist} and what we have shown imply  
\begin{theorem}\label{thm2.7} There is a unique vertex algebra $(M_0, Y, v_0, D)$
generated by $a(z)$ with $Y(a, z) = a(z)$ and $Da_n.v_0 = -na_{n-1}v_0$. 
\end{theorem}

\noindent
\begin{remark} \label{Remark_Mboson} In terms of operators on $M_0$, (\ref{CCR}) reads
$[p_m, q_n] = \delta_{m, n}$Id. These relations may be realized by taking $p_m = \frac{\partial}{\partial x_{-m}},
q_n = x_{-n}$ acting on the Fock space $\mathbb{C}[x_{-1}, x_{-2}, \hdots]$. This affords an alternate way to understand $M_0$.
\end{remark}

$M_0$ is variously called the (rank $1$) Heisenberg VA, Heisenberg algebra, or 
\emph{free boson}. In CFT it models a single free boson. (As opposed to standard mathematical usage,
 \emph{free} here means that the particle is not interacting with other particles.)

\subsection{Virasoro Algebra} \label{Subsect_VirasoroVOA}
The Virasoro algebra is the Lie algebra with underlying linear space
\begin{eqnarray*}
\mbox{Vir} = \bigoplus_{n \in \mathbb{Z}} \mathbb{C}L_n \oplus \mathbb{C}K
\end{eqnarray*}
and bracket relations
\begin{eqnarray}\label{Virrelns}
[L_m, L_n] = (m-n)L_{m+n} + \frac{m^3 - m}{12}\delta_{m, -n}K.
\end{eqnarray}
Set $\mbox{Vir}^{\geq} = \langle L_n, K \ | \ n \geq 0 \rangle, \ \mbox{Vir}^- = \langle L_n \ | \ n < 0 \rangle$, and let
$\mathbb{C}v_{c, h}$ be the $1$-dimensional $\mbox{Vir}^{\geq}$-module defined via
\begin{eqnarray*}
K.v_{c, h} &=& cv_{c, h}, \\
L_n.v_{c, h} &=& \delta_{n, 0}hv_{c, h} \ (n \geq 0).
 \end{eqnarray*}
 with arbitrary scalars $c, h$. 
The induced (Verma) module is then
\begin{eqnarray}\label{VirVerm}
M_{c, h} = \mathcal{U}(\mbox{Vir}) \otimes_{\mathcal{U}(\mbox{Vir}^{\geq})} \mathbb{C}v_{c, h} 
= \mathcal{U}(\mbox{Vir}^-)\otimes \mathbb{C}v_{c, h}.
\end{eqnarray}

\medskip
By analogy with Theorem \ref{thm2.7}, Exercise \ref{Exercise_2.5.2} (below) suggests that there is a VA with Fock space $M_{c, 0}$ and vacuum\footnote{As in the case of the Heisenberg algebra, we identify $v_{c, 0}$ and $1 \otimes v_{c, 0}$.}  $v_{c, 0}$, with $L_{-1}$ playing the r\^{o}le of $D$. This cannot be true as it stands because $\omega(z).v_{c, 0} = L_{-1}. v_{c, 0}z^{-1} + \hdots$ is \emph{not} creative.
 To cure this ill requires (at the very least) that we take a quotient of $M_{c, 0}$ by a $\mbox{Vir}$-submodule that contains $L_{-1}. v_{c, 0}$, and indeed it suffices to quotient out the cyclic $\mbox{Vir}$-submodule generated by this state.  We will abuse notation by identifying states, operators and fields associated with 
 $M_{c, 0}$ with the corresponding states, operators and fields induced on the quotient
 $M_{c, 0}/\mathcal{U}(\mbox{Vir})L_{-1}.v_{c, 0}$. We then arrive at

\begin{theorem}\label{VirVA} Set $\mbox{Vir}_c = M_{c, 0}/ \mathcal{U}(\mbox{Vir})L_{-1}.v_{c, 0}$,
$\omega = L_{-2}. v_{c, 0}$, and $Y(\omega, z) = \omega(z)$. Then $(\mbox{Vir}_c, Y, v_{c, 0}, L_{-1})$ is a vertex algebra generated by $Y(\omega, z)$.
\end{theorem}

\medskip
\noindent
$\mbox{Vir}_c$ is called the Virasoro VA of \emph{central charge} $c$.

\bigskip
\begin{exercise}\label{Exercise_2.5.1} 
Show that $L_{-1}, L_0$ and $L_1$ span a Lie subalgebra of
$\mbox{Vir}$. What Lie algebra is it?
\end{exercise}

\begin{exercise}\label{Exercise_2.5.2}  Identify elements of $\mbox{Vir}$ with the endomorphisms they
induce on $M_{c, h}$ and set $\omega(z) = \sum L_nz^{-n-2} \in \mbox{End}(M_{c, h})[[z, z^{-1}]]$. Prove that $\omega(z)$ is
a local field of order $4$, 
and $[L_{-1}, \omega(z)] = \partial \omega(z)$. 
\end{exercise}

\begin{exercise}\label{Exercise_2.5.3} Give the details of the proof of Theorem \ref{VirVA}.
\end{exercise}

\subsection{Axioms for a Vertex Operator Algebra}\label{Subsect_VOAaxioms}
There is no consensus as to nomenclature for the many variants 
of vertex algebra. Our definition of \emph{vertex operator algebra} (VOA) is the one used by many practitioners of the art, but not all.
\begin{definition} A VOA is a quadruple $(V, Y, \mathbf{1}, \omega)$ where 
$V = \oplus_{n \in \mathbb{Z}} V_n$ is a $\mathbb{Z}$-graded linear space and 
\begin{eqnarray*} 
&&Y: V \rightarrow \mathfrak{F}(V), \ v \mapsto Y(v, z) = \sum v_n z^{-n-1}\\
&& \mathbf{1}, \omega \in V, \ \mathbf{1} \not= 0. 
\end{eqnarray*}
The fields $Y(v, z)$ are mutually local and creative,
and the following hold:
\begin{eqnarray*}
&&Y(\omega, z) = \sum L_nz^{-n-2} \ \mbox{with a constant $c$ such that} \\
&&\ \ \ \  [L_m, L_n] = (m-n)L_{m+n} + \frac{m^3 - m}{12}\delta_{m, -n}c\mathrm{Id}_V \\
&&V_n = \{ v \in V_n \ | \ L_0v = nv \}\\
&&\dim V_n < \infty, \ V_n = 0 \ \mbox{for} \ n \ll 0\\
&&Y(L_{-1}u, z) = \partial Y(u, z) 
\end{eqnarray*}
\end{definition}

\medskip
In effect, a VOA is a vertex algebra with a dedicated Virasoro field.
This is the field determined by the distinguished state $\omega$, called the \emph{conformal} or \emph{Virasoro vector}.
 The modes of $\omega$ are operators 
$L_n$ satisfying the Virasoro relations (\ref{Virrelns}) with $K = c \mathrm{Id}_V$. As in Theorem \ref{VirVA} we call $c$ the \emph{central charge} of $V$. The mode $L_0$ of $\omega$, called the \emph{degree operator}, is required to be
semisimple, to have eigenvalues lying in a subset of $\mathbb{Z}$ that is bounded below,  
and to have
finite-dimensional eigenspaces. We often write $wt(v)=n$ if $v$ is an eigenstate for $L_0$ with eigenvalue $n$. It is not hard to see that
$[L_{-1}, Y(u, z)] = \partial Y(u, z)$, so that $(V, Y, \mathbf{1}, L_{-1})$ is a vertex algebra.

\medskip
It should come as no surprise that the vertex algebra $\mbox{Vir}_c$  has the structure
of a VOA of central charge $c$ with vacuum vector $v_{c, 0}$ and conformal vector $\omega$. To see this, note that $\{ L_{-n_1} \hdots L_{-n_k}. v_{c, 0} \ | n_1 \geq \hdots \geq n_k \geq 2\}$ is a basis of the Fock space. We have
 \begin{eqnarray*}
L_0. L_{-n_1} \hdots L_{-n_k}.v_{c, 0} &=& n_1L_{-n_1} \hdots L_{-n_k}.v_{c, 0}\\
&&+ L_{-n_1}L_0L_{-n_2} \hdots L_{-n_k}. v_{c, 0}.
\end{eqnarray*}
Now an easy induction shows that 
\begin{eqnarray*}
L_0. L_{-n_1} \hdots L_{-n_k}. v_{c, 0} = (\sum_i n_i).
L_{-n_1} \hdots L_{-n_k}. v_{c, 0}
\end{eqnarray*}
so that
\begin{eqnarray*}
wt(L_{-n_1}\hdots L_{-n_k}.v_{c, 0}) = \sum_i n_i.
\end{eqnarray*}
 The needed properties of $L_0$ required for the next result follow easily, and we obtain the following extension of Theorem \ref{VirVA}:
\begin{theorem}\label{VirVOA} $\mbox{Vir}_c$ is a VOA of central charge $c$.
\end{theorem}

\medskip
A pair of VOAs $V, V'$ are called \emph{isomorphic} if there is a linear isomorphism
$\varphi:V \rightarrow V', v \mapsto v'$ such that $\varphi(\omega) = \omega',$
and $fY(v, z)f^{-1} = Y(v, z)$.

\bigskip

In the following Exercises, $V$ is a VOA.
\medskip

\begin{exercise} \label{Exercise_2.6.1} Complete the proof of Theorem \ref{VirVOA}.
\end{exercise}

\begin{exercise} \label{Exercise_2.6.2} Prove the following: 
$Y(\mathbf{1}, z) = \mathrm{Id}, \mathbf{1} \in  V_0,  \omega \in V_2, L_n\mathbf{1} = 0 \ \mbox{for} \ n \geq -1, (L_{-1}v)_n = -nv_{n-1}$.
\end{exercise}

\begin{exercise} \label{Exercise_2.6.3} Suppose that $v \in V$ satisfies $L_{-1}v=0$. Prove that
$v \in V_0$.
\end{exercise}

\begin{exercise} \label{Exercise_2.6.4} Suppose that $V_0 = \mathbb{C}\bf{1}$
(cf. Exercise \ref{Exercise_2.6.2}). Prove that $V_n = 0$ for $n < 0$.
\end{exercise}

\begin{exercise} \label{Exercise_2.6.5} Show that $\dim V$ is \emph{finite} if, and only if, $\omega = 0$
(cf. Exercise \ref{Exercise_2.3.6}).
\end{exercise}

\begin{exercise} \label{Exercise_2.6.6} Show that the Heisenberg theory $M_0$ (cf. Subsection \ref{Subsect_HeisVOA}) is a VOA
with vacuum $\mathbf{1}= v_{0}, \omega = \frac{1}{2} a_{-1}^2\mathbf{1} =  \frac{1}{2}a_{-1}a$ and central charge 
$c= 1$. 
This is the  theory of \emph{one free boson}.
\end{exercise}

\begin{exercise} \label{Exercise_2.6.7} Let $U$, $V$ be linear spaces. Show that there is a natural injection
$\mathfrak{F}(U) \otimes \mathfrak{F}(V) \rightarrow \mathfrak{F}(U \otimes V)$. Suppose in addition that
$U$ and $V$ are Fock spaces for VOAs with vacuum vectors $\bf{1}, \bf{1}'$ and conformal vectors $\omega, \omega'$ respectively.
Show how to construct the \emph{tensor product VOA} $(U \otimes V, Y, \bf{1}\otimes \bf{1}', \omega
\otimes \omega')$. What is the central charge of this VOA?
\end{exercise}

\begin{exercise} \label{Exercise_2.6.8} Let $\varphi: V \rightarrow V'$ be an isomorphism of VOAs. Prove the following:
(i) $V$ and $V'$ have the \emph{same} central charge; (ii) $\varphi(\mathbf{1}) = \mathbf{1}'$.
\end{exercise}

\subsection{VOAs on the Cylinder and the Square Bracket Formalism}\label{Subsect_SquareVOA}
 There is a sense in which we may think of a VOA 
 as being `on the sphere'. This is closely related to the axiomatic approach via rationality (cf. \cite{FHL} and Subsection \ref{Subsect_MatrixElements}). Here we want to describe the corresponding VOA
 that lives `on the cylinder'. Roughly, this corresponds to a change of variable
 $z \rightarrow q_z-1$ which we call the \emph{square bracket formalism}. The main purpose is to construct vertex operators that are
automatically periodic in $z$ with period $2\pi i$.
Let $V=(V, Y, \mathbf{1}, \omega)$ be a VOA of central charge $c$.
For $v \in V$  introduce\footnote{We write modes in the square bracket formalism as
$v[n]$ rather than $v_{[n]}$.}
\begin{equation}
Y[v,z] =Y(q_{z}^{L_0}v,q_{z}-1) = \sum_{n\in \mathbb{Z}}v[n]z^{-n-1}.
\label{Ysquare}
\end{equation}
Here, $q_z^{L_0}$ is  the operator
\begin{eqnarray}\label{qopformalism}
q_z^{L_0}: V \rightarrow V[[z]], \ v \mapsto q_{kz}v \  (v \in V_k),
\end{eqnarray}
 and our $q$-convention (\ref{qconv}) is in force. Similar expressions will occur frequently in what follows. The $v[n]$ are new
operators on $V$, where for $v \in V_k$ and $m\ge 0$, are given by 
\begin{eqnarray}
v[m] &=&m!\sum\limits_{i\geq m}c(k,i,m)v_i  \label{square1} 
\end{eqnarray}
with
\begin{eqnarray}
\binom{k-1+x}{i} = \sum\limits_{m=0}^{i}c(k, i, m)x^{m}.
\label{square2}
\end{eqnarray}
From (\ref{square1}) and (\ref{square2}) we find for integer $n$ 
\begin{equation}
\sum_{i\geq 0}\binom{n}{i}v_i =\sum\limits_{m\geq 0}\frac{(n+1-k)^{m}}{m!%
}v[m].  \label{square4}
\end{equation}%
These identities are proved in the Appendix.
We also have a new conformal vector 
\begin{equation}\label{sqomega}
\tilde{\omega}=\omega -\frac{c}{24}\mathbf{1}, 
\end{equation}
with corresponding square bracket modes
\begin{eqnarray*}
Y[\tilde{\omega}, z] = \sum_n L[n]z^{-n-2}.
\end{eqnarray*}
In particular, $L[0]$ provides us with an alternative $\mathbb{Z}$%
-grading structure on $V$: 
\begin{eqnarray*}
V &=&\oplus _n V_{[n]}, \\
V_{[n]} &=&\{u\in V \ | \  L[0]u = nu\}.
\end{eqnarray*}%
We write $wt[v] = n$ if $v \in V_{[n]}$.
The following can be proved.

\begin{theorem} The quadruple $(V, Y[ \ , \ ], \mathbf{1}, \tilde{\omega})$ is a VOA of central charge $c$.
 \end{theorem}

Given a VOA $V$, we say that its alter ego $(V, Y[ \ , \ ], \mathbf{1}, \tilde{\omega})$ is  `on the cylinder'.
VOAs on the cylinder play an important r\^{o}le in forging the connections with modular forms.

\begin{example}
\label{ExampleH[]} In the square bracket formalism, the VOA $(M_0, Y[ \ , \ ], \mathbf{1}, \tilde{\omega})$ 
is generated by a state $a$ with $wt[a]=1$. It has a basis of
Fock vectors of the form $a[-n_{1}] \ldots a[-n_k]\mathbf{1}, \ n_1 \geq \ldots \geq n_k \geq 1$ satisfying
\begin{equation*}
\lbrack a[m],a[n]]=m\delta _{m+n,0}\mathrm{Id}.
\end{equation*}
\end{example}

\bigskip
\begin{exercise}\label{Exercise_2.7.1} Show that $L[-1]=L_{-1}+L_0$.
\end{exercise}

\begin{exercise}\label{Exercise_2.7.2} A state $v$ in a VOA $V$ is called \emph{primary} of weight $k$ with respect to the original Virasoro algebra $\{L_n\}$ if, and only if, it satisfies $L_nv = k\delta_{n, 0}v$
for $n \geq 0$. Prove that $v$  is primary of weight $k$ with respect to $\{L_n\}$ if, and only if, it is primary of weight $k$ with respect to $\{L[n]\}$.
\end{exercise}

\begin{exercise}\label{Exercise_2.7.3} Prove the assertions of Example \ref{ExampleH[]} in the more precise form that
$(M_0, Y, \mathbf{1}, \omega)$ and $(M_0, Y[\ , \ ], \mathbf{1}, \tilde{\omega})$
are \emph{isomorphic} Heisenberg VOAs.
\end{exercise}

\section{Modular and Quasimodular Forms}\label{Section_ModForms}
In this Section we compile some relevant background involving \emph{elliptic modular forms}. This is a standard part of analytic number theory, and there are many excellent texts dealing with the subject, e.g. \cite{Kn}, \cite{O}, \cite{Se}, \cite{Sc}. Because it is so central to our cause, we describe what we need here, referring the reader elsewhere for more details and further development.

\subsection{Modular Forms on $SL_2(\mathbb{Z})$}\label{Subsect_SL2ZModForms}
The (homogeneous) modular group is
\begin{eqnarray*}
\Gamma = SL_2(\mathbb{Z}) = \left\{\left(\begin{array}{cc}a & b \\ c & d\end{array}\right) |  \ a, b, c, d \in \mathbb{Z}, ad-bc = 1 \right\},
\end{eqnarray*}
with standard generators
$S = \left(\begin{array}{cc}0 & -1 \\ 1 & 0 \end{array}\right), \  T = \left(\begin{array}{cc}1 & 1 \\ 0 & 1 \end{array}\right)$.
The complex upper half-plane $\mathfrak{H}$ carries a left $\Gamma$-action by M\"{o}bius transformations
\begin{eqnarray}\label{GHaction}
(\gamma, \tau) &\mapsto& \gamma \tau = \frac{a\tau + b}{c\tau + d}, \ \gamma = \left(\begin{array}{cc}a & b \\ c & d\end{array}\right)\in \Gamma.
\end{eqnarray}
In particular, $T: \tau \mapsto \tau +1$ and $S: \tau \mapsto -1/\tau$. For $k \in \mathbb{Z}$, there is a right action of $\Gamma$ on meromorphic functions in $\mathfrak{H}$
given by
\begin{eqnarray}\label{slashaction}
f|_k \gamma (\tau) &=& (c\tau +d)^{-k}f(\gamma \tau).
\end{eqnarray}
A \emph{weak modular form} of weight $k$ on $\Gamma$ is an \emph{invariant} of this action. Thus
$f|_k \gamma(\tau) = f(\tau)$ for $\gamma \in \Gamma$, which amounts to
\begin{eqnarray*}
 f(\tau + 1) &=& f(\tau),\\
 f(-1/\tau) &=& \tau^k f(\tau).
\end{eqnarray*}
By a standard argument the first of these equalities implies that $f(\tau)$ has a \emph{$q$-expansion},
or \emph{Fourier expansion at $\infty$},
\begin{eqnarray}\label{genqexp}
f(\tau) = \sum_{n \in \mathbb{Z}} a_n q^n,
\end{eqnarray}
with constants $a_n$ called the \emph{Fourier coefficients} of $f(\tau)$. Here we are using our $q$-convention (\ref{qconv}).

\medskip
$f(\tau)$ is a \emph{meromorphic modular form of weight $k$} if its
$q$-expansion has the form
\begin{eqnarray}\label{meroqexp}
f(\tau) = \sum_{n \geq n_0} a_n q^n
\end{eqnarray}
for some $n_0$. Assume that $f(\tau)\not= 0$ with $a_{n_0}\not = 0$. We then say that $f(\tau)$ has a pole of order $n_0$ at $\infty$ if $n_0 \leq 0$ or a zero of order $n_0$ if $n_0\geq 0$. In the latter situation we also say that $f(\tau)$ is \emph{holomorphic at $\infty$}.
$f(\tau)$ is a \emph{holomorphic modular form of weight $k$} if it is holomorphic in 
$\mathfrak{H}\cup \{\infty\}$.  $f(\tau)$ is \emph{almost holomorphic} if it is holomorphic in
$\mathfrak{H}$ (the behaviour at $\infty$ being unspecified beyond being at worst a pole). 
Modular forms of weight $0$ are often called \emph{modular functions},  though we will not be consistent on this point. Let $\mathfrak{M}_k$ be the set of holomorphic modular forms of weight $k$. It is a $\mathbb{C}$-linear space, possibly equal to $0$.

\bigskip

\begin{exercise}\label{Exercise_3.1.1} Show that the \emph{kernel} of the $\Gamma$-action (\ref{GHaction}) is the 
\emph{center} of $\Gamma$ and consists of $\pm I$, where $I$ is the $2\times 2$ identity matrix. (The quotient group $PSL_2(\mathbb{Z}) = \tilde{\Gamma}= \Gamma/\{\pm I \}$ is the \emph{inhomogeneous} modular group.)
\end{exercise}

\begin{exercise}\label{Exercise_3.1.2} (a) Show that torsion elements in $\tilde{\Gamma}$
have order \emph{at most $3$}.  (b) Show that $\tilde{\Gamma}$ has a \emph{unique}
conjugacy class of subgroups of order $2$ or $3$.
\end{exercise}

\begin{exercise}\label{Exercise_3.1.3} Let $z \in \mathfrak{H}$ with Stab$_{\tilde{\Gamma}}(z)
= \{ \gamma \in \tilde{\Gamma} \ | \ \gamma.z = z \}$ the \emph{stabilizer} of $z$ in $\tilde{\Gamma}$.
Prove the following: (a)  Stab$_{\tilde{\Gamma}}(z)$ is a finite cyclic subgroup,
(b) each nontrivial torsion element in $\tilde{\Gamma}$ stabilizes a \emph{unique} point in $\mathfrak{H}$.
\end{exercise}

\begin{exercise}\label{Exercise_3.1.4} Show that $\tilde{\Gamma}$ acts \emph{properly discontinuously} on $\mathfrak{H}$ in the following sense: every $z \in \mathfrak{H}$ has an open neighborhood $N_z$ with the property that if $\gamma \in \tilde{\Gamma}$ then
$\gamma(N_z) \cap N_z = \phi$ if $\gamma \notin$ Stab$_{\tilde{\Gamma}}(z)$ and
$\gamma(N_z) \cap N_z = N_z$ otherwise. Conclude that the \emph{orbit space}
$\Gamma \setminus \mathfrak{H}$ is a topological $2$-manifold (a Hausdorff space
such that each point has an open neighborhood homeomorphic to $\mathbb{R}^2$).
\end{exercise}

\begin{exercise}\label{Exercise_3.1.5} Suppose that $f(\tau)$ is a nonzero weak modular form of weight $k$.
Show that $k$ is \emph{even}.
\end{exercise}

\begin{exercise}\label{Exercise_3.1.6} Let $E$ be the set of meromorphic modular functions of weight zero. Show that
$E$ is a field\footnote{Of course, field here is in the algebraic sense.} containing $\mathbb{C}$.
\end{exercise}

\begin{exercise}\label{Exercise_3.1.7} Show that pointwise multiplication defines a bilinear
product $\mathfrak{M}_k \otimes \mathfrak{M}_l \rightarrow \mathfrak{M}_{k+l}$, with respect to which
$\mathfrak{M} = \oplus_k \mathfrak{M}_k$ is a $\mathbb{Z}$-graded commutative $\mathbb{C}$-algebra.
\end{exercise}

\begin{exercise}\label{Exercise_3.1.8} Suppose that $f(\tau)$ is a meromorphic modular form of weight zero.
Show that $f'(\tau)$ is a meromorphic modular form of weight $2$.
\end{exercise}

\subsection{Eisenstein Series on $SL_2(\mathbb{Z})$}\label{Subsect_Eisenstein}
Beyond the fact that constants in $\mathbb{C}$ are modular functions
of weight $0$  (cf. Exercise \ref{Exercise_3.1.6}),  it is not so easy to construct nonconstant modular functions of weight $0$ or \emph{any} nonzero modular form of nonzero weight. We content ourselves 
with the description of some examples chosen because of their relevance to VOA theory.

\medskip
The most accessible nonconstant modular forms are the \emph{Eisenstein series}. For an integer $k\geq 2$,  set
\begin{eqnarray}
E_{k}(\tau) &=& -\frac{B_{k}}{k!}+\frac{2}{(k-1)!}\sum\limits_{n\geq 1}\frac{n^{k-1}q^{n}}{1-q^{n}} \notag \\
&=&  -\frac{B_{k}}{k!}+\frac{2}{(k-1)!}\sum\limits_{n\geq 1} \sigma_{k-1}(n)q^n. \label{Eisen1}
\end{eqnarray}
Here, $\sigma_{k-1}(n) = \sum_{d|n}d^{k-1}$ and $B_k$ is the $kth$ Bernoulli number defined by\footnote{Several different conventions are used to define Bernoulli number in the literature.}
\begin{equation}\label{bernoullidef}
\frac{z}{q_{z}-1}= \sum\limits_{k\geq 0} \frac{B_{k}}{k!}z^k = 1-\frac{1}{2}z+\frac{1}{12}z^2+ \hdots
\end{equation}
The following well-known identity of Euler 
\begin{eqnarray}\label{Eulerid}
 \zeta(k) = - \frac{(2\pi i)^{k}B_{k}}{2(k!)} \ (k \geq 2 \ \mbox{even}),
 \end{eqnarray}
permits us to reexpress the constant term of (\ref{Eisen1}) in terms of zeta-values.
The basic fact is this: \emph{Let $k\geq 3$. Then $E_{k}(\tau)$ is a holomorphic modular form of weight 
$k$; it is identically zero if, and only if, $k$ is odd}. We will see one way to prove this in Section \ref{Section_EllipticFuns}. We emphasize that $E_2(\tau)$ is \emph{not} a modular form.

\medskip
 The normalization employed in
(\ref{Eisen1}) is related to elliptic functions (Section \ref{Section_EllipticFuns}). In fact $B_{2k}$ never vanishes, so we can renormalize so that the $q$-expansion begins $1+ \hdots$. We single out the first three Eisenstein series corresponding to $k=2, 4, 6$ renormalized in this way, and rename them (following Ramanujan)
\begin{eqnarray*}
P &=& 1 - 24\sum_{n\geq 1} \sigma_1(n)q^n,  \\
Q&=& 1+240\sum_{n\geq 1} \sigma_3(n)q^n,  \\
R&=& 1 - 504 \sum_{n \geq 1} \sigma_5(n)q^n.
\end{eqnarray*}
$P, Q, R$ are \emph{algebraically independent}, so that they generate a weighted polynomial algebra, denoted 
\begin{eqnarray}\label{qmodformdef}
\mathfrak{Q} = \mathbb{C}[P, Q, R],
\end{eqnarray}
where $P,  Q, R$ naturally have weights (degree) $2, 4, 6$ respectively.
$\mathfrak{Q}$ is the algebra of \emph{quasimodular forms}. $\mathfrak{Q}$ contains every holomorphic modular form. Indeed, we have (cf. Exercise \ref{Exercise_3.1.7}) 
\begin{theorem}\label{thm4.1} The graded algebra $\mathfrak{M} = \oplus \mathfrak{M}_k$ of holomorphic modular forms on $\Gamma$ is the graded subalgebra 
$\mathbb{C}[Q, R]$ of $\mathfrak{Q}$. 
\end{theorem}
Theorem \ref{thm4.1} follows from a careful study of the singularities (zeros and poles) of modular forms, but we will not discuss this here. The Theorem contains a lot of information about holomorphic modular forms.
For example,  there are no such nonzero  forms of negative weight or weight $2$, holomorphic
forms of weight zero are necessarily constant, and $\dim \mathfrak{M}_k < \infty$. Indeed, inasmuch as $Q$ and $R$ are free generators in weights $4$ and $6$ respectively, the Hilbert-Poincar\'{e} series of $\mathfrak{M}$ is 
\begin{eqnarray}\label{HS1}
\sum_{k \geq 0} (\dim \mathfrak{M}_k) t^k &=& \frac{1}{(1-t^4)(1-t^6)} \\
&=& 1 + t^4+t^6+t^8+t^{10}+2t^{12}+ t^{14}+2t^{16}  \hdots \notag
\end{eqnarray}

\medskip
As we already mentioned, $E_2(\tau)$ is \emph{not} a modular form. Indeed, it satisfies the transformation law
\begin{equation}
E_{2}|_2 \gamma (\tau) =E_{2}(\tau )-\frac{c}{2\pi i(c\tau + d)}, 
\ \gamma = \left(\begin{array}{cc}a & b \\c & d\end{array}\right). \label{gammaE2}
\end{equation}
The importance of $E_2(\tau)$ for us stems from its relation to \emph{derivatives} of modular forms.
  Suppose that $f(\tau)$ is a meromorphic modular form of weight $k$. We define the \emph{modular derivative} of
  $f(\tau)$ by 
\begin{equation}
D_k f(\tau) = Df(\tau) = \left(\theta+kE_2(\tau) \right)f_{k}(\tau). \label{ModDer}
\end{equation}
where $\theta = qd/dq$. One can show without difficulty (cf. Exercise \ref{Exercise_3.2.4}) that $Df_k(\tau)$
is modular of weight $k+2$, and is holomorphic if $f_k(\tau)$ is.

\medskip

\bigskip

\begin{exercise}\label{Exercise_3.2.1} Prove that
\begin{equation*}
\frac{q_{z}}{(1-q_{z})^{2}}=\frac{1}{z^{2}}-\sum\limits_{k\geq 2}\frac{B_{k}}{k!}(k-1)z^{k-2}.
\end{equation*}
Deduce that $B_k = 0$ for odd $k \geq 3$.
\end{exercise}

\begin{exercise}\label{Exercise_3.2.2} Prove that $E_{8}=\frac{3}{7}E_{4}^{2}$ and $E_{10}=\frac{5}{11}E_{4}E_{6}$.
\end{exercise}

\begin{exercise}\label{Exercise_3.2.3} Show that (\ref{HS1}) is equivalent to the formula $\dim \mathfrak{M}_{2k} = [k/6]$ if $k \equiv 1 (\mbox{mod} \ 6)$
and $1+[k/6]$ otherwise.
\end{exercise}

\begin{exercise}\label{Exercise_3.2.4} Prove that $(D_kf)|_{k+2}(\tau) = D_k(f|_k \gamma)(\tau)$ for any meromorphic
function $f(\tau)$. Conclude that $D_k$ induces a linear map $\mathfrak{M}_k \rightarrow \mathfrak{M}_{k+2}$.
\end{exercise}

\begin{exercise}\label{Exercise_3.2.5} Prove that $DE_{4}=14E_{6}$ and $DE_{6}=\frac{60}{7}E_{4}^{2}$.
\end{exercise}

\begin{exercise}\label{Exercise_3.2.6} Let $D: \mathfrak{M} \rightarrow \mathfrak{M}$ be the linear map whose restriction to
$\mathfrak{M}_k$ is $D_k$. Prove that $D$ is a \emph{derivation} of the algebra $\mathfrak{M}$.
\end{exercise}

\subsection{Cusp-Forms and Modular Functions on $SL_2(\mathbb{Z})$}\label{Subsect_cusp} 
Thanks to Theorem \ref{thm4.1}, every holomorphic modular form of weight $k$ is equal to a unique
homogeneous polynomial in $Q$ and $R$. In this Subsection we describe some important examples 
of particular relevance to VOAs. We start with the \emph{discriminant}, defined via
\begin{eqnarray}\label{Delta1}
\Delta(\tau) = \frac{Q^3-R^2}{12^3} = q - 24q^2 + \hdots
\end{eqnarray}
$\Delta(\tau)$ is evidently a holomorphic modular form of weight $12$. It may alternatively be described
by a $q$-product which goes back to Kronecker, namely
\begin{eqnarray}\label{Delta2}
\Delta(\tau) = q \prod_{n\geq 1} (1-q^n)^{24}.
\end{eqnarray}
This formula finds its natural place in the theory of elliptic functions. From our present vantage point,
the fact that (\ref{Delta1}) and (\ref{Delta2}) coincide is miraculous. Beyond the product formula, the properties that make $\Delta(\tau)$ important for us are the following:
it \emph{does not vanish in $\mathfrak{H}$}, and (up to scalars) it is the unique nonzero holomorphic modular form of least weight that \emph{vanishes
at $\infty$}. The nonvanishing property has a natural explanation in the theory of elliptic functions. Concerning the second property, we introduce \emph{cusp-forms} defined by
  \begin{eqnarray*}
\mathfrak{S}_k &=& \{ f(\tau) \in \mathfrak{M}_k \ | \ f \ \mbox{vanishes at $\infty$}\}, \\
\mathfrak{S} &=& \oplus_k \mathfrak{S}_k.
\end{eqnarray*}
Our assertions then amount to the following:
$\mathfrak{S}_k = 0$ for $k<12$ and $\mathfrak{S}_{12}= \mathbb{C}\Delta(\tau)$. Using (\ref{Delta1}), it follows that
$\Delta(\tau)^{-1}$ is an almost holomorphic modular form of weight -$12$ with a pole of order $1$ at $\infty$. Applications of these facts are given in Exercise \ref{Exercise_3.3.1}.

\medskip
Closely related to $\Delta(\tau)$ is the \emph{Dedekind $\eta$-function}, whose $q$-expansion 
is the $24$th root of that for $\Delta(\tau)$:
\begin{eqnarray}\label{etadef}
\eta(\tau) = q^{1/24}\prod_{n \geq 1}(1-q^n).
\end{eqnarray}
$\eta(\tau)$ is not a modular form in the sense that we have defined it. Note that
\begin{eqnarray}\label{invetadef}
\eta(\tau)^{-1} &=& q^{-1/24} \sum_{n \geq 0} p(n)q^n \\
&=& q^{-1/24}(1 + q + 2q^2+3q^3+5q^4+ \hdots), \notag
\end{eqnarray}
an identity which goes back to Euler. ($p(n)$ is the \emph{unrestricted partition function}.)

\medskip
Our next example is the famous $j$-function, defined by
\begin{eqnarray}\label{jdef}
j(\tau) = \frac{Q^3}{\Delta(\tau)} = q^{-1} + 744 + 196884q + \hdots
\end{eqnarray}
As the quotient of two modular forms of weight $12$, $j(\tau)$ has weight zero, and because
$\Delta^{-1}$ is almost holomorphic then so too is $j(\tau)$. With the notation of Exercise
\ref{Exercise_3.1.6}, $j(\tau) \in E$. Let $\Gamma \setminus \mathfrak{H}$ be the orbit space for the action of $\Gamma$ on $\mathfrak{H}$ (cf. Exercise \ref{Exercise_3.1.4}). Because the weight is zero, we see from (\ref{slashaction})
that $j$ induces a map
\begin{eqnarray*}
j: \Gamma \setminus \mathfrak{H} \rightarrow \mathbb{C}
\end{eqnarray*}
which turns out to be a \emph{homeomorphism}.
   It can be shown that $E = \mathbb{C}(j)$  is exactly the field of rational
functions in $j$.

\bigskip

\begin{exercise}\label{Exercise_3.3.1} Considered as a subspace of the algebra $\mathfrak{M}$, show that $\mathfrak{S}$ is the principal ideal generated by $\Delta$. 
\end{exercise}

\begin{exercise}\label{Exercise_3.3.2} Show that $\mathfrak{M}_k = \mathfrak{S}_{2k} \oplus \mathbb{C}E_{2k}$ for $k \geq 2$.
\end{exercise}

\begin{exercise}\label{Exercise_3.3.3} Prove that $\theta \eta(\tau) = -\frac{1}{2} \eta(\tau)E_2(\tau)$. Conclude
that $D_{12}\Delta(\tau) = 0$. Give another proof of this by using Theorem \ref{thm4.1}. 
\end{exercise}

\begin{exercise}\label{Exercise_3.3.4} Regard $D$ as a derivation of $\mathfrak{M}$ as in Exercise \ref{Exercise_3.2.6}. Show that
the \emph{space of $D$-constants} (i.e., the subspace of $\mathfrak{M}$ \emph{annihilated} by $D$)
is the polynomial algebra $\mathbb{C}[\Delta]$. 
\end{exercise}

\begin{exercise}\label{Exercise_3.3.5}  Prove that the \emph{ring of almost holomorphic modular functions
of weight zero} on $\Gamma$ is the space $\mathbb{C}[j]$ of polynomials in $j(\tau)$.
\end{exercise}

\section{Characters of Vertex Operator Algebras}\label{Section_CharsVOAs}
Fix a VOA $(V, Y, \bf{1}, \omega)$ with $\mathbb{Z}$-graded Fock space $V = \oplus V_n$
and central charge $c$. In this Section we introduce the idea of the 
\emph{character} of $V$ as a sort of analog of the character of a group representation.
This is essentially the theory of \emph{$1$-point correlation functions} on $V$.

\subsection{Zero Modes}\label{Subsect_Zeromodes}
We start with a useful calculation. Suppose that $v \in V_k, w \in V_m, n \in \mathbb{Z}$. 
Remembering that $L_i = \omega_{i+1}$, we have
\begin{eqnarray*}
L_0v_nw &=& ([\omega_1, v_n] + v_nL_0)w \\
&=& \left(\sum_i \binom{1}{i} (L_{i-1}v)_{n+1-i} + v_nL_0\right)w \\
&=&(\left(L_{-1}v)_{n+1} +(L_0v)_n + v_nL_0\right)w \\
&=&(m+k-n-1)v_nw.
\end{eqnarray*}
Here, we used the commutator formula (cf. Section \ref{Subsect_VAaxioms}) for the second equality
and the last identity of Exercise \ref{Exercise_2.6.2} for the fourth equality. What we take from this is that
\emph{modes of homogeneous states are graded operators on $V$}:
\begin{eqnarray}\label{modeaction1}
v \in V_k \Rightarrow v_n: V_m \rightarrow V_{m+k-n-1}.\label{gradedmodes}
\end{eqnarray}
 In particular,
let us define\footnote{The zero mode of $v$
is generally not the \emph{zeroth} mode $v_0$ but rather the mode which has weight zero as an operator. However, in the convention used for modes in CFT as practiced by physicists, it is the zero mode.}  the \emph{zero mode} $o(v)$ of a state $v \in V_k$ to be $v_{k-1}$, and extend this definition to $V$ additively. From (\ref{modeaction1}) we then have for all integers $m$ and states $v$ that
\begin{eqnarray}
o(v): V_m \rightarrow V_m.\label{o(v)}
\end{eqnarray}

The point of all this is that as an operator on $V_m$ we can \emph{trace} the zero
mode and form a generating function (cf. (\ref{qopformalism}))
\begin{eqnarray}\label{Zdef}
Z(v, q) = \mathrm{Tr}_V o(v)q^{L_0-c/24} = q^{-c/24}\sum_n \mathrm{Tr}_{V_n}o(v)q^n.
\end{eqnarray}
(\ref{Zdef}) is to be regarded as a 
\emph{formal}  $q$-expansion at this point. Apart from memorializing the central charge,
the factor $q^{-c/24}$ may seem somewhat arbitrary. This feeling will pass.
Because the homogeneous spaces $V_n$ vanish for small enough $n$, we see that
\begin{eqnarray*}
Z(v, q) \in q^{-c/24}\mathbb{C}[[q]][q^{-1}].
\end{eqnarray*}
$Z=Z_V$ defines the \emph{character} of $V$, i.e., the linear map
\begin{eqnarray*}
&&Z: V \rightarrow q^{-c/24}\mathbb{C}[[q]][q^{-1}] \\
&&     \ \ \ \ \     v \mapsto Z(v, q).
\end{eqnarray*}

\bigskip

\begin{exercise}\label{Exercise_4.1.1} Let $U\otimes V$ be the tensor product of VOAs $U, V$ (Exercise 
\ref{Exercise_2.6.7}).  Prove that $Z_{U\otimes V} = Z_U Z_V$.
\end{exercise}

\begin{exercise}\label{Exercise_4.1.2} Suppose that $V$ is a VOA with $v \in V$. Prove the identity  $q_x^{L_0}Y(v, z)q_x^{-L_0} = Y(q_x^{L_0}v, q_xz)$.
\end{exercise}

\begin{exercise}\label{Exercise_4.1.3} Let $a$ be the generating state of weight $1$ for the Heisenberg VOA
(Subsection \ref{Subsect_HeisVOA}). Prove that the zero mode $o(a)$ is \emph{zero}.
\end{exercise}

\begin{exercise} \label{Exercise_truncweight}
Suppose that $V_0 = \mathbb{C}\bf{1}$
(cf. Exercise \ref{Exercise_2.6.4}). Prove using (\ref{gradedmodes}) that for $a\in V_{n_1}$ and $b\in V_{n_2}$ we have $a_n b=0$ for  all $n\ge n_1+n_2$. Using Exercise \ref{Exercise_local_Trunc},
deduce that $Y(a, z) \sim_k Y(b, z)$ with order of locality $k\le  n_1+n_2$. 
\end{exercise}

\subsection{Graded Dimension}\label{Subsect_gradeddim}
The most prominent $Z$-value is that obtained by tracing the zero mode of the vacuum.
From Exercise \ref{Exercise_2.6.1} we have $Y(\mathbf{1}, z) = \mathrm{Id}_V$ and $\mathbf{1} \in V_0$.
So the zero mode of $\mathbf{1}$
is $\mathrm{Id}_V$, whence
\begin{eqnarray}\label{0point}
Z_V(\mathbf{1}) = \mathrm{Tr}_Vq^{L_0-c/24} = q^{-c/24}\sum_n \dim V_n q^n.
\end{eqnarray}
This is variously called the \emph{graded dimension}, \emph{$q$-dimension}, \emph{$0$-point function}, or \emph{partition function} of $V$. 

\medskip
The graded dimensions of our two main examples $M_0$ and $\mbox{Vir}_c$ are readily computed.
This is because the Fock spaces are Verma modules, or closely related to them, and these are
easy to handle. Let us start with the Fock space $M_0$ for the free boson, which has central charge $c=1$ (Theorem \ref{thm2.7} and Exercise \ref{Exercise_2.6.7}).
In the notation of (\ref{MVerma}), $M_0$ (considered as a $\mathbb{Z}$-graded linear space) coincides with $\mathcal{U}(\hat{A}^-)$ equipped with the natural product grading for which
$a\otimes t^{-n}$ has weight $n$. Because of the PBW Theorem,
the universal enveloping algebra is itself isomorphic as graded space to the symmetric algebra
$S(\coprod_{n\geq 1} \mathbb{C}x_{-n})$ with $x_{-n}$ having weight $n$. (In other words, $M_0$ 
`is' a polynomial algebra in variables $x_{-n}$. Compare with Remark \ref{Remark_Mboson}.) As graded algebras, symmetric algebras are multiplicative over direct sums. It follows that
\begin{eqnarray*}
Z_{M_0}(\mathbf{1}) &=&q^{-1/24} \prod_{n=1}^{\infty} (q\mbox{-dimension of} \ \mathbb{C}[x_{-n}]) \\
&=&q^{-1/24} \prod_{n=1}^{\infty} (1+q^n+q^{2n} + \hdots) \\
&=& q^{-1/24} \prod_{n=1}^{\infty} (1- q^n)^{-1},
\end{eqnarray*}
which is none other than the inverse eta-function (\ref{etadef}), (\ref{invetadef}).
Thus we have 
\begin{eqnarray}\label{M0gdim}
Z_{M_0}(\mathbf{1}) = \eta(q)^{-1}.
\end{eqnarray}
For an integer $n \geq 1$, let $M_0^{\otimes n}$ be the $n$-fold tensor product of $M_0$ considered
as a VOA as described in Exercise \ref{Exercise_2.6.7}. This is the theory of \emph{$n$ free bosons}. Using Exercise \ref{Exercise_3.1.1} we deduce from (\ref{M0gdim}) that
\begin{eqnarray}\label{M0dgdim}
Z_{M_0^{\otimes n}}(\mathbf{1}) = \eta(q)^{-n}.
\end{eqnarray}
In particular, the graded dimension of the VOA $M_0^{\otimes 24}$ of $24$ free bosons
(the \emph{bosonic string})  is
the inverse discriminant $\Delta(\tau)^{-1}$.

\medskip
The calculation of the graded dimension of $\mbox{Vir}_c$  is similar. Indeed, the Fock space $M_{c, 0}$ (\ref{VirVerm}) is \emph{isomorphic} as $\mathbb{Z}$-graded linear space to $M_0$. We must quotient out  the graded submodule 
$\mathcal{U}(\mbox{Vir})L_{-1}. v_{c, 0}$, and this
is isomorphic to $M_{c, 0}[1]$, that is $M_{c, 0}$ with an overall shift of $+1$ in the grading, because
$L_{-1}. v_{c, 0}$ has weight $1$ as an element of $M_{c, 0}$. We find that
\begin{eqnarray}\label{vircqdim}
Z_{\mbox{Vir}_c}(\mathbf{1}) = \frac{q^{-c/24}}{\prod_{n\geq 2}(1-q^n)}, 
\end{eqnarray}
which is \emph{not} the $q$-expansion of a modular form.

\bigskip
Next we consider the character value $Z_V(\omega)$ for a VOA $V$. Because the zero mode
of the conformal vector is $L_0$, which acts on $V_n$ as multiplication by $n$, we have
\begin{eqnarray*}
Z_V(\omega) = q^{-c/24}\sum_n n\dim V_nq^n.
\end{eqnarray*}
 This is almost, but not quite, equal to $\theta Z_V(\mathbf{1})$ ($\theta$ as in (\ref{ModDer})).
If instead we use the square bracket conformal vector $\tilde{\omega} \in V_{[2]}$ (\ref{sqomega}) we find
\begin{eqnarray*}
Z_V(\tilde{\omega}) &=& q^{-c/24}\sum_n (n- \frac{c}{24})\dim V_n q^n \\
&=& \theta(Z_V(\mathbf{1})).
\end{eqnarray*}
In the case of $M_0$, for example, we obtain using Exercise \ref{Exercise_3.3.3} that
\begin{eqnarray*}
Z_{M_0}(\tilde{\omega})= \theta\eta(\tau)^{-1} = \frac{E_2(\tau)}{2\eta(\tau)}.
\end{eqnarray*}
This suggests that  `nicer'  character values obtain by evaluating $Z_V$
on states which are homogeneous in the square bracket formalism, i.e., lie in
$V_{[k]}$ for some $k$.

\subsection{The Character of the Heisenberg Algebra}\label{Subsect_CharHeisen}
It is generally a difficult problem to compute the $1$-point functions $Z_V(v)$ of a VOA $V$
for a complete basis of states. We describe the solution for the Heisenberg algebra $M_0$ (\cite{MT1}).
It well illustrates the principle suggested at the end of the previous Subsection.
\begin{theorem}\label{thm4.15} Let $M_{[0]} = \oplus_{n\geq 0} (M_0)_{[n]}$ be the Fock space for
$M_0$ equipped with the square bracket grading (cf. Section \ref{Subsect_SquareVOA}). Let $\mathfrak{Q}$ be the graded algebra of quasimodular forms (\ref{qmodformdef}). There is a \emph{surjection} of graded linear spaces
\begin{eqnarray*}
M_{[0]} \rightarrow \mathfrak{Q}, \ \ v \mapsto Q_v(\tau)
\end{eqnarray*}
such that $Z_{M_0}(v) = Q_v(\tau)/\eta(\tau)$.
\end{theorem}
Up to a normalizing factor $\eta(\tau)^{-1}$ then, every $1$-point function is a quasimodular form, and every quasimodular form of weight $k$ arises in this way from a state $v \in (M_0)_{[k]}$ (cf. Exercise \ref{Exercise_4.3.2}). There is an \emph{explicit} description of the quasimodular form
$Q_v(\tau)$ attached to a state $v$ with $wt[v]=k$ which goes as follows.
A basis of states for $(M_0)_{[k]}$ is given by
\begin{eqnarray}\label{[F]basis}
v_{\lambda} =a[-k_{1}]\ldots a[-k_{n}]\mathbf{1}
\end{eqnarray}
where $k = k_1+\hdots +k_n$ and $1\leq k_{1}\leq \ldots \leq k_{n}$ range over the parts of a partition $\lambda$ of $k$. The quasimodular form
$Q_{v_{\lambda}}(\tau)$ is given by
\begin{eqnarray}\label{Qvdef}
Q_{v_{\lambda}}(\tau) = \sum_{\varphi= ...(rs)...}\  \prod_{(rs)} (-1)^{r+1}\frac{(r+s-1)!}{(r-1)!(s-1)!}E_{r+s}(\tau),
\end{eqnarray}
where the notation is as follows. Let $\Phi = \{k_1, \hdots, k_n \}$ be the parts of the partition $\lambda$.
Then $\varphi$ ranges over all fixed-point-free involutions in the symmetric group
$\Sigma(\Phi)$, so that $\varphi$ can be represented as a product of transpositions
$...(rs)...$ with $(r, s)$ a pair of parts of $\lambda$. For each such $\varphi$, the product ranges over
the transpositions whose product (in $\Sigma(\Phi)$) is $\varphi$. We will indicate how 
(\ref{Qvdef}) can be proved in the next Section. A detailed proof appears later in Subsection \ref{Subsect_Heisen1pt}.   

\bigskip
\noindent
Assume formula (\ref{Qvdef}) in the following exercises.

\begin{exercise}\label{Exercise_4.3.1} Show that $Q_{v_{\lambda}}(\tau)$ vanishes if $\lambda$ has either an odd number of parts or an odd number of odd parts.
\end{exercise}

\begin{exercise}\label{Exercise_4.3.2}  Assume that $\lambda$ has both an even number of parts and an even number of odd parts. Prove that $Q_{v_{\lambda}}(\tau)$ has a \emph{nonzero} constant term, and in particular does not vanish.
\end{exercise}

\begin{exercise}\label{Exercise_4.3.3} With the same assumptions as the previous Exercise, prove that $Q_{v_{\lambda}}(\tau) \in \mathfrak{M}$ if, and only if, $\lambda$ has \emph{at most one part equal to $1$}.
\end{exercise}

\begin{exercise}\label{Exercise_4.3.4} Prove the assertion that \emph{every} quasimodular form arises as the trace of a state in $M_0$.
\end{exercise}

\section{Elliptic Functions and $2$-Point Functions}\label{Section_EllipticFuns}
There is an extension of the idea of  $1$-point functions to $n$-point functions for any
(nonnegative) $n$. We mainly restrict ourselves here to the case of $2$-point correlation 
functions, which are related to \emph{elliptic functions}.

\subsection{Elliptic Functions}\label{Subsect_Elliptic}
Throughout this Section,  \emph{lattice} means an additive subgroup 
$\Lambda \subseteq \mathbb{C}$ of rank $2$. As such it is the $\mathbb{Z}$-span
of an $\mathbb{R}$-basis $(\omega_1, \omega_2)$ of $\mathbb{C}$.
An \emph{elliptic function} is a function $f(z)$ which is meromorphic in
$\mathbb{C}$ and satisfies $f(z+\lambda) = f(z)$ for all $\lambda$
in some lattice $\Lambda$. Equivalently, $f(z+\omega_i) = f(z)$ for 
basis vectors $\omega_1, \omega_2$ of $\Lambda$. $\Lambda$ is the 
\emph{period lattice} of $f(z)$. Note that $\mathbb{C}/\Lambda$ has the structure of a
\emph{complex torus} (aka \emph{complex elliptic curve}) and that $f(z)$ induces a map
\begin{eqnarray*}
f: \mathbb{C}/\Lambda \rightarrow \mathbb{CP}^1= \mathbb{C}\cup \{\infty\}.
\end{eqnarray*}
The set of all meromorphic functions with period lattice $\Lambda$ is a
field $M_{\Lambda}$ (the function field of the torus). We have $\mathbb{C} \subseteq M_{\Lambda}$ 
where $\mathbb{C}$ is identified with the constants, moreover $f'(z) \in M_{\Lambda}$ whenever $f(z) \in M_{\Lambda}$.

\medskip
Two lattices $\Lambda_1, \Lambda_2$ are \emph{homothetic} if there is $\alpha \in \mathbb{C}$
with $\alpha\Lambda_1 = \Lambda_2$. It is usually enough to deal with some fixed lattice in a homothety class (the corresponding complex tori are  \emph{isomorphic}), and every $\Lambda$
is homothetic to a lattice with basis $(2\pi i, 2\pi i \tau)$ 
and  $\tau \in \mathfrak{H}$. We let $\Lambda_{\tau}$ denote this lattice.

\medskip
The classical \emph{Weierstrass $\wp $-function} is\footnote{Here and below,
a prime appended to a summation indicates that terms rendering the sum meaningless,
in this case $(m,n)=(0,0)$, are to be omitted.}

\begin{equation}
\wp (z,\tau )=\frac{1}{z^{2}}+\sum_{m,n\in \mathbb{Z}}^{\prime } \left\{\frac{1}{%
(z-\omega _{m,n})^{2}}-\frac{1}{\omega _{m,n}^{2}} \right\}.  \label{Weierstrass}
\end{equation}%
Here, $\omega _{m,n}=2\pi i(m\tau +n)$. The double sum is  independent of the order of
summation and absolutely convergent. It defines a function with the following properties: (a) double pole at each point of 
$\Lambda_{\tau}$, (b) holomorphic in $(\mathbb{C}\times \mathfrak{H}) \setminus \Lambda_{\tau}$,
(c) \emph{even} in $z$, (d) period lattice $\Lambda_{\tau}$.
In particular,  for fixed $\tau \in \mathfrak{H}$ the $\wp$-function $\wp(z, \tau)$ lies in the field
$M_{\Lambda_{\tau}}$. It turns out that 
\begin{eqnarray*}
M_{\Lambda_{\tau}} = \mathbb{C}(\wp(z, \tau), \wp'(z, \tau))
\end{eqnarray*}
is a \emph{function field in one variable}. Indeed the set of \emph{even} elliptic functions is a simple transcendental extension
$\mathbb{C}(\wp)$ and $M_{\Lambda_{\tau}} \supseteq \mathbb{C}(\wp)$ a \emph{quadratic extension}.

\bigskip
There is a natural left action of $\Gamma$ on $\mathbb{C}\times \mathfrak{H}$ extending
(\ref{GHaction}). It is given by
\begin{equation}
\gamma : (z, \tau) \mapsto \left(\frac{z}{c\tau +d},\frac{a\tau +b}{c\tau +d}\right), \ 
\gamma =\left( 
\begin{array}{cc}
a & b \\ 
c & d%
\end{array} \right) \in \Gamma,  \label{Modu}
\end{equation}
corresponding to a base change $2\pi i(\tau
,1)\mapsto 2\pi i(a\tau +b,c\tau +d)$ of $\Lambda_{\tau}$ followed by the homothety
(conformal rescaling) $z \mapsto z/(c\tau +d)$. Then it follows that%
\begin{equation}
\wp (\gamma (z, \tau))=(c\tau +d)^{2}\wp (z,\tau ).  \label{Weiermod}
\end{equation}
This says that $\wp(z, \tau)$ is \emph{Jacobi form} of weight $2$ (\cite{EZ}), though we will neither 
explain nor pursue this idea here. 

\medskip
What we need is that $\wp (z,\tau )$ is invariant under
$\tau \mapsto \tau +1$ as well as $z \mapsto z+2 \pi i$ (from the elliptic property).
 It follows
that $\wp (z,\tau )$ has a Fourier expansion in both $q$ and $q_{z}$. (Compare 
with the development in Subsection \ref{Subsect_SL2ZModForms}.)
To
describe this we define 
\begin{eqnarray}
P_{1}(z,\tau ) &=&\sum\limits_{n\in \mathbb{Z}}^{\prime }\frac{q_{z}^{n}}{%
1-q^{n}}-\frac{1}{2},  \label{P1def} \\
P_{2}(z,\tau ) &=&\frac{d}{dz}P_{1}(z,\tau
)=\sum\limits_{n\in \mathbb{Z}}^{\prime }\frac{nq_{z}^{n}}{1-q^{n}}.
\label{P2def}
\end{eqnarray}%
The extra term $-1/2$ in (\ref {P1def}) ensures that $P_{1}(z,\tau )$ is odd in $z$.
For nonzero $z$ in the fundamental parallelogram defined by the basis $(2\pi i, 2\pi i \tau)$
of $\Lambda_{\tau}$ we have $-2\pi \Im \tau <%
\Re(z)$ $<0$, so that $\left\vert q\right\vert <\left\vert
q_{z}\right\vert <1$. 
$P_1(z, \tau)$ and its $z$-derivatives are absolutely convergent in this domain. We can now give the Fourier expansion of the $\wp$-function, which reveals a fundamental relationship with the Eisenstein series of Subsection \ref{Subsect_Eisenstein}.
\begin{theorem}\label{thm5.1} We have
\begin{eqnarray*}
\wp(z, \tau) &=& P_2(z, \tau) - E_2(\tau) \\
&=&  \frac{1}{z^2} + \sum_{k \geq 2} (2k-1)E_{2k}(\tau)z^{2k-2}. 
\end{eqnarray*}
\end{theorem}
We sketch the proof, which uses a key identity (cf. Exercise \ref{Exercise_5.1.1} below):
\begin{equation}
\sum_{n\in \mathbb{Z}}\frac{1}{(x-2\pi in)^{2}}=\frac{q_{x}}{(1-q_{x})^{2}}
\label{sinhsum}
\end{equation}
for $x\not = 0$. In the exceptional case,
\begin{equation}
\sum'_{n\in \mathbb{Z}}\frac{1}{(2\pi in)^{2}}= \frac{2\zeta(2)}{(2\pi i)^2} = -\frac{1}{12}.
\label{zetasum}
\end{equation}
Now
\begin{equation}\label{wp1}
\wp (z,\tau )+\sum_{m \in \mathbb{Z}}\left[\sum_{n \in \mathbb{Z}} \frac{1}{\omega_{m, n}^2} \right]=\sum_{m\in \mathbb{Z}}\left[ \sum_{n\in \mathbb{Z}%
}\frac{1}{\left( z-\omega _{m,n}\right) ^{2}}\right] ,
\end{equation}
where the convergent nested double sums depend on the order of summation.
For the lhs, use (\ref{sinhsum}) with $x=2\pi i m\tau \not= 0$, (\ref{zetasum})
and $|q|<1$ to obtain
\begin{eqnarray*}
\sum_{m \in \mathbb{Z}}\left[\sum_{n \in \mathbb{Z}} \frac{1}{\omega_{m, n}^2} \right]
&=& -\frac{1}{12} +\sum_{0\not = m \in \mathbb{Z}}  \frac{q^m}{(1-q^m)^2}  \\
&=&-\frac{1}{12} +2\sum_{m, n >0}nq^{mn} = E_2(\tau)
\end{eqnarray*}
(cf. (\ref{Eisen1})). For the rhs of (\ref{wp1}), use (\ref{sinhsum}) with $x=z-2\pi im\tau$ and argue similarly using $|q_{z}q^{m}|,|q_{z}^{-1}q^{m}|<1$ for $m>0$ to get
\begin{eqnarray}
\sum_{m\in \mathbb{Z}}\frac{q_{z}q^{m}}{(1-q_{z}q^{m})^{2}} 
&=&\frac{q_{z}}{(1-q_{z})^{2}}+\sum_{m>0}\left( \frac{q_{z}q^{m}}{%
(1-q_{z}q^{m})^{2}}+\frac{q_{z}q^{-m}}{(1-q_{z}q^{-m})^{2}}\right)  \notag \\
&&\frac{q_{z}}{(1-q_{z})^{2}}+\sum_{m>0}%
\sum_{n>0}n(q_{z}^{n}+q_{z}^{-n})q^{nm}  \notag \\
&=&\frac{q_{z}}{(1-q_{z})^{2}}+\sum_{n>0}n\left( q_{z}^{n}+q_{z}^{-n}\right) 
\frac{q^{n}}{1-q^{n}} \\ \label{qzexp}
&=&\sum\limits_{n>0}n\left( \frac{q_{z}^{n}}{1-q^{n}}-\frac{q_{z}^{-n}}{%
1-q^{-n}}\right) =P_{2}(z,\tau ). \notag
\end{eqnarray}
This proves the first equality in the Theorem. From (\ref{Weierstrass}) we see that 
\begin{equation*}
\wp(z,\tau )=\frac{1}{z^{2}}+\sum\limits_{k\geq 3}(k-1)\tilde{E}_{k}(\tau )z^{k-2},
\end{equation*}%
with
\begin{equation}
\tilde{E}_{k}(\tau )=\sum_{m,n\in \mathbb{Z}}^{\prime }\frac{1}{\omega _{m,n}^{k}}=%
\frac{1}{(2\pi i)^{k}}\sum_{m,n\in \mathbb{Z}}^{\prime }\frac{1}{(m\tau
+n)^{k}}.  \label{EisenSum}
\end{equation}%
We can use
(\ref{sinhsum}) to identify $\tilde{E}_k$ with the corresponding Eisenstein series 
$E_k(\tau)$ (\ref{Eisen1}), in particular $\tilde{E}_{k}(\tau)$ is identically zero for $k$ odd.
This completes our discussion of Theorem \ref{thm5.1}.

\medskip
We note that $P_{1}$ is not an elliptic function (cf. Exercise \ref{Exercise_5.1.5}).  
Higher $z$-derivatives $P_{1}^{(m)}(z,\tau )$ for $m\geq 1$ are
elliptic functions, and are derivatives of $\wp(z, \tau)$ for $m\geq 2$.  We have
\begin{eqnarray}
P_{1}^{(m)}(z,\tau ) &=&\sum\limits_{n\in \mathbb{Z}}^{\prime }\frac{%
n^{m}q_{z}^{n}}{1-q^{n}}  \notag \\
&=&m!\left[ \frac{(-1)^{m+1}}{z^{m+1}}+\sum\limits_{k\geq m+1}\binom{k-1}{m}%
E_{k}(\tau )z^{k-m-1}\right].  \label{P1m}
\end{eqnarray}

\bigskip

\begin{exercise}\label{Exercise_5.1.1} Prove directly from the definition that an elliptic function
which is  \emph{holomorphic}  is necessarily constant. 
\end{exercise}

\begin{exercise}\label{Exercise_5.1.2} Verify (\ref{sinhsum}) by comparing poles.
\end{exercise}

\begin{exercise}\label{Exercise_5.1.3} Prove that for even $k\geq 4$, $\tilde{E}_k(\tau)$ coincides with $E_k(\tau)$.
(Use (\ref{Eulerid}).)
\end{exercise}

\begin{exercise}\label{Exercise_5.1.4} Deduce from (\ref{Weiermod}) that $E_{k}(\tau) \in \mathfrak{M}_k$ for even $k\geq 4$.
\end{exercise}

\begin{exercise}\label{Exercise_5.1.5} Prove that $P_{1}(z+2\pi i\tau ,\tau )=P_{1}(z,\tau )-1$.
\end{exercise}

\subsection{$2$-Point Correlation Functions}\label{Subsect_2point}
Let $(V, Y, \mathbf{1}, \omega)$ be a VOA of central charge $c$.
For an integer $n\geq 0$, the $n$-point correlation function for states $u^1, \hdots, u^n \in V$
is the formal expression
\begin{eqnarray}
&&\ \ \ \ \ \ \ F_{V}((u^1,z_{1}),  \hdots, (u^n,z_{n}), q)= \notag \\
&&\mathrm{Tr}%
_{V}Y(q_{1}^{L_0}u^1,q_{1})\hdots Y(q_{n}^{L_0}u^n,q_{n})q^{L_0-c/24},
\label{Fnpt}
\end{eqnarray}%
where $q_i = q_{z_i}$ for variables $z_1, \hdots, z_n$. For $n=0$ this reduces to
the graded dimension Tr$_Vq^{L_0-c/24}$ as discussed in Subsection \ref{Subsect_gradeddim}.
If $n=1$ and $u^1 \in V_k$, (\ref{Fnpt}) is 
\begin{eqnarray*}
\mathrm{Tr}_VY(q_1^{L_0}u^1, q_1)q^{L_0-c/24}
&=& q_1^k \sum_m\mathrm{Tr}_V u^1_m q_1^{-m-1}q^{L_0-c/24} \\
&=& \mathrm{Tr}_Vo(u^1)q^{L_0-c/24} = Z(u^1, q),
\end{eqnarray*}
where we used (\ref{modeaction1}) to get the second equality.  So for $n=1$, (\ref{Fnpt})
is the $1$-point function of Section \ref{Section_CharsVOAs}, as expected. There are similar modal expressions 
 for all $n$-point functions, but for $n\geq 2$ they are unhelpful. Here we focus on the 
 $2$-point function
\begin{equation}
F_{V}((u,z_{1}),(v,z_{2}),\tau )=\mathrm{Tr}%
_{V}Y(q_{1}^{L_0}u,q_{1})Y(q_{2}^{L_0}v,q_{2})q^{L_0-c/24}.
\label{Ftwopt}
\end{equation}%

We want to re-express the $2$-point function as a \emph{1-point function}, and for this we need be  able to manipulate vertex operators. More precisely, we need to manipulate expressions involving
vertex operators which are \emph{traced} over $V$. In such a context the locality of operators
(\ref{local1def}) 
simplifies in the sense that 
\begin{eqnarray*}
\mathrm{Tr}_V Y(u, z_1)Y(v, z_2)q^{L_0} = \mathrm{Tr}_V Y(v, z_2)Y(u, z_1)q^{L_0},
\end{eqnarray*}
where the additional factor $(z_1-z_2)^k$ (loc. cit.) has conveniently disappeared. Similar comments apply to
the associativity formula (\ref{associativity}), where we have
\begin{eqnarray*}
\mathrm{Tr}_V Y(u, z_1+z_2)Y(v, z_2)q^{L_0} = \mathrm{Tr}_V Y(Y(u, z_1)v, z_2)q^{L_0}.
\end{eqnarray*}
These and similar assertions fall under the heading of \emph{duality} in CFT, which is discussed in
\cite{FHL}.  We shall use them below without further comment.
Thus with some changes of variables together with Exercise \ref{Exercise_4.1.2},  we have
\begin{eqnarray}
&&F_{V}((u,z_{1}),(v,z_{2}),\tau ) \notag \\
&=& \mathrm{Tr}_{V} Y(Y(q_{1}^{L_0}u,q_{1}-q_{2})q_{2}^{L_0}v,q_{2})q^{L_0-c/24}  \notag
\\
&=&\mathrm{Tr}_{V}
Y(q_{2}^{L_0}
Y(q_{z_{12}}^{L_0}u, q_{z_{12}}-1)v,q_{2})
q^{L_0-c/24} \notag \\
&=&Z_{V}(Y[u,z_{12}]v,\tau ),  \label{FVz12}
\end{eqnarray}%
where $z_{12}=z_{1}-z_{2}$. This is the desired $1$-point function.  Similarly,  
\begin{eqnarray*}
&&F_{V}((u,z_{1}),(v,z_{2}+2\pi i\tau ),\tau ) \\
&=&q^{-c/24}\mathrm{Tr}_{V} Y(q_{1}^{L_0}u,q_{1})Y(q^{L_0}q_{2}^{L_0}v,qq_{2})q^{L_0} \\
&=&q^{-c/24}\mathrm{Tr}_{V} Y(q_{1}^{L_0}u,q_{1})q^{L_0}Y(q_{2}^{L_0}v,q_{2}) \\
&=&q^{-c/24}\mathrm{Tr}_{V} Y(q_{2}^{L_0}v,q_{2})Y(q_{1}^{L_0}u,q_{1})q^{L_0} \\
&=&F_{V}((u,z_{1}),(v,z_{2}),\tau ),
\end{eqnarray*}%
Thus $F_{V}$ is periodic in $z_2$
with period $2\pi i\tau $, and the same holds for $z_1$. It is obvious that
$F_V$ is also periodic in each $z_i$ with period $2\pi i$. It follows that, at least formally,
the $2$-point function $F_V$ (alias the $1$-point function (\ref{FVz12}))
is elliptic in the variable $z_{12}$ with period lattice $\Lambda_{\tau}$.

\subsection{Zhu Recursion Formula I}\label{Subsect_ZhuI}
We continue to pursue the ellipticity of the $2$-point function $F_V$. 
 It is the analytic of $F_V$ which needs to be established. To this end we develop a recursion formula of Zhu (\cite{Z}) which finds a number of applications.

\begin{theorem}
\label{TheoremZhu} We have
\begin{eqnarray}
&& \hspace{3.5cm} F_{V}((u,z_{1}),(v,z_{2}),\tau ) \notag \\
&&=\mathrm{Tr}_{V} o(u)o(v)q^{L_0-c/24} - \sum\limits_{m\geq 1}\frac{(-1)^{m}}{m!}P_{1}^{(m)}(z_{12},\tau
)Z_{V}(u[m]v,\tau ).  \label{ZhuRecuruv}
\end{eqnarray}
\end{theorem}

 The sum in (\ref{ZhuRecuruv}) is finite since $u[m]v=0$ for $m$ sufficiently
large, and from Subsection \ref{Subsect_Elliptic} $P_{1}^{(m)}(z_{12},\tau )$ is
elliptic for $m\geq 1$. Thus the ellipticity of $F_V$ is reduced to the \emph{convergence}
of Tr$_Vo(u)o(v)$ and the $1$-point functions $Z_{V}(u[m]v,\tau )$. This Theorem makes clear the deep connection between elliptic functions (and therefore also modular forms) and VOAs. There is an analogous recursion for all $n$-point functions.

\medskip
To prove Theorem \ref{TheoremZhu} we may assume that $u \in V_k$, whence
\begin{eqnarray}
&& \ \ \ \ \ \ \ \ \ \ \ F_{V}((u,z_{1}),(v,z_{2}),\tau )= \notag \\
&&\sum_{n\in \mathbb{Z}}q_{1}^{-n-1+k}
\mathrm{Tr}_{V}\left( u_nY(q_{2}^{L_0}v,q_{2})q^{L_0-c/24}\right) .
\label{FV2}
\end{eqnarray}
Using (\ref{modeaction1}), Exercise \ref{Exercise_4.1.2} and (\ref{square4}) we have
\begin{eqnarray*}
\left[ u_n,Y(q_{2}^{L_0}v,q_{2})\right]  &=&\sum\limits_{i\geq 0}\binom{n}{i}Y(u_iq_{2}^{L_0}v,q_{2})q_{2}^{n-i} \\
&=&q_{2}^{r}Y(q_{2}^{L_0}\sum\limits_{i\geq 0}\binom{n}{i}%
u_iv,q_{2}) \\
&=&q_{2}^{r}\sum\limits_{m\geq 0}\frac{r^{m}}{m!}Y(q_{2}^{L_0}u[m]v,q_{2}),
\end{eqnarray*}%
where $r=n+1-k$.
Hence 
\begin{eqnarray*}
&&\hspace{3.5cm} \mathrm{Tr}_{V}\left( u_nY(q_{2}^{L_0}v,q_{2})q^{L_0-c/24}\right)  \\
&=&\mathrm{Tr}_{V}\left( [u_n,Y(q_{2}^{L_0}v,q_{2})]q^{L_0-c/24}\right) +%
\mathrm{Tr}_{V}\left( Y(q_{2}^{L_0}v,q_{2})u_nq^{L_0-c/24}\right)  \\
&=&q_{2}^{r}\sum\limits_{m\geq 0}\frac{r^{m}}{m!}Z_{V}(u[m]v,\tau )+q^r
\mathrm{Tr}_{V}\left( Y(q_{2}^{L_0}v,q_{2})q^{L_0-c/24}u_n\right).
\end{eqnarray*}%
From this we obtain
\begin{equation*}
q_{2}^{r}\sum\limits_{m\geq 0}\frac{r^{m}}{m!}Z_{V}(u[m]v,\tau )  = 
(1-q^r)\mathrm{Tr}_V \left(u_nY(q_{2}^{L_0}v,q_{2})q^{L_0-c/24} \right), 
\end{equation*}%
so that for $r \neq 0$ we have 
\begin{equation*}
\mathrm{Tr}_{V}\left(u_nY(q_{2}^{L_0}v,q_{2})q^{L_0-c/24}\right) =\frac{%
q_{2}^{r}}{1-q^{r}}\sum\limits_{m\geq 1}\frac{r^{m}}{m!}Z_{V}(u[m]v,\tau ).
\end{equation*}%
Finally, the term corresponding to $r=0$ in (\ref{FV2}) is
$\mathrm{Tr}_V o(u)o(v)q^{L_0-c/24}$. Substituting into (\ref{FV2}),  we find 
\begin{eqnarray*}
F_{V}((u,z_{1}),(v,z_{2}),\tau ) &=&\mathrm{Tr}_{V}\left(
o(u)o(v)q^{L_0-c/24}\right)  \\
&&+\sum\limits_{m\geq 1}\frac{1}{m!}Z_{V}(u[m]v,\tau )\sum\limits_{n\in 
\mathbb{Z}}^{\prime }\frac{r^{m}q_{z_{12}}^{-n}}{1-q^{n}}, 
\end{eqnarray*}%
and the Theorem follows upon comparison with (\ref{P1m}).

\bigskip

\begin{exercise}\label{Exercise_5.3.1} 
Let $a$ be the generating state for the Heisenberg VOA $M_0$ (cf. Subsection \ref{Subsect_HeisVOA}). Prove that 
$F_{M_0}((a,z_{1}),(a,z_{2}),\tau )=P_{2}(z_{12},\tau)/\eta (\tau )$.
\end{exercise}

\begin{exercise}\label{Exercise_5.3.2} For states $u, v$ in a VOA $V$, show that $Z_V(u[0]v, q) = 0$.
\end{exercise}

\subsection{Zhu Recursion Formula II}\label{Subsect_ZhuII}
Theorem \ref{TheoremZhu} allows us to obtain a related recursion formula for
1-point functions.
\begin{theorem}
\label{Proposition_1pt} For $n\geq 1$, 
\begin{eqnarray}
Z_{V}(u[-n]v,\tau ) &=&\delta _{n,1}\mathrm{Tr}_{V}(o(u)o(v)q^{L_0-c/24}) 
\notag \\
&&+\sum\limits_{m\geq 1}(-1)^{m+1}\binom{n+m-1}{m}E_{n+m}(\tau
)Z_{V}(u[m]v,\tau ).  \notag \\
&&  \label{1ptrecur}
\end{eqnarray}
\end{theorem}
To see this, note from (\ref{FVz12}) that 
\begin{equation*}
F_{V}((u,z_{1}),(v,z_{2}),\tau )=\sum\limits_{n\in \mathbb{Z}}Z_{V}(u[-n]v,\tau )z_{12}^{-n-1}.
\end{equation*}%
Now compare this with the $z_{12}$-expansion of the rhs of (\ref{ZhuRecuruv}) using
(\ref{P1m}).  Taking $n\geq 1$ we obtain (\ref{1ptrecur}). (For $n \leq 0$ we get no information.)

\bigskip
One can apply Theorem \ref{Proposition_1pt} in a number of contexts. If we work with states $u, v, \hdots$ in $V$ that are homogeneous with respect to the square bracket Virasoro operator $L[0]$, then the $1$-point functions occurring on the rhs of (\ref{1ptrecur}) are those of states $u[m]v$ whose (square bracket) weight is \emph{strictly less} than that of $u[-n]v$
for $n \geq 1$. Thus one might hope to proceed inductively (with respect to square bracket weights) to show that $1$-point functions are holomorphic in $\mathfrak{H}$. To illustrate, we introduce the important class of VOAs $V$ of \emph{CFT-type} defined by the property that
the zero weight space $V_0$ is \emph{nondegenerate}, i.e., \emph{spanned} by the vacuum vector. This implies (Exercise \ref{Exercise_2.6.4}) that
\begin{eqnarray}\label{cfttypedef}
V = \mathbb{C}\mathbf{1} \oplus V_1 \oplus \hdots
\end{eqnarray}
Using Theorem \ref{Proposition_1pt} and the remarks following Theorem \ref{TheoremZhu}
 we obtain
\begin{lemma}\label{lemma5.6}
Suppose that $V$ is a VOA of CFT-type, and let $S$ be a generating set for $V$ as in Theorem \ref{thmexist}. Assume that Tr$_V o(u)o(v)q^{L_0-c/24}$ is holomorphic in $\mathfrak{H}$ for all $u \in S$ and $v \in V$,
and that the graded dimension $Z_V(\mathbf{1})$ is holomorphic in $\mathfrak{H}$. Then every $1$-point function for $V$ is holomorphic in $\mathfrak{H}$, and every $2$-point function for $V$ is elliptic.
\end{lemma}

\medskip
By way of example, consider the Heisenberg algebra $M_0$, which is certainly of CFT-type.  It is generated by a single state $a$ in weight $1$, and
$o(a)=0$ (cf. Exercise \ref{Exercise_4.1.3}). Furthermore $Z_{M_0}(\tau)$ is the inverse $\eta$-function 
(\ref{M0gdim}) and hence holomorphic in $\mathfrak{H}$. So the conditions of the Lemma apply to $M_0$,
so that all $1$- and $2$-point functions for $M_0$ have the desired analytic properties.
Indeed, the vanishing of the zero mode for $a$ means that in the recursion (\ref{1ptrecur}),
the anomalous first term on the rhs is not present (taking $u=a$, as we may). We get a recursion 
for $1$-point functions which may be solved with some effort, and this is how one proves 
Theorem \ref{thm4.15} and (\ref{Qvdef}). The details are described  in Subsection \ref{Subsect_Heisen1pt}.

\bigskip

\begin{exercise}\label{Exercise_5.4.1} Give the details for the Proof of Lemma \ref{lemma5.6}.
\end{exercise}

\begin{exercise}\label{Exercise_5.4.2} Show that the analysis of $1$-point and $2$-point functions associated to the 
Heisenberg algebra goes through with the same conclusions  for the Virasoro algebra $\mbox{Vir}_c$.
\end{exercise}

\newpage

\part{Modular-Invariance and Rational Vertex Operator Algebras}\label{PartII}
The representation theory of a VOA $V$, i.e., the study of $V$-modules and their characters (correlation functions) is fundamental. In this Section we introduce some of the ideas in this subject.

\section{Modules over a Vertex Operator Algebra}\label{Section_Modules}
 \subsection{Basic Definitions}\label{Subsect_ModuleDef}
 Let $V = (V, Y, \mathbf{1}, \omega)$ be a VOA of central charge $c$. 
As one might expect, a $V$-module is (roughly speaking) linear space $M$ admitting fields associated to states of $V$ which satisfy axioms analogous to those satisfied by the fields $Y(v, z)$. It is useful to introduce various types of modules, the most basic of which is the following.
\begin{definition} A weak $V$-module
 is a pair $(M, Y_M)$ where 
 \begin{eqnarray*}
&&Y_M: V \rightarrow \mathfrak{F}(M), \ v \mapsto Y_M(v, z) 
= \sum_n v^M_nz^{-n-1} \ \mbox{is a linear map},
\end{eqnarray*}
and the following hold for all $u, v \in V$, $w \in M$:
\begin{eqnarray*}
&& \ \ \ \ \  \mbox{vacuum}: Y_M(\mathbf{1}, z) = \mathrm{Id}_M\\
&&\ \ \ \ \ \ \mbox{locality}: Y_M(u, z) \sim Y_M(v, z) \\
&&  \mbox{associativity}: \mbox{for large enough} \  k, \\
&& (z_1+z_2)^kY_M(u, z_1+z_2)Y_M(v, z_2)w = (z_1+z_2)^kY_M(Y(u, z_1)v, z_2)w.
\end{eqnarray*}
\end{definition}
There is no notion of creativity or translation covariance \emph{per se} for $V$-modules.
It is not sufficient to assume only locality of operators here;  the associativity
axiom (the analog of (\ref{associativity})) is crucial. Locality and associativity
are jointly equivalent to the analog of (\ref{BI}), namely
\begin{eqnarray*}
\sum_{i \geq 0} \binom{p}{i} (u_{r+i}v)^M_{p+q-i} = \sum_{i \geq 0} (-1)^i \binom{r}{i} 
\left\{ u^M_{p+r-i}v^M_{q+i}
-(-1)^r v^M_{q+r-i}u^M_{p+i} \right\}.
\end{eqnarray*}
As before, this is the modal version of the \emph{Jacobi Identity}. For further details, see
\cite{FHL} and \cite{LL}.
 A weak $V$-module is essentially a module for a vertex algebra.
\begin{definition}\label{admVModdef} 
An \emph{admissible} $V$-module is a weak $V$-module $(M, Y_M)$
equipped with an $\mathbb{N}$-grading $M = \oplus_{n \geq 0} M_n$ such that
\begin{eqnarray}\label{modeaction2}
v \in V_k \Rightarrow v^M_n: M_m \rightarrow M_{m+k-n-1}. 
\end{eqnarray}
\end{definition}
Admissible modules are also called $\mathbb{N}$-gradable modules. Note that
(\ref{modeaction2}) is the analog of (\ref{modeaction1}). There is \emph{no} requirement
that the homogeneous spaces $M_n$ have finite dimension. An overall shift in the grading does not affect (\ref{modeaction2}), so we may, and usually shall, assume that $M_0 \not= 0$ if $M \not= 0$. We then refer to $M_0$ as the \emph{top level}.

\begin{definition}\label{Vmod} A $V$-module is a weak $V$-module $(M, Y_M)$ equipped with a grading
$M = \oplus_{\lambda \in \mathbb{C}}M_{\lambda}$ such that
\begin{eqnarray*}
&&\dim M_{\lambda} < \infty \\
&&\forall\  \lambda, M_{\lambda+n} = 0 \ \mbox{for} \ n \ll 0\\
&&L_0m = \lambda m, \ m \in M_{\lambda}
\end{eqnarray*}
\end{definition}

\medskip
We frequently call a $V$-module as in Definition \ref{Vmod} an \emph{ordinary} $V$-module if we want to emphasize that it is not merely a weak or admissible module.
There are containments 
\begin{eqnarray*}
\{\mbox{weak $V$-modules}\} \supseteq \{\mbox{admissible $V$-modules}\}
\supseteq \{\mbox{ordinary $V$-modules} \},
\end{eqnarray*}
which amounts to saying that ordinary $V$-modules can be equipped with an $\mathbb{N}$-grading making them admissible (cf. Exercise \ref{Exercise_6.1.3}).   A (weak, admissible, or ordinary) $V$-module $M$ is  \emph{irreducible} if no proper, nonzero subspace of $V$ is invariant under all modes 
$v^M_n$. More generally, we can define \emph{submodules} of $M$ in the usual way, though we will not go much  into this here.

\medskip
A VOA $V$ is \emph{ipso facto} an ordinary $V$-module in which $Y = Y_M$. It is called the \emph{adjoint} module. If the adjoint module is irreducible then we say that $V$ is \emph{simple}.
This is consistent with standard algebraic usage:  it can be shown using skew-symmetry
 that if $U \subseteq V$ is a submodule of the adjoint module $V$ then $U$ is a ($2$-sided) \emph{ideal}
in a natural sense, and that $V/U$ has a well-defined structure of VOA. See Exercise \ref{Exercise_6.1.8} for further details.

\medskip
We want to define the \emph{partition function} and \emph{character} of a $V$-module $M$ along the lines of  that for $V$ itself, as discussed in Section \ref{Section_CharsVOAs}. This makes no sense unless $M$ is equipped with a suitable grading. An important case when this can be carried through  is when $M$ is an irreducible, ordinary 
$V$-module. In this case, if $\lambda \in \mathbb{C}$ satisfies $M_{\lambda+n} \not= 0$ for some integer $n$ then  $\oplus_{n \in \mathbb{Z}} M_{\lambda + n}$ is invariant under all modes $v^M_n$ and hence coincides with $M$ thanks to irreducibility. Relabeling,  the grading on a (nonzero) irreducible, ordinary $V$-module $M$ takes the shape
\begin{eqnarray}\label{Vmodpartfunc}
M = \oplus_{n \geq 0} M_{h + n}.
\end{eqnarray}
$M_h$ is the \emph{top level} and $h$ a uniquely determined scalar  called the \emph{conformal weight} of $M$. It is an important numerical  invariant of the module.

\medskip
The 
\emph{zero mode} $o^M(v)$ for $v \in V$ is the mode of $Y_M(v, z)$ which has weight zero as an operator on $M$; it is defined because of (\ref{modeaction2}). 
We can now define the
\emph{character} $Z_M$ of an irreducible $V$-module $M$ of conformal weight $h$
 in the expected manner, namely
\begin{eqnarray}\label{chiMdef}
Z_M(v) = \mathrm{Tr}_M o^M(v)q^{L^M_0 - c/24} = q^{h-c/24}\sum_{n \geq 0} \mathrm{Tr}_{M_{h+n}} o^M(v) q^{n}.
\end{eqnarray}
The partition function of $M$ is
\begin{eqnarray*}
Z_M(\mathbf{1}) = \mathrm{Tr}_M q^{L^M_0 - c/24} = q^{h - c/24} \sum_{n \geq 0} \dim M_{h+n}q^n
\end{eqnarray*}
 where, naturally, $L^M_0$ is the corresponding zero mode for the Virasoro element.

\bigskip
The set of \emph{all} irreducible modules over the Heisenberg VOA $M_0$ is readily described.
As usual, let $a$ be the weight one state that generates $M_0$. In Subsection \ref{Subsect_HeisVOA} we 
defined, for each $h \in \mathbb{C}$, the Verma module $M_h$ and constructed a field
$a(z) \in \mathfrak{F}(M_h)$. It is more precise to denote this by $a^h(z)$. Much as in the case $h=0$,
one finds that each $M_h$ is an irreducible $M_0$-module of conformal weight $h$ with
$Y_{M_h}(a, z) = a^h(z)$. In particular, $M_0$ is a simple VOA. The Stone-von Neumann Theorem is essentially the converse:  for each $h$, $M_h$ is the unique (up to isomorphism\footnote{We have not defined morphisms of $V$-modules, but readers should be able to formulate it for themselves without diffuclty.})
irreducible module over $M_0$ of conformal weight $h$. See \cite{FLM} for a proof. The construction of the Verma module $M_h$ shows that
\begin{eqnarray}\label{Mhpartfunc}
Z_{M_h}(\mathbf{1}) = q^h/\eta(q).
\end{eqnarray}

\medskip
The characters $Z_{M_h}$ can be understood along the same lines as the special case
of $M_0$ that we described in Sections \ref{Section_CharsVOAs} and \ref{Section_EllipticFuns}. As illustrated by 
(\ref{Mhpartfunc}), results identical to Theorem \ref{thm4.15} and (\ref{Qvdef}) hold for $M_h$, except that an extra factor $q^h$ must be included. Our development of the theory of $1$- and $2$-point functions may be carried out, with essentially no change, for general $V$-modules rather than just adjoint modules. It should be pointed out, however, that the extra factor
\emph{spoils} the quasimodularity of the character values. 

\bigskip
While ordinary $V$-modules are perhaps natural, the reader may be wondering how and why admissible $V$-modules are relevant. Here we will limit ourselves here to a some general comments, and continue the discussion below. See Exercises \ref{Exercise_6.1.5}-\ref{Exercise_6.1.7} for some details, and \cite{DLM3}, \cite{Z} for complete proofs. One considers certain subspaces $O_0 \subseteq O_1 \subseteq \hdots \subseteq V$, the 
quotient spaces  $A_n(V) = V/O_n(V)$, and the inverse limit
\begin{eqnarray*}
A(V) = \lim_{\leftarrow} A_n(V).
\end{eqnarray*}
Each $A_n(V)$ has natural structure of \emph{associative algebra} such that 
the canonical projection $A_{n+1} \rightarrow A_n(V)$ is an algebra morphism. So
$A(V)$ is also an associative algebra. $A_0(V)$ is called the \emph{Zhu algebra} of $V$.

\medskip
 The representation theory of these algebras
  is intimately related to that of $V$ itself. There are functors
  \begin{eqnarray*}
\Omega_n: \mbox{Adm $V$-Mod} \rightarrow \mbox{$A_n(V)$-Mod}
\end{eqnarray*}
from the category of admissible $V$-modules to the category of $A_n(V)$-modules, and because
of the details of the construction the quotient functor $\Omega_{n}/\Omega_{n-1}$ makes sense
($\Omega_{-1}$ is trivial).
For an admissible $V$-module $M, \Omega_{n}(M)/\Omega_{n-1}(M)$ is an
$A_{n}(V)$-module that is \emph{not} the lift of an $A_{n-1}(V)$-module. 
$A_n(V)$ is designed in such a way that it acts naturally on the sum of the first $n$ graded pieces of an admissible $V$-module, and this is how the functor $\Omega_n$ is defined. It turns out that there is another functor
\begin{eqnarray}\label{Lfuncdef}
L_n: \mbox{$A_n(V)$-Mod} \rightarrow \mbox{Adm $V$-Mod}
\end{eqnarray}
which is a \emph{right inverse} of the functor $\Omega_n/\Omega_{n-1}$, and which is harder to describe. This is a key point. It is the 
existence of $L_n$ which motivates the introduction of admissible $V$-modules. $L_n$ and
$\Omega_n/\Omega_{n-1}$ induce \emph{bijections} between (isomorphism classes of) irreducible, admissible
$V$-modules and irreducible $A_n(V)$-modules which are not lifts of $A_{n-1}(V)$-modules.
For $n=0$, this is just the set of irreducible $A_0(V)$-modules. To a large extent these functors reduce the study of admissible $V$-modules to that of modules over the associative algebras
$A_n(V)$, which are more familiar objects, and they have led to a number of theoretical
advances. On the other hand, the computation of the Zhu algebra $A_0(V)$, not to mention the higher
$A_n(V)$'s,  is usually difficult. The complete structure has been elucidated in only a relatively few cases, and computer calculations have often been important.

\medskip
Needless to say, there is much more that can be said about modules over a VOA. There is a notion of
\emph{dual} module (\cite{B}, \cite{FHL} and Subsection \ref{Subsect_Adjoint}). There is also a theory of 
\emph{tensor products} of modules that is important. This is an extensive subject in its own right, and we can do no more than refer the reader to the literature (e.g., \cite{HL}, \cite{HLZ})  for further details.

\bigskip

\begin{exercise}\label{Exercise_6.1.1} Let $(M, Y_M)$ be a weak $V$-module. Prove that 
$Y_M(L_{-1}v, z) = \partial Y_M(v, z)$.
\end{exercise}

\begin{exercise}\label{Exercise_6.1.2} Show that $[L^M_m, L^M_n] = (m-n)L^M_{m-n} + (m^3-m)/12\delta_{m, -n}c\mathrm{Id}_M$.
(Thus, a weak module for $V$ is \emph{ipso facto} a module over the Virasoro algebra with the \emph{same} central charge as $V$.)
\end{exercise}

\begin{exercise}\label{Exercise_6.1.3} Show that an ordinary $V$-module $M$ is an admissible $V$-module as follows.
Let $\Lambda \subseteq \mathbb{C}$ consist of those $\lambda$ for which $M_{\lambda +k}=0$
whenever $k$ is a negative integer, and let $M_n = \oplus_{\lambda \in \Lambda} M_{\lambda +n}.$ 
Show that $M= \oplus_{n \geq 0} M_{n}$ is an $\mathbb{N}$-grading on $M$ satisfying 
(\ref{modeaction2}).
\end{exercise}

\begin{exercise}\label{Exercise_6.1.4} Give a complete proof that the Verma modules $M_h$ are
irreducible modules over the Heisenberg algebra $M_0$.
\end{exercise}

\begin{exercise}\label{Exercise_6.1.5} Let $M$ be an admissible $V$-module. Prove that for each $v \in V$,
the zero mode $o^M(L[-1]v)$ \emph{annihilates} $M$. (Use Exercises \ref{Exercise_6.1.1} and \ref{Exercise_2.7.1}.)
\end{exercise}

\begin{exercise}\label{Exercise_6.1.6} For $n \geq 0, u \in V_k, v \in V$ define
$u \circ_n v = \mathrm{Res}_z Y(u, z)v\frac{(1+z)^{k+n}}{z^{2n+2}}$.
Let $O_n(V)$ be the \emph{span} of all states $u \circ_n v$ and  $L[-1]u$.\\
(a) Prove that if $n=0$, the span of the states $u \circ_0 v$ already contains $L[-1]u$.\\
(b) Prove that $O_0(V) \subseteq O_1(V) \subseteq \hdots$.
\end{exercise}

\begin{exercise}\label{Exercise_6.1.7} With the notation of the previous Exercise, introduce the product
$u *_n v = \sum_{m=0}^n \binom{ m+n}{n}\mathrm{Res}_zY(u, z)v\frac{(1+z)^{k+n}}{z^{n+m+1}}.$\\
(a) Show that $O_n(V)$ is a $2$-sided ideal with respect to the product $*_n$.\\
(b) Show that $*_n$ induces a structure of associative algebra on the quotient space 
$A_n(V) = V/O_n(V)$.
\end{exercise}

\begin{exercise}\label{Exercise_6.1.8} $V$ is a VOA and $U \subseteq V$ a submodule of the adjoint module, so that
$v_n u \in U$ for all $u \in U, v \in V, n \in \mathbb{Z}$. Prove that $u_nv \in V$, and deduce that
if $U \not= V$ then $V/U$ inherits the structure of VOA.
\end{exercise}

 \subsection{$C_2$-Cofinite, Rational and Regular Vertex Operator Algebras}\label{Subsect_C2cofinite}
 We are going to focus  on some important classes of VOAs $V$ which have the property that they have only \emph{finitely many} (inequivalent) irreducible modules. The reader might well be surprised that there are any such VOAs at all beyond those of finite dimension (cf. Exercises \ref{Exercise_2.3.5} and \ref{Exercise_2.3.6}). We will also make the simplifying assumption that $V$ is of \emph{CFT-type} 
  (\ref{cfttypedef}) throughout the rest of these Notes, although for many of the results to be discussed this assumption is not necessary. 
 
 \begin{definition}\label{finiteVOAs} 
 (a) $V$ is \emph{rational} if every admissible
$V$-module is completely reducible, i.e. a direct sum of irreducible, admissible $V$-modules.\\
(b) $V$ is \emph{regular} if every weak $V$-module is a direct sum of irreducible, ordinary $V$-modules.\\
(c)  $V$ is  \emph{$C_2$-cofinite} if the graded subspace $C_2(V) = \langle u_{-2}v \ | \ u, v \in V \rangle$ has finite codimension in $V$.
\end{definition}

Based on what we said in the previous Subsection, it is easy to see that a regular VOA $V$ is a rational VOA. Indeed, an admissible $V$-module is a weak module,
hence a direct sum of irreducible, ordinary modules and \emph{ipso facto} a direct sum of
irreducible admissible modules. It is also known \cite{ABD}, \cite{Li} that regularity is
\emph{equivalent} to the conjunction of rationality and $C_2$-cofiniteness.

\medskip
While (a) and (b) of Definition \ref{finiteVOAs} both assert that certain module categories are semisimple,
(c) is rather different. (a) and (b) are \emph{external conditions} that can be difficult to verify,
whereas (c) is an \emph{internal} condition that is easier to deal with. On the other hand, regular VOAs have better modular-invariance properties than those which are $C_2$-cofinite.  

\begin{theorem}\label{thmC2VOA}
Suppose that $V$ is a $C_2$-cofinite VOA. Then the following hold.\\
(a) Each $A_n(V)$ is \emph{finite-dimensional}. \\
(b) Every weak $V$-module is an admissible module.\\
(c) $V$ has only finitely many isomorphism classes of irreducible, admissible modules.
\end{theorem}
Note that for a \emph{finitely generated} VOA $V$, (b) is \emph{equivalent} to $C_2$-cofiniteness.

\medskip
For further discussion of (a), see \cite{Z}, \cite{My1}, \cite{GN}, \cite{Bu}; (b) is proved in \cite{My1}. The approach in \cite{GN} produces a sort of weak analog of the PBW Theorem  in Lie theory (cf. Appendix) which applies to weak modules. This idea is very useful, and is used in \cite{ABD}, \cite{My1}, \cite{Bu} and elsewhere in the literature.
(c) follows from (a)  and the properties of the functors $L_n$ and $\Omega_n$ discussed in Subsection \ref{Subsect_ModuleDef}.

\medskip
The following omnibus result (\cite{DLM1}, \cite{DLM3}) collects some of the main facts about rational VOAs. 
\begin{theorem}\label{thmratVOA}
Suppose that $V$ is a rational VOA. Then the following hold.\\
(a) $A_0(V)$ is \emph{semisimple}. \\
(b) Each $A_n(V)$ is \emph{finite-dimensional}.\\
(c) $V$ has only finitely many isomorphism classes of irreducible, admissible $V$-modules, \\
(d) every irreducible, admissible $V$-module is an ordinary $V$-module.
\end{theorem}
Note that (b) is \emph{equivalent} to rationality (loc cit).

\medskip
 Whether a rational VOA is necessarily $C_2$-cofinite is presently one of the main open questions in the representation theory of VOAs. If this is so, then there would be no difference between rational and regular VOAs. In the early history of VOA theory it was possible to believe that rationality and $C_2$-cofiniteness were \emph{equivalent}.
 That, however, has turned out to be a chimera. There are VOAs which are $C_2$-cofinite but have admissible (in fact ordinary) modules which are \emph{not} completely reducible. These are 
 \emph{logarithmic field theories}, a name that we will justify in Section \ref{Section_VOAModInvariance}.

\bigskip

\begin{exercise}\label{Exercise_6.2.1} Give two proofs that the Heisenberg VOA $M_0$ is \emph{not} a rational VOA:
(a) by using Theorem \ref{thmratVOA}, and (b) by explicitly constructing an admissible
$M_0$-module that is \emph{not} completely reducible.
\end{exercise}

\begin{exercise}\label{Exercise_6.2.2} For any VOA $V$, show that the quotient space $P(V) = V/C_2(V)$ carries the structure of a \emph{Poisson algebra} in the following sense: the products $\{u, v \} =u_0 v, uv = u_{-1} v$ afflict
$P(V)$ with (well-defined) structures of Lie algebra and commutative, associative algebra respectively,
moreover $ \{uv, w\} = u\{v, w\} +  \{u, v\}w$.
\end{exercise}

\begin{exercise}\label{Exercise_6.2.3} Calculate the Poisson algebra $P(M_0)$ associated to the Heisenberg VOA.
\end{exercise}

\section{Examples of Regular Vertex Operator Algebras}\label{Section_ExamplesVOAS}
It is time to describe some further examples of VOAs beyond the Heisenberg and Virasoro theories.
In particular, we want to have available a selection of regular VOAs. Our examples are fairly standard, but require some effort to construct. For this reason, we will mainly limit ourselves to a description of the underlying Fock spaces and generating fields.  

\subsection{Vertex Algebras Associated to Lie Algebras}\label{Subsect_LieVOA}
The reader might want to look over Appendix $1$ before reading this Subsection. We can construct a VOA from a pair $(\mathfrak{g}, ( \ ,\ ))$
consisting of a Lie algebra $\mathfrak{g}$ equipped with a symmetric, invariant,  bilinear form
$( \ , \ ): \mathfrak{g} \otimes \mathfrak{g} \rightarrow \mathbb{C}$. The details amount to an elaboration of the
case of the Heisenberg algebra discussed in Subsection \ref{Subsect_HeisVOA}, which is the $1$-dimensional case. The \emph{affine Lie algebra} or \emph{Kac-Moody algebra}
associated to $(\mathfrak{g}, ( \ , \ ))$ is the linear space
\begin{eqnarray*}
\hat{\mathfrak{g}} = \mathfrak{g} \otimes \mathbb{C}[t, t^{-1}] \oplus \mathbb{C}K = \oplus_{n} \mathfrak{g} \otimes t^n \oplus \mathbb{C}K
\end{eqnarray*}
with brackets
\begin{eqnarray*}
[a \otimes t^m, b \otimes t^n] &=& [a, b] \otimes t^{m+n} + m (a, b) \delta_{m, -n}K, \ \ (a, b \in \mathfrak{g})\\
{[}\hat{\mathfrak{g}}, K {]} &=& 0.
\end{eqnarray*}

\medskip
$\hat{\mathfrak{g}}$ has a triangular decomposition with
$\hat{\mathfrak{g}}^{\pm} = \oplus_{\pm n > 0} \mathfrak{g} \otimes t^n$ and
$\hat{\mathfrak{g}}^0 = \mathfrak{g} \oplus \mathbb{C}K$. Here and below, we  identify
$\mathfrak{g}$ with $\mathfrak{g} \otimes t^0$. Fix a $\mathfrak{g}$-module $X$ and a scalar $l$. We 
extend $X$ to a $\hat{\mathfrak{g}}^+ \oplus \hat{\mathfrak{g}}^0$-module by letting $\hat{\mathfrak{g}}^+$  
\emph{annihilate} $X$ and letting $K$ act as multiplication by $l$ called the
\emph{level}.
We have the induced $\hat{\mathfrak{g}}$-module
\begin{eqnarray}\label{indgmoddef}
V_{\mathfrak{g}}(l, X) = \mbox{Ind}(X)
\end{eqnarray}
(notation as in (\ref{indmoddef})).
Following the Heisenberg case discussed in Subsection \ref{Subsect_HeisVOA}, 
we can define fields on $V_{\mathfrak{g}}(l, X)$ for each $a \in \mathfrak{g}$
by setting
\begin{eqnarray*}
Y_{V_{\mathfrak{g}}(l, X)}(a, z) = \sum_n a_n z^{-n-1}
\end{eqnarray*}
where $a_n$ is the induced action of $a \otimes t^n$. As in (\ref{locality1})
we obtain
\begin{eqnarray*}
&& \sum_{j = 0}^2 (-1)^j\binom{2}{j}[a_{2-j-r}, b_{j-s}] \\
&=& \sum_{j = 0}^2 (-1)^j\binom{2}{j}({[}a, b {]}_{2-r-s}+ (2-j-r)(a, b) 
l\delta_{2-j-r, s-j}\mathrm{Id})   \\
&=& \left\{ (2-r)
-2(1-r) -r  \right \}(a, b)l \delta_{r+s, 2}\mathrm{Id} = 0,
\end{eqnarray*}
so that the  fields $\{ Y_{V_{\mathfrak{g}}(l, X)}(a, z) \ | \  a \in \mathfrak{g} \}$
are mutually local of order two. Taking $X = \mathbb{C}\mathbf{1}$ to be the
\emph{trivial} $1$-dimensional $\mathfrak{g}$-module,  one shows via Theorem \ref{thmexist}
 that the corresponding fields generate a vertex algebra with Fock space $V_{\mathfrak{g}}(l, \mathbb{C}\mathbf{1})$.  Moreover, each $V_{\mathfrak{g}}(l, X)$ is
an admissible module.

\medskip
To describe a conformal vector  in $V_{\mathfrak{g}}(l, \mathbb{C}\mathbf{1})$ and thereby obtain the structure of VOA, it is convenient at this point to specialize to the case that $\mathfrak{g}$ is a \emph{finite-dimensional, simple Lie algebra} of dimension $d$, say. 
We will also take $( \ , \ )$ to be the  \emph{Killing form}, appropriately normalized.\footnote{The normalization is an important detail, of course, but we will not need it.}
 Note that this takes us out of the regime of the Heisenberg theory, to which we return in Subsection \ref{Subsect_LatticeVOA}. An approach that covers both cases is described in \cite{LL}. With our assumptions, one shows that
\begin{eqnarray}\label{Sugawaradef}
\omega =\frac{1}{2} \frac{1}{l+h^{\vee}} \sum_{i=1}^d u_i(-1)u^i
\end{eqnarray}
is the desired conformal vector with central charge $c=\frac{ld}{l+h^{\vee}}$. Here, $\{u_i\}$ is a basis of $\mathfrak{g}$, $\{u^i \}$ the basis dual to 
$\{u_i\}$ with respect to the form $( \ , \ )$, and $h^{\vee}$  the 
\emph{dual Coxeter number} of $\mathfrak{g}$. This is usually called the
 \emph{Sugawara construction}. Needless to say, we must also assume that $l+ h^{\vee} \not= 0$. 

\medskip
The $L_0$-grading on $V_{\mathfrak{g}}(l, \mathbb{C}\mathbf{1})$ that obtains from the Sugawara construction is the natural one in which
the state $a_n\mathbf{1}$ has weight $-n$ for $a \in \mathfrak{g}$ and $n \leq 0$. In particular
the zero weight space is $V_{\mathfrak{g}}(l, \mathbb{C}\mathbf{1})_0 = \mathbb{C}\mathbf{1}$, and
the VOA is of CFT-type.
Because an ideal in the adjoint module is a graded submodule  (cf. the discussion in Subsection \ref{Subsect_ModuleDef}), any \emph{proper} ideal necessarily lies in $\oplus_{n \geq 2} V_{\mathfrak{g}}(l, \mathbb{C}\mathbf{1})_n$. It follows that there is a \emph{unique maximal} proper ideal, call it $J$, and the quotient space
\begin{eqnarray*}
L_{\mathfrak{g}}(l, 0) = V_{\mathfrak{g}}(l, \mathbb{C}\mathbf{1})/J
\end{eqnarray*}
is a simple VOA. 

\medskip
More generally,  take $X$ to be a finite-dimensional irreducible $\mathfrak{g}$-module. As such it is a highest-weight module $L({\lambda})$ indexed by an element $\lambda$ in the weight lattice of $\mathfrak{g}$. The top level
of $V_{\mathfrak{g}}(l, L({\lambda}))$ is naturally identified with $L({\lambda})$, and because 
this is an irreducible $\mathfrak{g}$-module then there is a unique maximal proper
submodule $J \subseteq V_{\mathfrak{g}}(l, L({\lambda}))$ (considered as $\hat{\mathfrak{g}}$-module).
The quotient spaces
\begin{eqnarray*}
L_{\mathfrak{g}}(l, \lambda) = V_{\mathfrak{g}}(l, L({\lambda}))/J
\end{eqnarray*}
are ordinary, irreducible $V_{\mathfrak{g}}(l, \mathbb{C}\mathbf{1})$-modules, and they are inequivalent for distinct choices of highest weight $\lambda$. Thus the VOA 
$V_{\mathfrak{g}}(l, \mathbb{C}\mathbf{1})$ has infinitely many inequivalent ordinary, irreducible modules, and in particular it cannot be rational (Theorem \ref{thmratVOA}).
Concerning the question of regularity of these VOAs, we collect the main facts  (\cite{FZ}, \cite{DL}, \cite{DM1}, \cite{DLM2}, \cite{DLM4}):
\begin{theorem}\label{thmaffgVOA}
Let $\mathfrak{g}$ be a finite-dimensional simple Lie algebra. The simple VOA
$L_{\mathfrak{g}}(l, 0)$ is \emph{rational} if, and only if, $l$ is a \emph{positive integer}. In this case
it is \emph{regular}, and the ordinary, irreducible modules are the spaces
$L_{\mathfrak{g}}(l, \lambda)$ where $\lambda$ satisfies $\lambda(\theta) \leq l$ and $\theta$ is the longest
positive root.
\end{theorem}
These theories are called WZW models in the physics literature.

\subsection{Discrete Series Virasoro Algebras}\label{Subsect_DiscreteVirasoro}
Here we discuss some quotients of Virasoro VOAs $\mbox{Vir}_c$ (cf. Theorems \ref{VirVA} and \ref{VirVOA})
that turn out to be regular. As in the last Subsection, it is the underlying Lie structure that makes the calculations manageable. The details are quite different, however, and depend on the \emph{Kac determinant} (e.g. \cite{KR}) and the structure of the Verma modules (\ref{VirVerm}) $M_{c, h}$ over the Virasoro algebra (these are $\mbox{Vir}_c$-modules) (\cite{FF}).  There is no space to describe these results systematically here, although we discuss some examples of Kac determinants in Subsection 
\ref{Subsect_Invariant bilinear}. So we give less detail compared to the WZW models. The theories we are going to describe in this Subsection find important applications in the physics of phase transitions and critical phenomena.
See \cite{FMS} for further background.

\medskip
 The Virasoro VOA $\mbox{Vir}_c$ may, or may not, be a simple VOA, but there is a \emph{unique} maximal proper submodule $J$ and $L_c = \mbox{Vir}_c/J$ is a simple vertex operator algebra  of central charge $c$. It turns out that $\mbox{Vir}_c$ is \emph{never} rational (cf. Exercise \ref{Exercise_7.2.2}). As for the rationality of $L_c$, we have the following omnibus result:
 
 \begin{theorem}\label{thmratVir} The following are equivalent:\\
 (a)  $L_c$ is a rational VOA. \\
 (b) $J \not= 0$. \\
 (c) $c$ lies in the so-called \emph{discrete series}, i.e., there are coprime integers $p, q \geq 2$
such that
\begin{eqnarray}
c = c_{pq} = 1- \frac{6(p-q)^2}{pq}.\label{cpq}
\end{eqnarray}
In this case $L_c$ is regular, the conformal weights of the ordinary irreducible modules are
\begin{equation*}
h_{r, s} = \frac{(pr - qs)^2-(p-q)^2}{4pq}, \ 1 \leq r \leq q-1, 1\leq s  \leq p-1
\end{equation*}
(taking only one value of $h$ for each pair $h_{r, s}, h_{q-r, p-s}$), and two ordinary irreducible modules are isomorphic if, and only if, they have the same conformal weight.\footnote{Generally, a VOA may have inequivalent irreducible modules with the \emph{same} conformal weight.} Thus there are just
$(p-1)(q-1)/2$ inequivalent ordinary irreducible modules over $L_c$.
\end{theorem}
See \cite{Wa} for the proof of rationality (also \cite{DMZ}), where the idea is to compute the Zhu algebra
$A_0(L_c)$. Regularity is shown in \cite{DLM2}. The origin of the values $c_{p, q}$ is discussed in
Subsection \ref{Subsect_Invariant bilinear}

\medskip
Apart from the trivial case when $p=2, q=3$, the two `smallest' cases, i.e., those with the fewest number of ordinary irreducible modules,
correspond to $(p, q) = (2,5)$ and $(3,4)$. In the first case (the \emph{Yang-Lee} model in physics)
we have $c=-22/5$ and conformal weights $0, -1/5$.
In the second case (the \emph{Ising model}) $c =1/2$ with conformal weights
$0, 1/2, 1/16$.

\bigskip

\begin{exercise}\label{Exercise_7.2.1} Prove that $\mbox{Vir}_c$ has a unique maximal proper submodule $J$.
\end{exercise}

\begin{exercise}\label{Exercise_7.2.2} Suppose that $J \not= 0$. Prove that  $\mbox{Vir}_c$ 
is \emph{not} a rational VOA.
\end{exercise}

\begin{exercise}\label{Exercise_7.2.3} Opine on the statement that the case $p=2, q=3$ is `trivial'.
\end{exercise}

\medskip

\subsection{Lattice Theories}\label{Subsect_LatticeVOA}
Lattice theories (\cite{B}, \cite{FLM}) are VOAs whose connections with Lie algebras are of lesser importance compared to the examples in the last two Subsections. Their basic properties are of a more combinatorial nature, and reflect features that one may expect in  general rational theories. Because of this and the fact that they are amenable to computation, lattice theories occupy a central position in current VOA theory.

\medskip
Let $d$ be a positive integer and $\mathfrak{h} = \mathbb{C}^d$ a rank $d$ linear space equipped with a nondegenerate symmetric bilinear form $( \ , \ )$. Consideration of $\mathfrak{h}$ as an \emph{abelian} Lie algebra leads to the affine algebra $\hat{\mathfrak{h}}$
as in Subsection \ref{Subsect_LieVOA}. Let
\begin{eqnarray}\label{rankdHeisdef}
M_0^d = V_{\mathfrak{h}}(1, \mathbb{C}\mathbf{1})
\end{eqnarray}
be the corresponding vertex algebra of \emph{level $1$}. The conformal vector $\omega$ is defined as in (\ref{Sugawaradef}) with $l= h^{\vee} = 1$. The resulting VOA has central charge $c=d$.
This is nothing more than a slightly different approach to the \emph{rank $d$ Heisenberg VOA}, as discussed in Subsection \ref{Subsect_HeisVOA} (cf.  Exercise \ref{Exercise_7.3.1}).

\medskip
The irreducible $\mathfrak{h}$-modules are $1$-dimensional and indexed by a 
weight  in the dual space of $\mathfrak{h}$. Identifying $\mathfrak{h}$ with its dual via $( \ , \ )$,
we obtain $M_0^d$-modules (\ref{indgmoddef}) with underlying linear space
\begin{eqnarray*}
V_{\mathfrak{h}}(1, \alpha) = S(\hat{\mathfrak{h}}^-) \otimes e^{\alpha} \ (\alpha \in \mathfrak{h}).
\end{eqnarray*}
Here,  $\mathbf{1} \otimes e^{\alpha}$ (or just $e^{\alpha}$) is notation for the spanning vector of the ($1$-dimensional) top level
of $V_{\mathfrak{h}}(1, \alpha)$ and
\begin{eqnarray}\label{betaaction}
\beta \otimes e^{\alpha} = \beta_0.e^{\alpha} = (\beta, \alpha)1 \otimes e^{\alpha}, \ \beta \in \mathfrak{h} = \mathfrak{h} \otimes t^0.
\end{eqnarray}

\medskip
In order to describe the Fock spaces of lattice theories we need a bit more structure.
Namely, we assume that $(\mathfrak{h}, ( \ , \ ))$ is the scalar extension of a Euclidean space.
Thus, $E = \mathbb{R}^d = \mathfrak{h}_{\mathbb{R}}$ is a real space equipped with a \emph{positive-definite quadratic form} $Q: E \rightarrow \mathbb{R}$, $\mathfrak{h} = \mathbb{C} \otimes_{\mathbb{R}} E$,
and $( \ , \ )$ is the $\mathbb{C}$-linear extension of the bilinear form on $E$ defined by $Q$, also denoted by $( \ , \ )$. In particular, $Q(\alpha) = (\alpha, \alpha)/2$ for $\alpha \in E$.
A \emph{lattice} $L \subseteq E$ is the additive subgroup spanned by a
\emph{basis} of $E$. $L$ is an \emph{even} lattice if $(\alpha, \alpha) \in 2\mathbb{Z}$ for all
$\alpha \in L$, i.e., the restriction of $Q$ to $L$ is integral. 

\medskip
For an even lattice $L \subseteq E$ we introduce the linear space
\begin{eqnarray*}
V_L = \oplus_{\alpha \in L} V_{\mathfrak{h}}(1, \alpha).
\end{eqnarray*}
Identifying $\oplus_{\alpha} \mathbb{C}e^{\alpha}$ with the \emph{group algebra}\footnote{We 
only explicitly use the linear structure of $\mathbb{C}[L]$, although the algebra structure also
plays a r\^{o}le.}
$\mathbb{C}[L]$ of the lattice, we can write (\ref{VLdef}) more compactly as
\begin{eqnarray}\label{VLdef}
V_L = S(\hat{\mathfrak{h}}^-) \otimes \mathbb{C}[L].
\end{eqnarray}

There is a natural grading on $V_L$ that turns out to be the one defined by the
$L_0$ operator. We take the tensor product grading on (\ref{VLdef}) in which
$S(\hat{\mathfrak{h}}^-)$ has the grading of the Fock space of the rank $d$ Heisenberg algebra that it is,
and where
$e^{\alpha}$ has weight $Q(\alpha)$. Using (\ref{M0dgdim}), the partition function of $V_L$ is
\begin{eqnarray}\label{VLpartfuncdef}
Z_{V_L}(\mathbf{1}) = \frac{\sum_{\alpha \in L} q^{Q(\alpha)}}{\eta(q)^d}.
\end{eqnarray}
The numerator here is the \emph{theta function} of $L$, a topic to which we shall return in
Section \ref{Section_VectorModForm}.

\medskip
So far then, we have described the Fock space $V_L$ as a sum of Heisenberg modules.
We define $Y(v, z)$ for $v \in M_0^d$ to be the operator
whose restriction to $V_{\mathfrak{h}}(1, \alpha)$ is just $Y_{V_{\mathfrak{h}}(1, \alpha)}(v, z)$. 
In order to impose the structure of VOA
on $V_L$, we must construct fields for all of the states in the Fock space (\ref{VLdef}).
Because of Theorem \ref{thmexist} it suffices to define $Y(e^{\alpha}, z)$ for $\alpha \in L$ and establish locality, but nothing that has come so far has prepared us for this. The generating fields we have considered in detail for the Heisenberg, WZW and Virasoro theories have modes $a_n$ that are closely related to some Lie algebra, but in theories such as $V_L$ this will generally not be the case. We 
content ourselves with the prescription for $Y(e^{\alpha}, z)$, referring the reader to \cite{FLM}, \cite{K1} for further background and motivation:
\begin{eqnarray}\label{Yealphadef}
Y(e^{\alpha}, z) = \exp\left(\sum_{n > 0} \frac{\alpha_{-n}}{n}z^{n} \right) \exp\left( \sum_{n < 0} \frac{\alpha_{-n}}{n}z^{n} \right)e^{\alpha}z^{\alpha}.
\end{eqnarray}
Beyond the modes $\alpha_n$ of $Y(\alpha, z)$, $z^{\alpha}$ is a shift operator
$z^{\alpha}: v \otimes e^{\beta} \mapsto z^{(\alpha, \beta)}v \otimes e^{\beta} \ (v \in \hat{\mathfrak{h}}^-)$,
and $e^{\alpha}: v \otimes e^{\beta} \mapsto \epsilon(\alpha, \beta) v \otimes e^{\alpha + \beta}$
for a certain bilinear $2$-cocycle $\epsilon: L \otimes L \rightarrow \{ \pm 1 \}$ (loc cit).

\medskip
The ordinary, irreducible modules over $V_L$ are constructed in \cite{Do}. The underlying Fock spaces
are very similar to (\ref{VLdef}), and are indexed by the \emph{cosets} of $L$ in its
$\mathbb{Z}$-dual $L^0$ (cf. Exercise \ref{Exercise_7.3.7}). Precisely, they are 
\begin{eqnarray*}
V_{L+ \lambda} = \oplus_{\alpha \in L} V_{\mathfrak{h}}(1, \alpha + \lambda) = \hat{\mathfrak{h}}^- \otimes \mathbb{C}[L+\lambda]
\end{eqnarray*}
for $\lambda \in L^0$, with partition functions
\begin{eqnarray}\label{VLMpartfuncdef}
Z_{V_{L+\lambda}}(\mathbf{1}) = \frac{\sum_{\alpha \in L} q^{Q(\alpha + \lambda)}}{\eta(q)^d}.
\end{eqnarray}
The fields $Y_{V_{L+ \lambda}}(v, z)$ are similarly analogous to (\ref{Yealphadef}) (loc cit).
Indeed, one can usefully combine \emph{all} of these fields and Fock spaces into a bigger and better edifice. For this, see \cite{DL}. For rationality and $C_2$-cofiniteness, see \cite{Do}
and \cite{DLM4} respectively. Summarizing,

\begin{theorem}\label{thmVLVOA} Let $L$ be an even lattice. Then $V_L$ is
a regular VOA, and its ordinary, irreducible modules are the Fock spaces $V_{L+ \lambda}$.
It thus has just $|L^0:L|$ distinct ordinary, irreducible modules.
\end{theorem}

\bigskip
In the following Exercises, $L \subseteq E$ is an even lattice in Euclidean space as above.

\medskip
\begin{exercise}\label{Exercise_7.3.1} Show that the VOA  (\ref{rankdHeisdef}) is \emph{isomorphic}
to the tensor product $M_0^{\otimes d}$ of $d$ copies of the Heisenberg VOA $M_0$
(cf. Exercise \ref{Exercise_2.6.7}).
\end{exercise}

\begin{exercise}\label{Exercise_7.3.2} Show that $M_0^d$ is a \emph{simple} VOA.
\end{exercise}

\begin{exercise}\label{Exercise_7.3.3} Verify that if $\alpha \in L$ then $\mathbf{1}\otimes e^{\alpha}$ has
$L_0$-weight $Q(\alpha)$.
\end{exercise}

\begin{exercise}\label{Exercise_7.3.4} In the definition of $V_L$, what is the purpose of requiring $L$ to
be an \emph{even} lattice? What about positive-definiteness?
\end{exercise}

\begin{exercise}\label{Exercise_7.3.5} Let $\alpha \in L$.\\
(a) Prove that $Y(e^{\alpha}, z)$ is a creative field in $\mathfrak{F}(V_L)$. \\
(b) Prove that $Y(e^{\alpha}, z)$ and $Y(v, z)$ are mutually local ($v \in \hat{\mathfrak{h}}^-$).
\end{exercise}

\begin{exercise}\label{Exercise_7.3.6} Let $L$ be an even lattice with $L_0 = \{ \alpha \in L \ | \ Q(\alpha) = 1\}$.
Prove that $L_0$ is a \emph{semisimple root system} with components
of type $ADE$.
\end{exercise}

\begin{exercise}\label{Exercise_7.3.7} The \emph{dual lattice} of $L$ is defined via 
\begin{eqnarray*}
L^0 = \{ \beta \in E \ | \ (\alpha, \beta) \in \mathbb{Z} \ 
\forall\  \alpha \in L\}.
\end{eqnarray*}
Prove that $L \subseteq L^0$ is a subgroup of \emph{finite index}.
\end{exercise}

\begin{exercise}\label{Exercise_7.3.8} Let $\mathfrak{g}$ be a finite-dimensional simple Lie algebra of type
$ADE$. \\
(a) Show that the WZW model $L_{\mathfrak{g}}(1, 0)$ of level $1$ is (isomorphic to)
the lattice theory $V_L$ where $L$ is the root lattice associated to $\mathfrak{g}$. \\
(b) Compute the number of inequivalent ordinary, irreducible modules over
$L_{\mathfrak{g}}(1, 0)$ both by using Theorem \ref{thmaffgVOA}, and by using Theorem \ref{thmVLVOA}.
\end{exercise}

\begin{exercise}\label{Exercise_7.3.9} Let $L_1, L_2$ be a pair of even lattices.\\
(a) Show that the \emph{orthogonal direct sum} $L_1 \perp L_2$ is an even lattice.\\
(b) Prove that $V_{L_1\perp L_2} \cong V_{L_1} \otimes V_{L_2}$ (cf. Exercise \ref{Exercise_2.6.7}).
\end{exercise}

\section{Vector-Valued Modular Forms}\label{Section_VectorModForm}
In order to formulate \emph{modular-invariance} for $C_2$-cofinite and regular
VOAs,  the idea of a \emph{vector-valued modular form} is useful. This generalizes the theory of
modular forms that we discussed in  Section \ref{Section_ModForms}, and includes as a special case the theory of modular forms on a finite-index subgroup of $SL_2(\mathbb{Z})$. We use the notation of Section \ref{Section_ModForms}.

\subsection{Basic Definitions}\label{Subsect_VVModformDef}
Fix an integer $k$ and let $\mathfrak{F}_k$ be the space of holomorphic functions\footnote{We could equally well deal with \emph{meromorphic functions}.} in $\mathfrak{H}$ regarded as a right $\Gamma$-module
with respect to the action defined in (\ref{GHaction}), (\ref{slashaction}). A \emph{weak  vector-valued modular form of weight $k$} may be taken to be a finite-dimensional $\Gamma$-submodule 
$V \subseteq \mathfrak{F}_k$.
Let\footnote{Superscript $t$ denotes \emph{transpose}.} 
$F(\tau) = (f_1(\tau), \hdots, f_p(\tau))^t$ where the component functions $f_i(\tau)$ are a set of (not necessarily linearly independent) generators for $V$. There is then a representation
$\rho: \Gamma \rightarrow GL_p(\mathbb{C})$ such that
\begin{eqnarray}\label{vvformdef}
\rho(\gamma)F(\tau) = F|_k \gamma(\tau), \ \gamma \in \Gamma,
\end{eqnarray}
where $|_k$ is the obvious extension of the stroke operator to vector-valued functions.
We also call the pair $(F, \rho)$ a weak vector-valued modular form of weight $k$. Given a pair 
$(F, \rho)$ satisfying (\ref{vvformdef}),  we recover $V$ as the span of the component functions of $F(\tau)$. The classical modular forms of Section \ref{Section_ModForms} correspond to the case when $\rho$ is the trivial $1$-dimensional representation of $\Gamma$.

\medskip
To describe the extension of (\ref{genqexp}) to the vector-valued case, decompose
$V$ into a direct sum of $T$-invariant indecomposable subspaces
\begin{eqnarray*}
V = V_1 \oplus \hdots \oplus V_r
\end{eqnarray*}
corresponding to the Jordan decomposition of the action $T: f(\tau) \mapsto f(\tau +1)$.
The characteristic polynomial on $V_i$ is $(x-e^{2 \pi i \mu_i})^{\dim V_i}$.
 The basic fact is 
 \begin{theorem}\label{thmlogqexp}
  There are $q$-expansions 
 $g_j(\tau) = q^{\mu_i}\sum_{n \in \mathbb{Z}}a_{ijn}q^n, (0 \leq j \leq n_i -1)$ such that
 the functions
 \begin{eqnarray}\label{logqexpdef}
g_0(\tau) + g_1(\tau)\log q + \hdots + g_m(\tau)(\log q)^m, \ 0 \leq m \leq n_i-1,
\end{eqnarray}
are a \emph{basis} of $V_i$. In particular, $V$ has a basis of functions of this form.
\end{theorem}
 We call (\ref{logqexpdef}) a \emph{logarithmic}, 
or \emph{polynomial}\footnote{We may rewrite (\ref{logqexpdef})
using powers of $\tau$, or other polynomials in $\tau$, instead of powers of $\log q$.}
$q$-expansion. 

\medskip
Suppose that $(F, \rho)$ is a weak vector-valued modular form. Then the component functions of $F(\tau)$
are linear combinations of polynomial $q$-expansions (\ref{logqexpdef}). We say that
$(F, \rho)$, or simply $F(\tau)$, is \emph{almost holomorphic} if the component functions are holomorphic in $\mathfrak{H}$ and if the $q$-expansions $g_j(\tau)$ are \emph{left-finite} 
or \emph{meromorphic at $\infty$}, i.e., for all
$i, j$ the Fourier coefficients $a_{ijn}$ \emph{vanish} for $n \ll 0$.  Similarly, $F(\tau)$ is \emph{holomorphic} if it is almost
holomorphic and if $a_{ijn}=0$ whenever $\Re(\mu_i)+n<0$. 
These definitions are independent of  the choice of $g_j(\tau)$.

\medskip
Fix an integer $N \geq 1$.  We set
\begin{eqnarray*}
\Delta(N) = \langle \gamma T^N \gamma^{-1} \ | \ \gamma \in \Gamma \rangle.
\end{eqnarray*}
This is the smallest normal subgroup of $\Gamma$ that contains $T^N$. We say that a
subgroup $G \subseteq \Gamma$ has level $N$ if $\Delta(N) \subseteq G$. A representation 
$\rho: \Gamma \rightarrow GL_p(\mathbb{C})$ has \emph{level $N$} if ker$\rho$ 
has level $N$ (equivalently, $\rho(T)$ has finite order dividing $N$). 
A vector-valued modular form
$(F, \rho)$ has level $N$ if $\rho$ has level $N$. Now recall that finite-order operators are
\emph{diagonalizable}. It follows from Theorem \ref{thmlogqexp} that if $(F, \rho)$ has level $N$
  then the component functions of $F(\tau)$ have $q$-expansions that are free of logarithmic terms.
Indeed, the eigenvalues of $\rho(T)$ are $N$th. roots of unity, so that the $q$-expansions
(\ref{logqexpdef}) reduce to a single $q$-expansion of the form
\begin{eqnarray}\label{gjdef}
g_j(\tau) = q^{r/N}\sum_{n \geq 0} a_{jn}q^n
\end{eqnarray}
for some integer $r$.

\medskip
The \emph{principle congruence subgroup of level $N$} is the subgroup of $\Gamma$ defined by
\begin{eqnarray*}
\Gamma(N) = \{ \gamma \in \Gamma \ | \ \gamma \equiv I_2 \ (\mbox{mod} \ N)\}.
\end{eqnarray*}
We have $\Delta(N) \unlhd \Gamma(N) \unlhd \Gamma$. While $\Gamma(N)$ always has finite index in
$\Gamma$, $\Delta(N)$ has finite index if, and only if $N \leq 5$  (\cite{KLN}, \cite{Wa}).  A subgroup $G \subseteq \Gamma$ is
a \emph{congruence subgroup} if
$\Gamma(N) \subseteq G$ for some $N$; $\rho$ and $(F, \rho)$ are called \emph{modular} if
ker$\rho$ is a congruence subgroup. It follows that $(F, \rho)$ is modular 
if, and only if, the component functions 
$g_j(\tau)$ of $F(\tau)$ are such that $g_j|_k \gamma(\tau)$ has a $q$-expansion  
 of  shape (\ref{gjdef}) for every $\gamma \in \Gamma$. This is precisely the definition of a \emph{classical
modular form} of weight $k$ and level $N$ (we are assuming holomorphy in $\mathfrak{H}$
for convenience).
The case of level $1$ again reduces to the theory discussed in Section \ref{Section_ModForms}. 

\medskip
Because $\Gamma(N)$ has finite index in $\Gamma$ it follows that 
the \emph{image} $\rho(\Gamma)$ is \emph{finite} whenever $\rho$ is modular.
 However, the \emph{converse} is \emph{false}: it may be that
the image $\rho(\Gamma)$ is \emph{finite}, so that ker$\rho$ has finite index in $\Gamma$
and therefore has some finite level, yet it  is not a congruence subgroup. The existence of such subgroups goes back  to Klein and Fricke.
In this case, a vector-valued modular form $(F, \rho)$ will have some finite level $N$
and its component functions have $q$-expansions (\ref{gjdef}),
however not all of them will be classical modular forms in the previous sense. 
This is essentially the theory of modular forms on \emph{noncongruence subgroups}.
Modular forms on noncongruence subgroups, and more generally component functions
of vector valued modular forms, share many properties in common with classical modular forms
and the differences between them can be subtle. It can be difficult to determine whether a
given vector-valued modular form $(F, \rho)$ is modular.
A fundamental problem in this direction is the following: 

\medskip
\noindent
\emph{Conjecture: Let $(F, \rho)$ be a vector-valued modular form of level $N$ and weight $k$,
and suppose that the component functions of $F(\tau)$ are linearly independent\footnote{This condition is harmless in practice, but is necessary to avoid trivial counterexamples, e.g. when $F=0$.} and have 
\emph{rational integers} Fourier coefficients. Then $(F, \rho)$ is modular.}

\medskip
\noindent
We shall see how this fits into VOA theory in the next Section.

\bigskip

\begin{exercise}\label{Exercise_8.1.1} Prove the following: (a) $\Delta(N) \subseteq \Gamma(N) \unlhd \Gamma$,
(b) if $G \subseteq \Gamma$ is a subgroup of finite index then
$\Delta(N) \subseteq G$ for some $N$.
\end{exercise}

\begin{exercise}\label{Exercise_8.1.2} Let $\rho: \Gamma \rightarrow GL_p(\mathbb{C})$ be a representation of level $N$. Show that  $\rho$ is modular if, and only if, $\Gamma(N) \subseteq$ ker$\rho$.
\end{exercise}

\begin{exercise}\label{Exercise_8.1.3}: Let $\tilde{\Gamma}$ be the inhomogeneous modular group
(Exercise \ref{Exercise_3.1.1}) and let $\tilde{\Gamma}(N)$ be the image of $\Gamma(N)$ under the natural projection
$\Gamma \rightarrow \tilde{\Gamma}$. Prove that $\tilde{\Gamma}(N)$ is \emph{torsion-free}
if, and only if, $N \geq 2$.
\end{exercise}

\begin{exercise}\label{Exercise_8.1.4} It is known that $\Gamma$ can be abstractly defined by generators and relations
$\langle x, y \ | \ x^4 =y^6 = x^2y^{-3} = 1 \rangle$.
Use this to prove the following:
(a) $\Gamma/\Gamma' \cong \mathbb{Z}_{12}$, (b) $\Gamma'$ is a congruence subgroup of level $12$.
($\Gamma'$ is the \emph{commutator subgroup} of $\Gamma$.)
\end{exercise}

\begin{exercise}\label{Exercise_8.1.5} Let $V \subseteq \mathfrak{F}_k$ be a finite-dimensional $\Gamma$-submodule
and $(f_1, \hdots, f_p)$ a sequence of functions in $V$ that contains a basis. Prove the existence
of a representation $\rho$ satisfying (\ref{vvformdef}). (Hint: first do the case that $(f_1, \hdots, f_p)$
is a linearly independent set.)
\end{exercise}

\subsection{Examples of Vector-Valued Modular forms}\label{Subsect_ExamplesVVMF}
One can construct a slew of almost holomorphic vector-valued modular forms using \emph{modular linear differential equations} (MLDE) \cite{M}. We briefly explain this. Let $k, n$ be integers with $n$ positive.
The $n$th iterate $D_k^n$ of the differential operator (\ref{ModDer}) is the intertwining map
\begin{eqnarray*}
D_k^n = D_{k+2n-2} \circ \hdots D_{k+2} \circ D_k: \mathcal{F}_k \rightarrow \mathcal{F}_{k+2n}.
\end{eqnarray*}
For justification of the notation, see Exercise \ref{Exercise_3.2.4}. A modular linear differential equation
is a differential equation of the form
\begin{eqnarray}\label{MLDEdef}
(D_k^n + g_2(\tau) D_k^{n-2} + \hdots + g_{2n}(\tau))f = 0 , \ \  g_i(\tau) \in \mathfrak{M}_{2i}.
\end{eqnarray}
Using (\ref{ModDer}) one can write (\ref{MLDEdef}) as an ordinary differential equation with coefficients in the algebra of quasimodular forms $\mathfrak{Q}$. We can also write everything in terms of the variable $q$
(in the interior of the unit disk in the $q$-plane)
\begin{eqnarray}\label{MLDEqdef}
(\theta^n + h_1(q)\theta^{n-1} + \hdots + h_{2n}(q))f=0, \ \ h_i(q) \in \mathfrak{Q}_{2i},
\end{eqnarray}
where we recall that $\theta = q d/dq$. Then one sees that $q=0$ is a regular singular point
(\cite{H}, \cite{I}). By the theory of ODE, the space of solutions is an $n$-dimensional linear space, and because the coefficients are holomorphic in $\mathfrak{H}$, so too are the solutions. One sees that the space of solutions is a $\Gamma$-submodule of $\mathfrak{F}_{k+2n}$, and the theory of
Frobenius-Fuchs (loc cit) shows that the solutions have $q$-expansions which are meromorphic
at $\infty$ in the sense of Subsection \ref{Subsect_VVModformDef}. A disadvantage of this approach is that it is hard to get information about the representation of $\Gamma$ furnished by the space of solutions.

\medskip
We have seen that vector-valued modular forms naturally incorporate
the classical theory of level $N$ modular forms. We complete this Subsection with a discussion of an important class of such forms, namely \emph{theta functions}.
 Let $L$ be an even lattice of rank $d$
with associated positive-definite quadratic form $Q$ (Subsection \ref{Subsect_LatticeVOA}). The theta function of
$L$  is defined by
\begin{eqnarray*}
\theta_L(\tau) = \sum_{\alpha \in L} q^{Q(\alpha)} = \sum_{n \geq 0} |L_n|q^n
\end{eqnarray*}
where $L_n = \{ \alpha \in L \ | \ Q(\alpha) = n\}$ (cf. (\ref{VLpartfuncdef})).
Hecke and Schoeneberg proved (\cite{O}, \cite{Se}, \cite{Sc}) that if $d$ is \emph{even}
then $\theta_L(\tau)$ is a holomorphic modular form of weight $d/2$
and a certain level $N$. A precise description of the level would take us too far afield, but it 
divides twice the exponent of the finite abelian group $L^0/L$ (cf. Exercise \ref{Exercise_7.3.7}). In particular,
suppose that $L$ is \emph{self-dual} in the sense that $L=L^0$. Then the level is $1$, and as we explained this means that $\theta_L(\tau)$ is a holomorphic modular form of weight $d/2$
on the full group $\Gamma$.

\medskip
There are various ways to prove the modularity of $\theta_L(\tau)$. One method that is useful in many 
other contexts is that of \emph{Poisson summation} (\cite{O}, \cite{Se}). 
The approach in (\cite{Sc}) shows that the space 
spanned by the theta functions corresponding to the \emph{cosets}
of $L$ in $L^0$, i.e., the numerators of the expressions on the rhs of (\ref{VLMpartfuncdef}),
is a $\Gamma$-submodule of $\mathfrak{F}_{d/2}$. 
Note that the theta functions of such cosets arise as the numerator in the expression
(\ref{VLMpartfuncdef}) of the character of an ordinary, irreducible module over a lattice VOA.

\medskip
The reader may be wondering about the case when the rank $d$ of $L$ is \emph{odd}. One still has holomorphic theta functions as above, however they are of \emph{half-integral weight} and do not qualify as modular forms as we have defined them. Odd powers of the eta function also have half-integral weight. These and other examples demonstrate the significance of half-integer weight (vector-valued) modular forms to our subject, but there is no time to develop the subject here.

\bigskip

\begin{exercise}\label{Exercise_8.2.1} Use your knowledge of the theory of ODEs to verify the details of the
assertions following (\ref{MLDEdef}) leading to the result that the solution space is a
$\Gamma$-submodule of $\mathfrak{F}_{k+2n}$.
\end{exercise}

\begin{exercise}\label{Exercise_8.2.2} Why is there no term $g_1(\tau)D_k^{n-1}$ in (\ref{MLDEdef})?
\end{exercise}

\begin{exercise}\label{Exercise_8.2.3} For a positive-definite, even lattice $L$ of rank $d$, prove the estimate $|L_n| = O(n^{d/2})$, and deduce that $\theta_L(\tau)$
is \emph{holomorphic} in $\mathfrak{H}$.
\end{exercise}

\begin{exercise}\label{Exercise_8.2.4} Show that $E_8$ is the \emph{only} finite dimensional simple Lie algebra
whose root lattice is  even and self-dual. 
\end{exercise}

\begin{exercise}\label{Exercise_8.2.5} Show that the theta function $\theta_{E_8}(\tau)$ of the $E_8$ root lattice
coincides with the Eisenstein series $Q$ of Section \ref{Section_ModForms}.
\end{exercise}

\begin{exercise}\label{Exercise_8.2.6} Show that the partition functions of a lattice theory $V_L$ and its ordinary, irreducible modules are \emph{classical, almost holomorphic, modular functions of weight zero} of some level $N$.
\end{exercise}

\section{Vertex Operator Algebras and Modular-Invariance}\label{Section_VOAModInvariance}
In this Section we describe some of the main results concerning the
connections between (vector-valued) modular forms and VOAs.
We are concerned here exclusively with regular and $C_2$-cofinite VOAs as discussed in Sections
\ref{Section_Modules} and \ref{Section_ExamplesVOAS}. We recall that $V$ is always assumed to be of CFT-type.

 \subsection{The Regular Case}\label{Subsect_Regular}
It is convenient to assume at the outset that $V$ is a $C_2$-cofinite (but not necessarily rational) VOA of central charge $c$.
By Theorem \ref{thmC2VOA} there are only finitely many inequivalent, ordinary, irreducible $V$-modules, and we denote them $V=M^1, M^2, \hdots, M^r$. Let the conformal weight of $M^i$ be $h^i$ (cf. (\ref{Vmodpartfunc}) and  attendant discussion), and let $Z_i$ be the character of $M^i$ (\ref{chiMdef}).

\medskip
The first basic fact is that  \emph{$1$-point functions are
holomorphic in $\mathfrak{H}$}. For example, it follows from this and Theorems \ref{TheoremZhu} and
\ref{Proposition_1pt} that the $2$-point functions $F_V((u_1, z_1), (u_2, z_2))$ are elliptic functions.
 There are two approaches to the holomorphy of $1$-point functions. The first
 (\cite{Z}) is to find a modular linear differential equation (\ref{MLDEqdef}) satisfied
by $f=Z_i(v, q)$. In this case, because the coefficients of the MLDE are holomorphic in $\mathfrak{H}$, 
then so are the solutions
(cf. Exercise \ref{Exercise_8.2.1}). The second approach (\cite{GN}) uses the PBW-type bases that we already mentioned in Subsection \ref{Subsect_C2cofinite}.

\medskip
We now take $V$ to be regular. The main properties \emph{vis-\`{a}-vis} modular-invariance are as follows:
\begin{theorem}\label{thmregmodinv} Let the notation be as above, and assume that $V$ is \emph{regular}. Then the following hold.\\
\noindent
(a) The central charge $c$ and conformal weights $h^i$ are \emph{rational} numbers.\\
(b) There is a representation $\rho: \tilde{\Gamma} \rightarrow GL_r(\mathbb{C})$ of the inhomogeneous modular group (cf. Exercise \ref{Exercise_3.1.1}) with the following property:
if $v \in V$ has $L[0]$-weight $k$ and we set
$F_v = (Z_1(v), \hdots, Z_r(v))$, then $(F_v, \rho)$ is an 
\emph{almost holomorphic vector-valued modular form of weight $k$ and finite level $N$}.
\end{theorem}

 We have already discussed the holomorphy of $Z_i(v)$.
 The heart of the matter - that there is $\rho$ such that $(F_v, \rho)$ is a vector-valued modular form of weight $k$ - 
 is more difficult. It ultimately depends on the complete reducibility of admissible $V$-modules into ordinary irreducible $V$-modules.
 See \cite{Z}, \cite{DLM4} for details. The argument  shows that the representation $\rho$ is \emph{independent of the state $v$}. Once the vector-valued modular form is available, one uses the theory of ODEs with regular singular points (\cite{MA}) to show that (a) holds. The argument, which is arithmetic in nature, makes use of the fact that if $v$ is taken to be the vacuum vector then the component functions $Z_i(\mathbf{1})$ of $F_{\mathbf{1}}$ are just the partition functions of the ordinary irreducible modules over $V$, and as such have \emph{integral} Fourier coefficients. 
 Also, because $F_{\mathbf{1}}$ has weight zero (because $\mathbf{1} \in V_{[0]}$), ker$\rho$ contains $\pm I_2$ and so $\rho$ descends to a representation of $\tilde{\Gamma}$. The rationality of conformal weights and central charge implies that
 $(F_v, \rho)$ has finite level $N$ (e.g., one can take $N$ to be the gcd of the denominators of
 the rational numbers $h_i - c/24$). There is a basic open problem here, namely:
 
 \medskip
 \noindent
 \emph{Modularity Conjecture: In the context of Theorem \ref{thmregmodinv}, $(F_v, \rho)$ is \emph{modular}}.
 
 \medskip
 This is an article of faith in the physics literature. There are compelling arguments
 (e.g., \cite{Ba1}, \cite{Ba2}, \cite{FMS}) which, however, are not (yet) mathematically rigorous.
  Note that this Conjecture follows from the conjectured modularity of vector-valued modular forms
  of level $N$ with integral Fourier coefficients stated at the end of Subsection \ref{Subsect_VVModformDef}. There are 
  other avenues via which the modularity of $(F_v, \rho)$ might be established, in particular using the theory of tensor products of modules over a VOA and tensor categories (cf. \cite{HL}).
  
  \medskip
  It hardly needs to be said that all known regular VOAs satisfy the Modularity Conjecture. The case
  of lattice theories follows from Exercise \ref{Exercise_8.2.6}. The case of
  WZW models was studied prior to the advent of VOA theory using Lie theory (cf. \cite{KP}, \cite{K2}).
  A discussion of this case as well as that of the simple Virasoro VOAs $L_c$ in the discrete series may be found in \cite{FMS}.

  \subsection{The $C_2$-Cofinite Case}\label{Subsect_C2}
 One desires an analog of Theorem \ref{thmregmodinv} for the more general case of
 $C_2$-cofinite VOAs, but any generalization must deal with the fact that the span of the
 partition functions $Z_i(\mathbf{1})$ of the ordinary irreducible modules is generally \emph{not} a $\Gamma$-module unless $V$ is a regular VOA. Miyamoto's solution \cite{My1} (see also \cite{Fl})
 involves generalized or \emph{pseudo trace functions}. The idea is to utilize the admissible $V$-modules $L_n(X)$ constructed
 from a finite-dimensional module $X$ over the algebra $A_n(V)$ (\ref{Lfuncdef}). $C_2$-cofiniteness implies that 
 $A_n(V)$ is finite-dimensional (Theorem \ref{thmC2VOA}), and this leads to the fact that each of the homogeneous pieces $L_n(X)_m$ are also finite-dimensional. Because $L_n(X)$ is admissible then
 the zero mode $o(\omega)=L_0$ of the conformal vector operates on these homogeneous pieces
 (\ref{modeaction2}). However, in the present context $L_0$ may not be the degree operator, indeed $L_0$ \emph{may not be a semisimple operator}. 
 
 \medskip
 We decompose $L_n(X)_m$ into a direct sum of
 Jordan blocks for the action of $L_0$. On such a block $B$ there is an $L_0$-eigenvector with
 eigenvalue $m+\lambda, \ \lambda \in \mathbb{C}, \ L_0 - (m + \lambda)I$ is
 \emph{nilpotent}, and the exponential operator 
 \begin{eqnarray}\label{nilpqtrace}
q^{L_0} = q^{m+\lambda} \sum_{t \geq 0} \frac{(2\pi i \tau(L_0-m-\lambda))^t}{t!}
\end{eqnarray}
on $B$ reduces to a \emph{finite} sum. If $X$ is indecomposable, $\lambda$ is determined by
the action of $\omega$, which (when regarded as an element of $A_n(V)$) turns out to be a central element and thus acts on $X$ as a scalar. One can piece together the exponentials
(\ref{nilpqtrace}) and incorporate zero modes $o(v)$ of other states as before. However, the details are not quite straightforward as one needs pseudotraces (\cite{My1}), which is a type of symmetric function on $A_n(V)$ which replaces the usual trace. 
  
  \medskip
  The upshot of the analysis sketched above is this: we can define\footnote{$\phi$ denotes `pseudo'. } (pseudo) trace functions
  $\mathrm{Tr}^{\phi}_{L_n(X)} o(v)q^{L_0-c/24}$. Once these gadgets are introduced, one can use the arguments in the regular case described in the previous Subsection together with additional arguments (to account for the failure of $A_n(V)$ to be semisimple)
   to show that for each $n$ and for $v \in V_{[k]}$, the pseudo trace functions define a 
   (finite-dimensional) almost holomorphic vector-valued modular form of weight $k$. Alternatively, they span a finite-dimensional $\Gamma$-submodule of $\mathfrak{F}_k$ (notation as in Subsection \ref{Subsect_VVModformDef}). In particular, the pseudo characters 
  $\mathrm{Tr}^{\phi}_{L_n(X)} q^{L_0-c/24}$ are seen to be linear combinations of
  characters of ordinary, irreducible $V$-modules with coefficients in $\mathbb{C}[\tau]$.
  That is, they are polynomial $q$-expansions in the sense of Subsection \ref{Subsect_VVModformDef}. This is, of course,
  fully consistent with Theorem \ref{thmlogqexp}.  Furthermore, one finds as in the regular case
  that the central charge and conformal weights of the ordinary, irreducible $V$-modules again lie in 
  $\mathbb{Q}$.

 \medskip
 It would take as too far afield to try to describe any VOAs for which the pseudo trace functions actually involve log terms. Such theories are, naturally, called \emph{logarithmic field theories} in the physics literature. For some examples, see e.g.,  \cite{GK}, \cite{A} and references therein.
  
  \bigskip
  
\begin{exercise}\label{Exercise_9.2.1} Prove that the (image of) the conformal vector $\omega$
  is a central element of $A_n(V)$ (cf. Exercises \ref{Exercise_6.1.6}, \ref{Exercise_6.1.7}).
 \end{exercise}
 
 \subsection{The Holomorphic Case}\label{Subsect_Holomorphic}
We call a simple, regular VOA $V$ \emph{holomorphic} if it has a
 \emph{unique} irreducible module, namely the adjoint module $V$. It seems likely that a simple 
 $VOA$ with a unique ordinary irreducible module is necessarily regular, and therefore holomorphic, but this appears to be unknown. Be that as it may, in the case of holomorphic
 VOAs Theorem \ref{thmregmodinv} can be refined, and in particular the Modularity Conjecture
 of Subsection \ref{Subsect_Regular} holds in this case. This is because if a vector-valued modular form 
 of weight $k$ has a single component $f(\tau)$ then it affords a $1$-dimensional representation of $\Gamma$ and
 so there is a character $\alpha: \Gamma \rightarrow \mathbb{C}^*$ such that
 \begin{eqnarray}\label{holmodinv}
f|_k \gamma( \tau) = \alpha(\gamma) f(\tau), \ \ \gamma \in \Gamma.
\end{eqnarray}
 Since $\Gamma'$ is a congruence subgroup of level $12$ (Exercise \ref{Exercise_8.1.4}) it follows that $f(\tau)$ is
 a classical modular form of level dividing $12$. Thanks to Theorem \ref{thmregmodinv} all of this applies  with $f = Z_V(v, q)$, indeed a bit more is true in this case: the group of characters of $\Gamma$ is cyclic of order $12$ (Exercise \ref{Exercise_8.1.4}) hence that of $\tilde{\Gamma}$ is cyclic of order $6$; and one can argue (cf. Exercise \ref{Exercise_9.3.1}) that
 $S \in$ ker$\alpha$, so that in fact $\alpha$ has order dividing $3$ and each $\alpha(\gamma)$
 in (\ref{holmodinv}) is a cube root of unity.
 We thus arrive at
  \begin{theorem}\label{holVOA} Suppose that $V$ is a holomorphic VOA of central charge $c$. Then the following hold:\\
(a) If $v \in V_{[k]}$ then $Z_V(v, \tau)$ is an almost holomorphic modular form of  weight $k$ and level $1$ or $3$.\\
(b) $c$ is an \emph{integer divisible by $8$}. It is divisible by $24$ if, and only if, $Z_V(v, q)$ has level $1$.
 \end{theorem}
 
 \medskip
 Lattice theories provide a large number of holomorphic VOAs. From Theorem \ref{thmVLVOA}
 it is immediate that $V_L$ is holomorphic if, and only if, $L=L^0$ is self-dual. 
 The partition function is
$\theta_L(\tau)/\eta^{c}(\tau)$
where $c$ is the rank of $L$ (\ref{VLpartfuncdef}), and in this case the modularity of
the partition function follows directly from comments in Subsection \ref{Subsect_ExamplesVVMF}. 
 
 \medskip
 We also mention that the modules over a tensor product $U \otimes V$ of VOAs 
 (Exercise \ref{Exercise_2.6.7}) are just the tensor products $M \otimes N$ of modules $M$ over $U$ and $N$ over $V$
( \cite{FHL}). In particular, if $U, V$ are holomorphic then so too is $U \otimes V$.
 
 \bigskip

\begin{exercise}\label{Exercise_9.3.1} Let $V$ be a holomorphic VOA, and let $\alpha$ be the character of $\Gamma$ satisfying  $(*) \  Z_V(\mathbf{1})|_0 \gamma (\tau) = \alpha(\gamma) Z_V(\mathbf{1})$.  
 Prove that $\alpha(S) = 1$. (Hint: take $\gamma = S$ and evaluate $(*)$ at $\tau = i$.) Using this, give the details of the proofs of (a) and (b) in Theorem \ref{holVOA}. 
\end{exercise}

\begin{exercise}\label{Exercise_9.3.2} Let $V$ be a holomorphic VOA of central charge $c$, and let $v \in V_{[k]}$. Prove that
 $Z_V(v, \tau) = g(\tau)/\eta^c(\tau)$ where $g(\tau)$ is an almost holomorphic
 modular form on $\Gamma$ of weight $k+c/2$.
\end{exercise}

 \subsection{Applications of Modular-Invariance}\label{Subsect_Applications}
 Theorem \ref{thmregmodinv} places strong conditions on the $1$-point trace functions of a regular VOA,  and in particular on the partition function. If $V$ is a holomorphic VOA the conditions are even stronger. In this Subsection we give a few illustrations of how modular-invariance can be used to study the structure of holomorphic VOAs. 

\medskip
By Exercise \ref{Exercise_9.3.2}, $Z_V(\mathbf{1}) = g(\tau)/\eta^c(\tau)$ where $g(\tau) = 1 + \hdots \in \mathfrak{M}_{c/2}$
 is a holomorphic modular form on $\Gamma$ of weight $c/2$. There are no (nonzero)
 such forms of negative weight, so we have $c\geq 0$. If $c=0$ then $g(\tau)=1=Z_V(\tau)$,
 corresponding to the $1$-dimensional VOA $\mathbb{C}\mathbf{1}$ (cf. Exercise \ref{Exercise_2.6.5}) which is indeed holomorphic. 
 
 \medskip
 Since $8|c$, the next two cases are $c=8, 16$, when $g(\tau)$ has weight
 $4$ and $8$ respectively. Because of the structure of the algebra $\mathfrak{M}$ of modular forms
 on $\Gamma$ (Theorem \ref{thm4.1} and (\ref{HS1})) there is only one choice for
 $g(\tau)$ in these cases, namely $g(\tau) = Q$ or $Q^2$, so the partition function is
 \emph{uniquely determined} as
 $Z_V(\mathbf{1}) = Q/\eta^8(\tau)$ or $Q^2/\eta^{16}(\tau) = \left(Q/\eta^8(\tau) \right)^2$ (Exercise \ref{Exercise_8.2.5} is relevant here).  We have already seen holomorphic VOAs with these partition functions in Subsection \ref{Subsect_Holomorphic}, namely the lattice theories $V_{E_8}$ and $V_{E_8 \perp E_8} \equiv V_{E_8}^{\otimes 2}$ ($E_8$ refers to the root lattice of type $E_8$).  In fact, there is a second even, self-dual lattice 
 $L_2$ of rank $16$ not isometric to $E_8 \perp E_8$ and we obtain in this way a second holomorphic
 VOA $V_{L_2}$. It turns out that these are the \emph{only} holomorphic VOAs (up to isomorphism) with $c=8$ or $16$. This result requires additional techniques based on applications
 of the recursion in Theorem \ref{Proposition_1pt} and analytic properties of vector-valued modular forms (\cite{DM2}, \cite{DM3}). To summarize:
 \begin{theorem} Suppose that $V$ is a holomorphic VOA of central charge $c \leq 16$. Then one of the following holds:\\
 (a) $c=0$ and $V = \mathbb{C}\mathbf{1}$.\\
 (b) $c=8$ and $V = V_{E_8}$ is the $E_8$-lattice theory. \\
 (c) $c=16$ and $V=V_{E_8 \perp E_8}$ or $V_{L_2}$ is a lattice theory.
 \end{theorem}
 
\medskip
We now consider holomorphic VOAs $V$ of central charge $c=24$. In some ways, this is the most interesting case. If $c \geq 32$ the number of isometry classes of even, self-dual
lattices of rank $c$ is very large  (see \cite{Se} for further comments), so there are a correspondingly large number of isomorphism classes of holomorphic VOAs. For rank $24$ there are just $24$ isometry classes of even, self-dual lattices (cf. \cite{CS}, \cite{Se}), so one might hope that there are not too many holomorphic VOAs with $c=24$. In fact, Schellekens has conjectured that there are just $71$ such theories (\cite{Sch}).
Now
$Z_V(\mathbf{1}) = q^{-1}+ \hdots$ is an almost holomorphic modular function of weight zero and level $1$
by Theorem \ref{holVOA}. As such it is a polynomial in the modular function 
$j(\tau) = q^{-1}+744+ \hdots$
(cf. (\ref{jdef}) and Exercise \ref{Exercise_3.3.5}). So there is an integer $d$ such that
\begin{eqnarray*}
Z_V(\mathbf{1}) = j(\tau) + (d-744) = q^{-1}+d+196884q+ \hdots
\end{eqnarray*}
and the partition function is \emph{determined uniquely} by $d$. Obviously $d=\dim V_1$, so it is a nonnegative integer, but one cannot say more 
about $d$ on the basis of modular-invariance alone because $j(\tau)+c'$ is a modular function
for \emph{any} constant $c'$. It can in fact be proved that there are only finitely many choices  of $d$ that correspond to possible holomorphic VOAs\footnote{No more than a few hundred.}. The arguments use Lie algebra theory, starting with the Lie algebra structure on $V_1$ (Exercise \ref{Exercise_9.4.1}) as well as modular forms (\cite{DM1}-\cite{DM3}, \cite{Sch}). Of the $71$ conjectured holomorphic $c=24$ VOAs, it seems that
only $39$ are known to exist. Beyond the $24$ lattice theories, the other $15$ are constructed as so-called
$\mathbb{Z}_2$-orbifold models of lattice theories \cite{DGM}. The first construction of this type
\cite{FLM}  leads to  the famous \emph{Moonshine Module},  about which we will shortly say a bit more.
It is a major problem to decide whether the others also exist, and to develop construction techniques when they do.

 \medskip
As a final example, we mention some recent work of E. Witten \cite{Wi} where certain holomorphic vertex operator algebras $V^{(k)}$ are posited to exist which are related, via the \emph{AdS-CFT correspondence},  to  phenomena concerning gravity with a negative cosmological constant. 
$V^{(k)}$ has central charge $c_k = 24k, k = 1, 2, \hdots$ and a \emph{minimal structure} 
compatible with the requirements of modular-invariance imposed by Theorem \ref{holVOA}. 
To explain what this is supposed to mean, recall (cf. Theorem \ref{thmratVir}) that 
$\mbox{Vir}_{c_k}=L_{c_k}$ is simple, and the $L_{c_k}$-submodule of
 $V^{(k)}$ generated by $\mathbf{1}$ is a graded subspace $U$ naturally identified as the Fock space for
 $L_{c_k}$. By (\ref{vircqdim}), the graded dimension of $U$
is 
\begin{eqnarray*}
q^{-k}\prod_{n \geq 2}(1-q^n)^{-1} = q^{-k}\sum_{n=0}^k d_nq^n + O(q)
\end{eqnarray*}
for integers $d_0, \hdots d_k$. The posited minimal structure of $V^{(k)}$ means  that 
the partition function of $V^{(k)}$ also satisfies
\begin{eqnarray*}
Z_{V^{(k)}}(\mathbf{1}) = q^{-k}\sum_{n=0}^kd_nq^n + O(q).
\end{eqnarray*}
In other words, the first $k+1$ graded subspaces $V^{(k)}_n (0 \leq n \leq k)$ of
$V^{(k)}$ coincide with the corresponding graded pieces of $U$, so that they are as small as they can be. We know that $Z_{V^{(k)}}(\mathbf{1})$ is a monic polynomial $\Phi_k(j)$
of degree $k$ in $j(\tau)$, and it is clear that $\Phi_k$ is \emph{uniquely determined}
by $d_0, \hdots, d_k$, and hence by $k$. 

\medskip
As in the case of the `missing' holomorphic $c=24$ theories, the main question here for the VOA theorist is whether $V^{(k)}$ exists or not.  The answer is unknown for any $k$ with the notable exception of the
\emph{Moonshine module} $V^{\natural}$  (\cite{B}, \cite{FLM}, \cite{DGM}, \cite{My2})  corresponding to $k=1$. In this case the graded dimension of $U$ is $q^{-1} + O(q)$,
the partition function of $V^{\natural}$ is
\begin{eqnarray*}
Z_{V^{\natural}}(q) = j(q) - 744 = q^{-1} + 0 + 196884q+ \hdots,
\end{eqnarray*}
and the minimal structure is  reflected in the vanishing of the constant term. In this case the Lie algebra
   structure on the weight $1$ subspace is absent, and one must exploit instead the 
   \emph{Griess algebra}, i.e., the commutative
   algebra structure on $V^{\natural}_2$ (cf. Exercise \ref{Exercise_9.4.3}).
   
   \medskip
   One of the main features of the Moonshine Module is its automorphism group, which is the Monster sporadic simple group (\cite{FLM}, \cite{G1}, \cite{G2}). In order to develop this aspect of $V^{\natural}$ as well as the $\mathbb{Z}_2$-orbifold construction that we mentioned above and other features of VOAs, it would be necessary to develop the theory of \emph{automorphism groups} of VOAs. This will have to wait for another day. A brief description of some of the connections between automorphism groups and generalized modular forms can be found in \cite{KM}.
   
\bigskip

\begin{exercise}\label{Exercise_9.4.1} Let $V$ be a VOA of CFT-type. Prove the following:\\
(a) The product
$[a, b] = a_0b$ equips $V_1$ with the structure of a \emph{Lie algebra}.\\
(b) $\langle a, b \rangle = a_1b$ defines a symmetric, invariant, bilinear form on $V_1$.
\end{exercise}

\begin{exercise}\label{Exercise_9.4.2} Prove that the Fourier coefficients of the $q$-expansion of  $\Phi_k(j)$
are nonnegative integers (a necessary condition for the existence of $V^{(k)}$).
\end{exercise}

\begin{exercise}\label{Exercise_9.4.3} Show that the product $a_1b \ (a, b \in V^{\natural}_2)$ equips the weight $2$
subspace of $V^{\natural}$ with the structure of a \emph{commutative, nonassociative}
algebra.
\end{exercise}

\newpage

\part{Two Current Research Areas}\label{PartIII}

\section{Some Preliminaries}
\label{Section_Preliminary}

\subsection{VOAs and Rational Matrix Elements}\label{Subsect_MatrixElements}

As noted in Section \ref{Subsect_VOAaxioms} there are a number of equivalent sets of axioms for VOA theory. 
Here we discuss one of these equivalent approaches wherein the properties of a VOA are expressed in terms of the properties of matrix elements which turn out to be rational functions of the formal vertex operator parameters. In many ways, this is the closest approach to CFT (e.g. \cite{FMS}) in that the formal parameters can
be taken to be complex numbers with the matrix elements considered as
rational functions on the Riemann sphere.
\bigskip

We begin by defining matrix elements. In order to simplify the discussion,
we always assume that the VOA is of CFT-type (\ref{cfttypedef}). This condition is satisfied in all examples we consider.   We define the \emph{restricted dual space} of $V$ by \cite{FHL} 
\begin{equation}
V^{\prime }=\oplus _{n\geq 0}V_{n}^{\ast },  \label{eq: Vdual}
\end{equation}%
where $V_{n}^{\ast }$ is the dual space of linear functionals on the finite
dimensional space $V_{n}$. Let $\langle \ ,\rangle _{d}$ denote the
canonical pairing between $V^{\prime }$ and $V$. Define \emph{matrix
elements} for $a^{\prime }\in V^{\prime }$, $b\in V$ and vertex operators 
$Y(u^{1},z_{1}),\ldots Y(u^{n},z_{n})$ by 
\begin{equation}
\langle a^{\prime },Y(u^{1},z_{1})\ldots Y(u^{n},z_{n})b\rangle _{d}.
\label{eq: Fmatrix}
\end{equation}%
In particular, choosing $b=\mathbf{1}$ and $a^{\prime }=\mathbf{1}^{\prime }$
we obtain the (genus zero) $n$-point correlation function 
\begin{equation}
F_{V}^{(0)}((u^{1},z_{1}),\ldots ,(u^{n},z_{n}))=\langle \mathbf{1}^{\prime
},Y(u^{1},z_{1})\ldots Y(u^{n},z_{n})\mathbf{1}\rangle _{d}.
\label{eq: genus0npt}
\end{equation}

\medskip
One can show in general that every matrix element is a homogeneous rational
function of $z_{1},\ldots ,z_{n}$ \cite{FHL}, \cite{DGM}. Thus the formal
parameters of VOA theory can be replaced by complex parameters on
(appropriate subdomains of) the genus zero Riemann sphere $\mathbb{CP}^{1}$. 
We illustrate this by considering matrix elements containing one or two
vertex operators. Recall from (\ref{gradedmodes}) that for $u\in V_{n}$ 
\begin{equation}
u_{k}:V_{m}\rightarrow V_{m+n-k-1}.  \label{eq: vnhom}
\end{equation}%
Hence it follows that for $a^{\prime }\in V_{m^{\prime }}^{\prime }$, $b\in
V_{m}$ and $u\in V_{n}$ we obtain a monomial 
\begin{equation}
\langle a^{\prime },Y(u,z)b\rangle _{d}=
C_{a^{\prime }b}^{u}z^{m^{\prime}-m-n},  \label{eq: 1ptmatrix}
\end{equation}%
where $C_{a^{\prime }b}^{u}=\langle a^{\prime },u_{m+n-m^{\prime}-1}b\rangle _{d}$.

\bigskip

We next consider the matrix element of two vertex operators to find
(recalling convention (\ref{xypower})):

\begin{theorem}
\label{theorem: RationalYY}Let $a^{\prime }\in V_{m^{\prime }}^{\prime }$, 
$b\in V_{m}$, $u^{1}\in V_{n_{1}}$ and $u^{2}\in V_{n_{2}}$. Then
\begin{eqnarray}
\langle a^{\prime },Y(u^{1},z_{1})Y(u^{2},z_{2})b\rangle _{d} &
=&\frac{f(z_{1},z_{2})}{z_{1}^{m+n_{1}}z_{2}^{m+n_{2}}(z_{1}-z_{2})^{n_{1}+n_{2}}},
\label{eq: Fu1u2} \\
\langle a^{\prime },Y(u^{2},z_{2})Y(u^{1},z_{1})b\rangle _{d} &
=&\frac{f(z_{1},z_{2})}{z_{1}^{m+n_{1}}z_{2}^{m+n_{2}}(-z_{2}+z_{1})^{n_{1}+n_{2}}},
\label{eq: Fu2_u1}
\end{eqnarray}%
where $f(z_{1},z_{2})$ is a homogeneous polynomial of degree 
$m+m^{\prime}+n_{1}+n_{2}$.
\end{theorem}

\begin{remark}
\label{Remark ComplexDomain}The matrix elements (\ref{eq: Fu1u2}), (\ref{eq:
Fu2_u1}) are thus determined by a unique homogeneous rational function which
can be evaluated on $\mathbb{CP}^{1}$ in the domains $\left\vert
z_{1}\right\vert >\left\vert z_{2}\right\vert $ and $\left\vert
z_{2}\right\vert >\left\vert z_{1}\right\vert $ respectively.
\end{remark}

\noindent\textbf{Proof.} Consider 
\begin{equation*}
\langle a^{\prime },Y(u^{1},z_{1})Y(u^{2},z_{2})b\rangle _{d}=\sum_{k\geq
0}\sum_{c\in V_{k}}\langle a^{\prime },Y(u^{1},z_{1})c\rangle _{d}\langle
c^{\prime },Y(u^{2},z_{2})b\rangle _{d},
\end{equation*}%
where $c$ ranges over any basis of $V_{k}$ and $c^{\prime }\in V_{k}^{\ast }$
is dual to $c$. From (\ref{eq: 1ptmatrix}) it follows that 
\begin{equation*}
\langle a^{\prime },Y(u^{1},z_{1})Y(u^{2},z_{2})b\rangle _{d}=\frac{%
z_{1}^{m^{\prime }-n_{1}}}{z_{2}^{m+n_{2}}}G(\frac{z_{2}}{z_{1}}),
\end{equation*}%
for infinite series 
\begin{equation*}
G(x)=\sum_{k\geq 0}\sum_{c\in V_{k}}
C_{a^{\prime }c}^{u^{1}}C_{c^{\prime}b}^{u^{2}}x^{k}.
\end{equation*}%
Hence the matrix element is homogeneous of degree $m^{\prime
}-m-n_{1}-n_{2}$. Similarly 
\begin{equation*}
\langle a^{\prime },Y(u^{2},z_{2})Y(u^{1},z_{1})b\rangle _{d}=
\frac{z_{2}^{m^{\prime }-n_{2}}}{z_{1}^{m+n_{1}}}H(\frac{z_{1}}{z_{2}}),
\end{equation*}%
for infinite series 
\begin{equation*}
H(y)=\sum_{k\geq 0}\sum_{c\in V_{k}}
C_{a^{\prime }c}^{u^{2}}C_{c^{\prime}b}^{u^{1}}y^{k}.
\end{equation*}%
But $Y(u^{2},z_{2})$ and $Y(u^{1},z_{1})$ are local of order at most 
$n_{1}+n_{2}$ (cf. Exercise \ref{Exercise_truncweight}) and hence
\begin{equation}
\frac{(z_{1}-z_{2})^{n_{1}+n_{2}}}{z_{1}^{m+n_{1}}z_{2}^{m+n_{2}}}%
z_{1}^{m^{\prime }+m}G(\frac{z_{2}}{z_{1}})=\frac{(z_{1}-z_{2})^{n_{1}+n_{2}}%
}{z_{1}^{m+n_{1}}z_{2}^{m+n_{2}}}z_{2}^{m^{\prime }+m}H(\frac{z_{1}}{z_{2}}).
\label{eq: FG}
\end{equation}%
Then it follows that 
\begin{equation*}
f(z_{1},z_{2})=z_{1}^{m^{\prime }+m}(z_{1}-z_{2})^{n_{1}+n_{2}}G(\frac{z_{2}%
}{z_{1}})=z_{2}^{m^{\prime }+m}(z_{1}-z_{2})^{n_{1}+n_{2}}H(\frac{z_{1}}{%
z_{2}}),
\end{equation*}%
is a homogeneous polynomial of degree $m+m^{\prime }+n_{1}+n_{2}$. 

\medskip

Properties (\ref{eq: Fu1u2}) and (\ref{eq: Fu2_u1}) are equivalent to
locality of $Y(u^{1},z_{1})$ and $Y(u^{2},z_{2})$ so that the axioms of a
VOA can be alternatively formulated in terms of rational matrix elements 
\cite{DGM}, \cite{FHL}. Theorem \ref{theorem: RationalYY} can also be
generalized for all matrix elements. Furthermore, using the vertex
commutator property (cf. Exercise \ref{Exercise_2.3.7}) one can also derive a recursive
relationship in terms of rational functions between matrix elements for $n$
vertex operators and $n-1$ vertex operators that is the genus zero version
of Zhu's recursion formula I of Theorem \ref{TheoremZhu} \cite{Z}.

\begin{exercise}

\label{Exercise Cabu} Prove (\ref{eq: 1ptmatrix}).
\end{exercise}

\begin{exercise}
\label{Exercise fhomog} Show (\ref{eq: FG}) implies that $f(z_{1},z_{2})$ is
a polynomial.
\end{exercise}

\subsection{Genus Zero Heisenberg Correlation Functions}
\label{Subsect_HeisenCorrFun}

We illustrate these structures by considering the example of the rank one
Heisenberg VOA $M_0$ generated by a weight one vector $a$. Let 
\begin{equation}
G_{n}^{(0)}(z_{1},\ldots ,z_{n})\equiv F_{M_0}^{(0)}((a,z_{1}),\ldots
,(a,z_{n})).  \label{eq: Fan}
\end{equation}%
denote the $n$-point correlation function for $n$ Heisenberg vectors. This
must be a symmetric rational function in $z_{i}$ with poles of order two at $%
z_{i}=z_{j}$ for all $i\neq j$ from locality. We now determine its exact
form. Since $a_{0}\mathbf{1}=0$ it follows that $G_{1}^{(0)}(z_{1})=0$. The
2-point function is%
\begin{equation*}
G_{2}^{(0)}(z_{1},z_{2})=\sum_{m\geq 0}z_{1}^{-m-1}\langle \mathbf{1}%
^{\prime },a_{m}Y(a,z_{2})\mathbf{1}\rangle _{d},
\end{equation*}%
where (\ref{eq: vnhom}) implies that there is no contribution for $m<0$.
Commuting $a_{m}$ we find 
\begin{equation*}
G_{2}^{(0)}(z_{1},z_{2})=\sum_{m\geq 0}z_{1}^{-m-1}\langle \mathbf{1}%
^{\prime },[a_{m},Y(a,z_{2})]\mathbf{1}\rangle _{d},
\end{equation*}%
using $a_{m}\mathbf{1}$ $=0$ for $m\geq 0$. But the Heisenberg commutation
relations imply 
\begin{equation*}
\lbrack a_{m},Y(a,z_{2})]=mz_{2}^{m-1},
\end{equation*}%
so that 
\begin{equation}
G_{2}^{(0)}(z_{1},z_{2})=\sum_{m\geq 0}mz_{1}^{-m-1}z_{2}^{m-1}=\frac{1}{%
(z_{1}-z_{2})^{2}}.  \label{eq: F2Heisenberg}
\end{equation}

\medskip
The general $n$-point function is similarly given by 
\begin{equation*}
G_{n}^{(0)}(z_{1},\ldots ,z_{n})=
\sum_{m\geq 0}z_{1}^{-m-1}\sum_{i=2}^{n}\langle \mathbf{1}^{\prime },Y(a,z_{2})\ldots
\lbrack a_{m},Y(a,z_{i})]\ldots Y(a,z_{n})\mathbf{1}\rangle _{d},
\end{equation*}%
leading to a recursive identity 
\begin{equation}
G_{n}^{(0)}(z_{1},\ldots ,z_{n})=\sum_{i=2}^{n}\frac{1}{(z_{1}-z_{i})^{2}}%
G_{n-2}^{(0)}(z_{2},\ldots ,\hat{z}_{i}\ldots ,z_{n}),  \label{eq: Gniter}
\end{equation}%
where $\hat{z}_{i}$ is deleted. Thus we may recursively solve to find 
$G_{n}^{(0)}=0$ for $n$ odd whereas for $n$ even, $G_{n}^{(0)}$ is expressed
as multiples of rational terms of the form $\frac{1}{(z_{i}-z_{j})^{2}}$ for
all possible pairings $z_{i},z_{j}$. This can be equivalently described in
terms of the subset, denoted by $F(\Phi )$, of the permutations of the label
set $\Phi =\{1,\ldots n\}$ consisting of \emph{fixed-point-free involutions}%
. Thus a typical element $\varphi \in F(\Phi )$ is given by $\varphi =\ldots
(ij)\ldots $, a product of $n/2$ disjoint cycles. We then find (\ref{eq:
Gniter}) implies

\begin{theorem}
\label{theorem: HeisenGenFun}$G_{n}^{(0)}$ vanishes for $n$ odd whereas for $%
n$ even%
\begin{equation}
G_{n}^{(0)}(z_{1},\ldots ,z_{n})=\sum_{\varphi \in F(\Phi )}\prod \frac{1}{%
(z_{i}-z_{j})^{2}},  \label{eq: FHeisen}
\end{equation}%
where the product ranges over all the cycles of $\varphi =\ldots (ij)\ldots $%
. 
\end{theorem}

\begin{remark}
\label{Remark Gen G}Using associativity one can show that $%
G_{n}^{(0)}(z_{1},\ldots ,z_{n})$ is in fact a generating function for all
matrix elements of the Heisenberg VOA.
\end{remark}

\begin{exercise}
\label{Exercise F(Phi)}Show that $\left\vert F(\Phi )\right\vert =(n-1)!!=(n-1).(n-3).(n-5)\ldots$.
\end{exercise}

\begin{exercise}
\label{Example G4}For $n=4$ show that 
$F(\Phi)=\{(12)(34),(13)(24),(14)(23)\}$ and $G_{4}^{(0)}(z_{1},z_{2},z_{3},z_{4})$
is given by%
\begin{equation*}
\frac{1}{(z_{1}-z_{2})^{2}(z_{3}-z_{4})^{2}}+\frac{1}{%
(z_{1}-z_{3})^{2}(z_{2}-z_{4})^{2}}+\frac{1}{%
(z_{1}-z_{4})^{2}(z_{2}-z_{3})^{2}}.
\end{equation*}
\end{exercise}

\subsection{Adjoint Vertex Operators}
\label{Subsect_Adjoint}

The Virasoro subalgebra $\{L_{-1},L_{0},L_{1}\}$ generates a natural action
on vertex operators associated with $SL(2,\mathbb{C})$ M\"{o}bius
transformations  on $z$ (cf. \cite{B}, \cite{DGM}, \cite{FHL}, \cite{K1} and Exercise \ref{Exercise_2.5.1}). Thus under the translation $z\mapsto z+\lambda $ generated by $L_{-1}$ we have (cf. Exercise \ref{Exercise_translation}) 
\begin{equation}
q_{\lambda }^{L_{-1}}Y(u,z)q_{\lambda }^{-L_{-1}}=Y(u,z+\lambda ).
\label{eq: Y_T}
\end{equation}%
Under  $z\mapsto q_{\lambda }z$ generated by $L_{0}$ we have
(cf. Exercise \ref{Exercise_4.1.2}) 
\begin{equation}
q_{\lambda }^{L_{0}}Y(u,z)q_{\lambda }^{-L_{0}}=Y(q_{\lambda
}^{L_{0}}u,q_{\lambda }z).  \label{eq: Y_D}
\end{equation}%
Finally, under the transformation $z\mapsto \frac{z}{1-\lambda z}$ generated
by $L_{1}$ we find 
\begin{equation}
q_{\lambda }^{L_{1}}Y(u,z)q_{\lambda }^{-L_{1}}=Y(q_{\lambda (1-\lambda
z)}^{L_{1}}(1-\lambda z)^{-2L_{0}}u,\frac{z}{1-\lambda z}).  \label{eq: Y_L1}
\end{equation}
Combining these it follows that the transformation 
$z\mapsto -\lambda^{2}z^{-1}$ is described by 
$T_{\lambda }\equiv q_{\lambda}^{L_{-1}}q_{\lambda ^{-1}}^{L_{1}}q_{\lambda }^{L_{-1}}$ with%
\begin{equation}
T_{\lambda }Y(u,z)T_{\lambda }^{-1}=
Y(q_{-z\lambda ^{-2}}^{L_{1}}(\lambda^{-2}z^{2})^{-L_{0}}u,-\lambda ^{2}z^{-1}).  \label{eq: Y_U}
\end{equation}
Taking $\lambda =\sqrt{-1}$ in (\ref{eq: Y_U}) corresponding to the
inversion $z\mapsto z^{-1}$ we find 
\begin{equation}
Y^{\dagger }(u,z)\equiv T_{\sqrt{-1}}Y(u,z)T_{\sqrt{-1}%
}^{-1}=Y(q_{z}^{L_{1}}(-z^{2})^{-L_{0}}u,z^{-1}).  \label{eq: adj op}
\end{equation}%
We call $Y^{\dagger }(u,z)$ the \emph{adjoint} vertex operator\footnote{%
This terminology differs from that of \cite{FHL}}. For $u$ of weight $wt(u)$
it follows that $Y^{\dagger }(u,z)=\sum_{n}u_{n}^{\dagger }z^{-n-1}$ has
modes 
\begin{equation}
u_{n}^{\dagger }=(-1)^{wt(u)}\sum_{k=0}^{wt(u)}\frac{1}{k!}%
(L_{1}^{k}u)_{2wt(u)-n-k-2}.  \label{eq: udagger}
\end{equation}

\medskip 

For a quasi-primary state $u$ (\ref{eq: udagger}) simplifies to 
\begin{equation}
u_{n}^{\dagger }=(-1)^{wt(u)}u_{2wt(u)-n-2}.  \label{eq: adj op qp}
\end{equation}%
Thus for a weight one Heisenberg vector $a$ we find%
\begin{equation}
a_{n}^{\dagger }=-a_{-n}.  \label{eq: adagger}
\end{equation}%
and for the\ weight two Virasoro vector $\omega $ we find that for $%
L_{n}^{\dagger }\equiv \omega _{n+1}^{\dagger }$  
\begin{equation}
L_{n}^{\dagger }=L_{-n}.  \label{eq: LLdagger}
\end{equation}
\medskip

We also note that the adjoint vertex operators can be used to construct a
canonical $V$-module as follows. Define vertex operators $Y_{V^{\prime
}}:V\rightarrow \mathcal{F}(V^{\prime })$ by%
\begin{equation}
\langle Y^{\prime }(u,z)a^{\prime },b\rangle _{d}=\langle a^{\prime
},Y^{\dagger }(u,z)b\rangle _{d},  \label{eq: VOAdual}
\end{equation}%
for $a^{\prime }\in V^{\prime }$ and $b,u\in V$. Then $(V^{\prime
},Y_{V^{\prime }})$ can be shown to be a $V$-module called the \emph{dual or
contragradient module} \cite{FHL}.
\bigskip 

\begin{exercise}
\label{Exercise Tlambda} Prove (\ref{eq: Y_U}).
\end{exercise}

\begin{exercise}
\label{Exercise SL2C} Show for a quasi-primary state $u$ (i.e. $L_{1}u=0$)
of weight $wt(u)$ that under a M\"{o}bius transformation $z\rightarrow \phi
(z)=\frac{az+b}{cz+d}$ 
\begin{equation}
Y(u,z)\rightarrow \left( \frac{d\phi }{dz}\right) ^{wt(u)}Y(u,\phi (z)).
\label{eq: Y_Mobius}
\end{equation}
\end{exercise}

\begin{exercise}
\label{Exercise F0form} Hence show for $n$ quasi-primary vectors $u^{i}$ of
weight $wt(u^{i})$ that the rational $n$-point function (\ref{eq: genus0npt}%
) is associated with a (formal) M\"{o}bius-invariant differential form on 
$\mathbb{CP}^{1}$ 
\begin{equation}
\mathcal{F}_{V}^{(0)}(u^{1},\ldots ,u^{n})=F_{V}^{(0)}((u^{1},z_{1}),\ldots
,(u^{n},z_{n}))\prod_{1\leq i\leq n}dz_{i}^{wt(u^{i})}.
\label{eq: Fcal0form}
\end{equation}
\end{exercise}

\begin{remark}
\label{Remark Fcal0form}$\mathcal{F}_{V}^{(0)}(u^{1},\ldots ,u^{n})$ is a
conformally invariant global meromorphic differential form on $\mathbb{CP}%
^{1} $ if $u^{1},\ldots ,u^{n}$ are primary vectors i.e. $L_{n}u^{i}=0$ for
all $n>0$.
\end{remark}

\begin{exercise}
\label{Exercise udagger}Prove (\ref{eq: udagger}).
\end{exercise}

\begin{exercise}
\label{Exercise udagdag}Show that $(Y^{\dagger }(u,z))^{\dagger }=Y(u,z)$.
\end{exercise}

\subsection{Invariant Bilinear Forms }\label{Subsect_Invariant bilinear}
\label{subsect_LiZ} 
In this Subsection we consider the construction of a canonical
bilinear form on $V$ motivated by (\ref{eq: VOAdual}). We say a bilinear form%
\emph{\ }$\langle \ ,\rangle :V\times V{\longrightarrow }\mathbb{C}$ is 
\emph{invariant }if for all $a,b,u\in V$ 
\begin{equation}
\langle Y(u,z)a,b\rangle =\langle a,Y^{\dagger }(u,z)b\rangle ,
\label{eq: inv bil form}
\end{equation}%
with $Y^{\dagger }(a,z)$ the adjoint operator of (\ref{eq: adj op}). In
terms of modes, (\ref{eq: inv bil form}) reads 
\begin{equation}
\langle u_{n}a,b\rangle =\langle a,u_{n}^{\dagger }b\rangle ,
\label{eq: abuudagger}
\end{equation}%
using (\ref{eq: udagger}). Applying (\ref{eq: LLdagger}) it follows that 
\begin{equation}
\langle L_{0}a,b\rangle =\langle a,L_{0}b\rangle .  \label{eq: L0LiZ}
\end{equation}%
Thus for homogeneous $a$ and $b$ then $\langle a,b\rangle =0$ for 
$wt(a)\not=wt(b)$.

\bigskip

Next consider $a,b$ with $wt(a)=wt(b)$. Invariance and skew-symmetry (cf. Exercise \ref{Exercise_skew}) give 
\begin{eqnarray*}
\langle \mathbf{1},Y^{\dag }(a,z)b\rangle &=&(-z^{2})^{-wt(a)}\langle 
\mathbf{1},Y(q_{z}^{L_{1}}a,z^{-1})b\rangle \\
&=&(-z^{2})^{-wt(b)}\langle \mathbf{1}%
,q_{z^{-1}}^{L_{-1}}Y(b,-z^{-1})q_{z}^{L_{1}}a\rangle \\
&=&\langle \mathbf{1},q_{z^{-1}}^{L_{-1}}Y^{\dag
}(q_{z}^{L_{1}}b,-z)q_{z}^{L_{1}}a\rangle .
\end{eqnarray*}%
But (\ref{eq: LLdagger}) implies this is%
\begin{equation*}
\langle q_{z^{-1}}^{L_{1}}\mathbf{1},Y^{\dag
}(q_{z}^{L_{1}}b,-z)q_{z}^{L_{1}}a\rangle =\langle \mathbf{1},Y^{\dag
}(q_{z}^{L_{1}}b,-z)q_{z}^{L_{1}}a\rangle .
\end{equation*}%
Using invariance this becomes 
\begin{equation*}
\langle Y(q_{z}^{L_{1}}b,-z)\mathbf{1},q_{z}^{L_{1}}a\rangle .
\end{equation*}%
Finally, using Exercise \ref{Exercise_2.3.2} and (\ref{eq: LLdagger}) this is 
\begin{equation*}
\langle q_{-z}^{L_{-1}}q_{z}^{L_{1}}b,q_{z}^{L_{1}}a\rangle =\langle
b,q_{z}^{L_{-1}}a\rangle =\langle b,Y(a,z)\mathbf{1}\rangle .
\end{equation*}%
Thus we have shown 
\begin{equation*}
\langle Y(a,z)\mathbf{1},b\rangle =\langle b,Y(a,z)\mathbf{1}\rangle .
\end{equation*}%
In particular, considering the $z^{0}$ term, this implies that the bilinear
form is symmetric:%
\begin{equation}
\langle a,b\rangle =\langle b,a\rangle .  \label{eq: LiZsym}
\end{equation}

\bigskip

Consider again $a,b$ with $wt(a)=wt(b)$. Using the creation axiom $a_{-1}%
\mathbf{1}=a$ we obtain 
\begin{equation}
\langle a,b\rangle =\langle \mathbf{1},a_{-1}^{\dagger }b\rangle .
\label{eq: modeLiZ}
\end{equation}%
with $a_{-1}^{\dagger }b\in V_{0}$. Thanks to the assumption that $V$ is of CFT-type\footnote{The general situation  is discussed in \cite{Li}). } we have 
$a_{-1}^{\dagger}b=\alpha \mathbf{1}$ for some $\alpha \in \mathbb{C}$ 
with $\langle a,b\rangle =\alpha \langle \mathbf{1},\mathbf{1}\rangle $. Hence either 
$\langle \mathbf{1},\mathbf{1}\rangle =0$ so that $\langle a,b\rangle =0$ for
all $a,b$ or else $\langle a,b\rangle $ is non-trivial and is uniquely
determined up to the value of $\langle \mathbf{1},\mathbf{1}\rangle \neq 0$
in which case we choose the normalization $\langle \mathbf{1},\mathbf{1}%
\rangle =1$.

\bigskip

It is straightforward to show that if $\langle \mathbf{1},\mathbf{1}\rangle
\neq 0$ then $L_{1}V_{1}=0$ (cf. Exercise \ref{Exercise L1V1}). Li has shown \cite{Li} that the converse is also true i.e. for a VOA of CFT-type, then $\langle \mathbf{1},\mathbf{1}%
\rangle \neq 0$ if and only if $L_{1}V_{1}=0$. We say that a VOA is of 
\emph{Strong CFT-type} if it is of CFT-type and $L_{1}V_{1}=0$. Such a VOA
therefore has a \emph{unique} normalized invariant bilinear form.

\bigskip

$\langle \ ,\ \rangle $ determines a standard map from $V$ to the restricted dual space 
$V^{\prime }$ defined by 
\begin{equation}
a\mapsto \langle a\ ,\cdot \ \rangle .  \label{eq: dualmap}
\end{equation}%
Let $\mathcal{K}$ denote the kernel of this map. $\langle \ ,\ \rangle $ is
nondegenerate with $\mathcal{K}$ trivial if, and only if,  $V$ is isomorphic to $V^{\prime
}$ i.e. $V$ is self-dual. In this case, we may identify $\langle \ ,\
\rangle $ with the canonical pairing $\langle \ ,\ \rangle _{d}$ and the
dual module (\ref{eq: VOAdual}) is isomorphic to the original VOA.
\medskip

The nondegeneracy of $\langle \ ,\ \rangle $ is also related to the
simplicity of the VOA $V$ in much that same way that nondegeneracy of the
Killing form determines semi-simplicity in Lie theory \cite{Li}. Let $%
\mathcal{I}\subset V$ denote the maximal proper ideal of $V$ (cf. Exercise \ref{Exercise_6.1.8}),
so that 
\begin{equation}
u_{n}b\in \mathcal{I},  \label{eq: unI}
\end{equation}%
for all $b\in \mathcal{I}$, $u\in V$. $V$ is  simple if $%
\mathcal{I}$ is trivial (cf. Section \ref{Section_Modules}). We now show that assuming $V$ is of
strong CFT-type then $\mathcal{I}=$ $\mathcal{K}$ and hence $V$ is simple
if, and only if, \ $\langle \ ,\ \rangle $ is nondegenerate.

We firstly note that $\mathbf{1}\notin \mathcal{I}$ (otherwise $u=u_{-1}%
\mathbf{1}\in \mathcal{I}$ for all $u\in V$). Because $V$ is of CFT-type, then
for all $b\in \mathcal{I}$ it follows $b\notin V_{0}$ and so 
\begin{equation}
\langle \mathbf{1},b\rangle =0.  \label{eq: 1bzero}
\end{equation}%
Consider $u\in V$ and $b\in \mathcal{I}$. Then $u_{-1}^{\dag }b$ $\in 
\mathcal{I}$ from (\ref{eq: udagger}) and so 
\begin{equation*}
\langle u,b\rangle =\langle \mathbf{1},u_{-1}^{\dag }b\rangle =0,
\end{equation*}%
for all $u$ from (\ref{eq: 1bzero}). Hence we find $\mathcal{I}\subseteq $ $%
\mathcal{K}$. Conversely, suppose that $c\in \mathcal{K}$. Then 
\begin{equation*}
\langle Y^{\dag }(u,z)v,c\rangle =0,
\end{equation*}%
for all $u,v\in V$. Invariance implies $\langle v,Y(u,z)c\rangle =0$ and
hence $u_{n}c\in \mathcal{K}$ for all $u_{n}$. But given $V$ is of strong
CFT-type then $\langle \ ,\ \rangle $ is nontrivial so that $\mathcal{K}\neq
V$ and hence $\mathcal{K}\subseteq \mathcal{I}$. \ Thus we conclude $%
\mathcal{I}=\mathcal{K}$.

\bigskip

Altogether we may summarize these results as follows:

\begin{theorem}\label{theorem: LiZ}
Let $V$ be a VOA. An invariant bilinear form $\langle \ ,\ \rangle $
on $V$ is
symmetric and diagonal with respect to the canonical $L_0$-grading. Furthermore, if
$V$ is of strong CFT-type, $\langle \ ,\ \rangle $ is unique up to
normalization and is nondegenerate if and only if $V$ is simple. 
\end{theorem}

The invariant bilinear form is equivalent to the chiral part of the
Zamolodchikov metric in CFT (\cite{BPZ, FMS, P}) where (abusing notation) 
\begin{eqnarray}
\langle a,b\rangle &=&\lim_{z_{1}\rightarrow 0}\lim_{z_{2}\rightarrow
0}\langle Y(a,z_{1})\mathbf{1},Y(b,z_{2})\mathbf{1}\rangle  \notag \\
&=&\lim_{z_{1}\rightarrow 0}\lim_{z_{2}\rightarrow 0}\langle \mathbf{1}%
,Y^{\dagger }(a,z_{1})Y(b,z_{2})\mathbf{1}\rangle  \notag \\
&=&"\langle \mathbf{1},Y(a,w_{1}=\infty )Y(b,z_{2}=0)\mathbf{1}\rangle ",
\label{Zam}
\end{eqnarray}%
for $w_{1}=1/z_{1}$ following (\ref{eq: adj op}). We thus refer to the
nondegenerate bilinear form as the \emph{Li-Zamolodchikov metric} on $V$ or
LiZ-metric for short\footnote{%
Although we use the term metric here, the bilinear form is not necessarily
positive definite.}.

\bigskip

Consider the rank one Heisenberg VOA $M_0$ generated by a weight one state $a$
with $V$ spanned by Fock vectors 
\begin{equation}
v=a_{-1}^{e_{1}}a_{-2}^{e_{2}}\ldots a_{-p}^{e_{p}}\mathbf{1},
\label{eq: Fockstate}
\end{equation}%
for non-negative integers $e_{i}$. Using (\ref{eq: adagger}), we find that
the Fock basis consisting of vectors of the form (\ref{eq: Fockstate}) is
orthogonal with respect to the LiZ-metric with 
\begin{equation}
\langle v,v\rangle =\prod_{1\leq i\leq p}(-i)^{e_{i}}e_{i}!.
\label{eq: inner prod}
\end{equation}%
Clearly $\langle \ ,\ \rangle $ is non-degenerate so by Theorem \ref{theorem: LiZ} it follows that $M_0$ is a simple VOA (as already discussed  in Section \ref{Section_Modules}).

\bigskip 

Consider the Virasoro VOA $\mbox{Vir}_{c}$ generated by the Virasoro vector $\omega $
of central charge $c$. Using (\ref{eq: L0LiZ}) it is sufficient to consider
the non-degeneracy of $\langle \ ,\ \rangle $ on each homogeneous space $%
V_{n}$. In particular, let $M_{n}(c)=(\langle a,b\rangle )$ be the Gram
matrix of $(\mbox{Vir}_{c})_{n}$ with respect to some basis. The \emph{Kac determinant}
is $\det M_{n}(c)$ (\cite{KR}), which is conveniently considered as a polynomial in $c$.
By Theorem \ref{theorem: LiZ}, $\mbox{Vir}_{c}$ is simple if, and only if, $\det M_{n}(c)\neq 0$ for all $n$. 
For $n=2$ we have $V_{2}=\mathbb{C}\omega $  with
Kac determinant 
\begin{equation}
\det M_{2}(c)=\langle \omega ,\omega \rangle =\langle \mathbf{1},L_{2}L_{-2}%
\mathbf{1}\rangle =\frac{c}{2},  \label{eq: detM2}
\end{equation}%
with a zero at $c=0$. For $n=4$ we have 
$V_{4}=\mathbb{C} \langle L_{-2}^{2}\mathbf{1},L_{-4}\mathbf{1} \rangle$ with%
\begin{equation}
M_{4}(c)=\left[ 
\begin{array}{cc}
c(4+\frac{1}{2}c) & 3c \\ 
3c & 5c%
\end{array}%
\right] ,  \label{eq M4}
\end{equation}%
and Kac determinant 
\begin{equation}
\det M_{4}(c)=\frac{1}{2}c^{2}\left( 5c+22\right)  \label{eq: detM4}
\end{equation}%
with zeros at $c=0,-\frac{22}{5}$. 

\medskip
There is a general formula for the Kac
determinant $\det M_{n}(c)$ which turns out to have zeros for central charge 
\begin{equation}
c = c_{p,q}=1-\frac{6(p-q)^{2}}{pq},  \label{eq: Cpq}
\end{equation}%
where $(p-1)(q-1)=n$ for coprime $p,q\geq 2$. Thus $\mbox{Vir}_{c}$ is a simple VOA
iff $c\neq c_{p,q}$ for some coprime $p,q\geq 2$ (cf. Theorem \ref{thmratVir}).

\begin{exercise}
\label{Exercise L1V1} Show that if $\langle \mathbf{1},\mathbf{1}\rangle
\neq 0$ then $L_{1}V_{1}=0$.
\end{exercise}

\begin{exercise}
\label{Exercise YaYb} Suppose that $a\in V_{m}, b\in V_{n}$ and at least one of $a$ or $b$
is quasi-primary. Prove that the $2$-point correlation function is given by%
\begin{equation*}
\langle \mathbf{1},Y(a,z_{1})Y(b,z_{2})\mathbf{1}\rangle =\frac{\langle
a,b\rangle }{(z_{1}-z_{2})^{2m}}\delta _{m,n}.
\end{equation*}%
(The Zamolodchikov metric\ in CFT is often introduced in this way.)
\end{exercise}

\begin{exercise}
\label{Exercise vv Fock}Verify (\ref{eq: inner prod}).
\end{exercise}

\section{The Genus Two Partition Function for the Heisenberg VOA}
\label{Section_GenusTwo}

In this section we will discuss some recent research by the authors wherein
we develop a theory of partition and $n$-point correlation functions on a
Riemann surface of genus two \cite{T1, MT2, MT3, MT4}. The basic idea is to
construct a genus two Riemann surface by specific sewing schemes where we
either sew together two once punctured tori or self-sew a twice punctured
torus (i.e. attach a handle). The partition and $n$-point functions on the
genus two surface are then defined in terms of correlation functions on the
lower genus surfaces combined together in an appropriate way. We will not
explore the full details entailed in this programme. Instead we will
consider the example of the Heisenberg VOA $M_0$ and compute the partition
function on the genus two surface formed from two tori.

\subsection{Genus One Heisenberg 1-Point Functions}
\label{Subsect_Heisen1pt}

We first discuss the genus one 1-point correlation function for all elements
of the Heisenberg VOA $M_0$ generated by the weight one Heisenberg vector $a$ 
\cite{MT1}. We make heavy use of the Zhu recursion formulas I and II of
Theorems \ref{TheoremZhu} and \ref{Proposition_1pt}. In particular, we prove Theorem \ref{thm4.15} by
considering the $1$-point function $Z_{M_0}(v,\tau )$ for a Fock vector in the
square bracket formulation 
\begin{equation}
v=a[-k_{1}]\ldots a[-k_{n}]\mathbf{1},  \label{eq: vFovkki}
\end{equation}%
for $k_{i}\geq 1$. The Fock vector $v$ is of square bracket weight $%
wt[v]=\sum_{i}k_{i}$. We want to show that%
\begin{equation}
Z_{M_0}^{(1)}(v,\tau )=\frac{Q_{v}(\tau )}{\eta (q)},
\label{eq: genus1 ptn func}
\end{equation}%
for $Q_{v}(\tau )$ $\in \mathfrak{Q}$, the algebra of quasimodular forms. $%
Q_{v}(\tau )$ is of weight $wt[v]$ and is expressed in terms of 
\begin{equation}
C(k,l)=C(k,l,\tau )=(-1)^{l+1}\frac{(k+l-1)!}{(k-1)!(l-1)!}E_{k+l}(\tau ),
\label{eq: Ckldef}
\end{equation}%
for $k,l\geq 1$. Here $E_{n}(\tau )$ is the Eisenstein series of (\ref{Eisen1}). We
recall that $E_{n}=0$ for $n$ odd, $E_{2}(\tau )$ is a quasimodular form of
weight 2 and $E_{n}$ is a modular form of weight $n$ for even $n\geq 4$.
Thus $C(k,l,\tau )$ is a quasi-modular form of weight $k+l$. We also note
that $C(k,l)=C(l,k)$.
\medskip

Each Fock vector $v$ is described by a label set $\Phi _{\lambda
}=\{k_{1},\ldots ,k_{n}\}$ which corresponds in a natural 1-1 manner with
unrestricted partitions $\lambda =\{1^{e_{1}}, 2^{e_{2}}, \ldots \}$ of $wt[v]$
(where $e_{i}\geq 0$). We write $v=v(\lambda )$ to indicate this
correspondence, which will play a significant r\^{o}le later on. Define $F(\Phi
_{\lambda })$ to be the subset of all permutations on $\Phi _{\lambda }$
consisting only of \emph{fixed-point-free involutions}. Let $\varphi =\ldots
(k_{i}k_{j})\ldots $, a product of disjoint cycles, denote a typical element
of $F(\Phi _{\lambda })$.
\medskip

We can now describe the 1-point function $Z_{M_0}^{(1)}(v(\lambda ),\tau )$ of
(\ref{eq: genus1 ptn func}) \cite{MT1}

\begin{theorem}
\label{theorem: HeisenGen1}For even $n$ 
\begin{eqnarray}
Q_{v}(\tau ) &=&\sum_{\phi \in F(\Phi _{\lambda })}\Gamma (\phi ,\tau ),
\label{eq: Qvtau} \\
\Gamma (\phi ,\tau ) &=&\prod_{(k_{i}k_{j})}C(k_{i},k_{j},\tau ).
\label{eq: Gammaphi}
\end{eqnarray}%
for $C$ of (\ref{eq: Ckldef}) and where the product ranges over all the
cycles of $\varphi =\ldots (k_{i}k_{j})\ldots \in F(\Phi _{\lambda })$. $%
Q_{v}(\tau )$ $\in \mathfrak{Q}$ and is of weight $wt[v]$. For $n$ odd $%
Q_{v}(\tau )$ vanishes.
\end{theorem}

\noindent \textbf{Proof.} Let $v(\lambda
)=a[-k_{1}]w$ for $w=a[-k_{2}]\ldots a[-k_{n}]\mathbf{1}$ and use the Zhu
recursion formula II of Theorem  \ref{Proposition_1pt} to find%
\begin{eqnarray*}
Z_{M_0}^{(1)}(a[-k_{1}]w,\tau ) &=&\delta _{k_{1},1}\mathrm{Tr}%
_{M_0}(o(a)o(w)q^{L_{0}-1/24}) \\
&&+\sum_{m\geq 1}(-1)^{m+1}\left( 
\begin{array}{c}
k_{1}+m-1 \\ 
m%
\end{array}%
\right) E_{k_{1}+m}(\tau )Z_{M_0}^{(1)}(a[m]w,\tau ).
\end{eqnarray*}%
But $o(a)u=0$ for all $u\in M$ and the Heisenberg commutation relations
imply 
\begin{eqnarray*}
Z_{M_0}^{(1)}(a[-k_{1}]w,\tau ) &=&0+\sum_{j=2}^{n}(-1)^{k_{j}+1}\left( 
\begin{array}{c}
k_{1}+k_{j}-1 \\ 
k_{j}%
\end{array}%
\right) E_{k_{1}+k_{j}}(\tau )k_{j}Z_{M_0}^{(1)}(\hat{w},\tau ) \\
&=&\sum_{j=2}^{n}C(k_{1},k_{j},\tau )Z_{M_0}^{(1)}(\hat{w},\tau ),
\end{eqnarray*}%
where $\hat{w}$ denotes the Fock vector with label set $\{k_{2},\ldots ,\hat{%
k}_{j}\ldots ,k_{n}\}$ with the index $k_{j}$ deleted. The result follows by
repeated application of this recursive formula until we obtain $\hat{w}=%
\mathbf{1}$ for which $Z_{M_0}^{(1)}(\mathbf{1},\tau )=$ $\frac{1}{\eta (q)}$.
The resulting expression for $Q_{v}(\tau )$ is clearly a quasimodular form
of weight $wt[v]=\sum_{i}k_{i}$. Thus Theorem \ref{theorem: HeisenGen1} follows. 

\bigskip

Some further insight into the combinatorial structure of $Q_{v}(\tau )$ can
be garnered by a consideration of the $n$-point function for $n$ Heisenberg
vectors which we denote by 
\begin{equation}
G_{n}^{(1)}(z_{1},\ldots ,z_{n},\tau )\equiv F_{M_0}^{(1)}((a,z_{1}),\ldots
(a,z),\tau ).  \label{eq: Gndef}
\end{equation}%
This is a symmetric function in $z_{i}$ with a pole of order two at $%
z_{i}=z_{j}$ for all $i\neq j$ (from locality). For $n=1$ we immediately
find 
\begin{equation*}
G_{1}^{(1)}(z_{1},\tau )=\mathrm{Tr}_{M_0}o(a)q^{L_{0}-1/24}=0.
\end{equation*}%
The 2-point function is easily computed via the Zhu recursion formula I of
Theorem  \ref{TheoremZhu} to obtain%
\begin{eqnarray}
G_{2}^{(1)}(z_{1},z_{2},\tau ) &=&\mathrm{Tr}_{M_0}o(a)o(a)q^{L_{0}-1/24}-%
\sum_{m\geq 1}\frac{(-1)^{m}}{m!}P_{1}^{(m)}(z_{12},\tau
)Z_{M_0}^{(1)}(a[m]a,\tau )  \notag \\
&=&0+P_{2}(z_{12},\tau )\frac{1}{\eta (q)},  \label{eq: G2mod}
\end{eqnarray}%
since $a[m]a=\mathbf{1}\delta _{m,1}$ and where, from Theorem \ref{thm5.1}, we recall 
\begin{equation*}
P_{2}(z,\tau )=\frac{d}{dz}P_{1}(z,\tau )=\frac{1}{z^{2}}+\sum_{n=2}^{\infty
}(n-1)E_{n}(\tau )z^{n-2}.
\end{equation*}%
(\ref{eq: G2mod}) is the elliptic analogue of the genus zero formula (\ref%
{eq: F2Heisenberg}) and reflects a deeper geometrical structure underlying
the Heisenberg VOA e.g. \cite{MT4}. \ 

\bigskip

Using the $n$-point correlation function version of the Zhu recursion
formula I we can similarly obtain the genus one analogue of Theorem \ref%
{theorem: HeisenGenFun} to find \cite{MT1}

\begin{theorem}
\label{theorem: HeisenGenFunGenusOne}For $n$ even%
\begin{equation}
G_{n}^{(1)}(z_{1},\ldots ,z_{n},\tau )=\frac{1}{\eta (q)}\sum_{\varphi \in
F(\Phi )}\prod_{(ij)}P_{2}(z_{ij},\tau ),  \label{eq: FHeisenGenusone}
\end{equation}%
where the product ranges over all the cycles of $\varphi =\ldots (ij)\ldots $
\ for $\Phi =\{1,2,\ldots ,n\}$ whereas for $n$ odd $G_{n}^{(1)}$ vanishes. 
\end{theorem}

We may use this result to compute \emph{any} genus one $n$-point correlation
function for $M_0$ by a consideration of an appropriate analytic expansion of $%
G_{n}^{(1)}(z_{1},\ldots ,z_{n},\tau )$ \cite{MT1}. In particular, we can
re-derive (\ref{eq: Qvtau}) by making use of the following identity 
\begin{eqnarray}
G_{n}^{(1)}(z_{1},\ldots ,z_{n},\tau ) &=&Z_{M_0}^{(1)}(Y[a,z_{1}]\ldots
Y[a,z_{n}]\mathbf{1},\tau )  \notag \\
&=&\sum_{k_{1},\ldots k_{n}}Z_{M_0}^{(1)}(v,\tau )z_{1}^{k_{1}-1}\ldots
z_{n}^{k_{n}-1},  \label{eq: Gn1pt}
\end{eqnarray}%
for Fock vector $v=a[-k_{1}]\ldots a[-k_{n}]\mathbf{1}$ for all $k_{i}$. We
may extract the non-negative values of $k_{1},\ldots ,k_{n}$ from the
expansion 
\begin{equation}
P_{2}(z_{ij},\tau )=\frac{1}{(z_{i}-z_{j})^{2}}+\sum_{k_{i},k_{j}\geq
1}^{\infty }C(k_{i},k_{j},\tau )z_{i}^{k_{i}-1}z_{j}^{k_{j}-1},
\label{eq: P2Ckl}
\end{equation}%
for $C$ of (\ref{eq: Ckldef}). Thus (\ref{eq: FHeisenGenusone}) implies the
formula (\ref{eq: Qvtau}) of Theorem \ref{theorem: HeisenGen1} found for $Q_{v}(\tau )$.

\bigskip 

It is very useful to recast Theorem \ref{theorem: HeisenGen1} in terms of
graph theory as follows. Consider a Fock vector $v(\lambda )$ with label set 
$\Phi _{\lambda }=\{k_{1},\ldots ,k_{n}\}$ and let $\phi \in F(\Phi
_{\lambda })$ be a fixed-point-free involution of $\Phi _{\lambda }$ leading
to a contribution $\Gamma (\phi ,\tau )$ to $Q_{v}(\tau )$ in (\ref{eq:
Gammaphi}). We may then associate to each $\phi \in F(\Phi _{\lambda })$ a $%
\phi $\emph{-graph} $\gamma _{\phi \text{ }}$consisting of $n$ vertices
labelled by $\Phi _{\lambda }$ of unit valence with $n/2$ unoriented edges
connecting the pairs of vertices $(k_{i},k_{j})$ determined by $\varphi
=\ldots (k_{i}k_{j})\ldots $. Following Exercise \ref{Exercise F(Phi)} there
are $(n-1)!!$ such graphs for a given label set $\Phi _{\lambda }$. Thus, in
Exercise \ref{Example VFock} with $v=a[-1]^{3}a[-2]^{2}a[-5]\mathbf{1}$
there are 15 independent $\phi$-graphs (cf. Exercise \ref{Exercise Fockgraphs}). 
A $\phi -$graph for a fixed point
involution $\phi =(11)(22)(15)$ is shown below\footnote{%
Note that there are 3 distinct fixed point involutions notated by $%
(11)(22)(15)$.}.

\begin{center}
\begin{picture}(250,80)

\put(100,50){\line(1,2){10}}
\put(90,50){\makebox(0,0){1}}
\put(100,50){\circle*{4}}

\put(105,78){\makebox(0,0){1}}
\put(110,70){\circle*{4}}

\put(130,70){\line(1,-2){10}}
\put(139,78){\makebox(0,0){2}}
\put(130,70){\circle*{4}}

\put(150,50){\makebox(0,0){2}}
\put(140,50){\circle*{4}}

\put(110,30){\line(1,0){20}}
\put(135,20){\makebox(0,0){5}}
\put(130,30){\circle*{4}}

\put(105,20){\makebox(0,0){1}}
\put(110,30){\circle*{4}}
\end{picture}

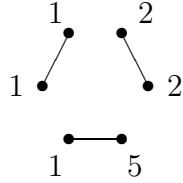
{\small Fig.~1 }
\end{center}

\bigskip 

Given a $\phi $-graph $\gamma _{\phi}$ we define a 
\emph{weight function} 
\begin{equation*}
\kappa :\{\gamma _{\phi \text{ }}\}\longrightarrow \mathfrak{Q},
\end{equation*}%
as follows: for every edge $E$ labeled as \ $\overset{k}{\bullet }-\overset{l%
}{\bullet }$ define 
\begin{equation}
\kappa (E,\tau )=C(k,l,\tau ),  \label{eq: kappaE}
\end{equation}%
with 
\begin{equation}
\kappa (\gamma _{\phi },\tau )=\prod \kappa (E,\tau ),  \label{eq: kappagam}
\end{equation}%
where the product is taken over all edges of $\gamma _{\phi}$. Thus
the $\phi$-graph of Fig.~1 has weight $C(1,1)C(2,2)C(1,5)=-30E_{2}(\tau
)E_{4}(\tau )E_{6}(\tau )$.

\medskip
Clearly Theorem \ref{theorem: HeisenGen1} can now be restated in terms of
graphs:

\begin{theorem}
\label{theorem: HeisenGraph}For a Fock vector $v(\lambda )$ with label set $%
\Phi _{\lambda }=\{k_{1},\ldots ,k_{n}\}$ and even $n$ 
\begin{equation}
Q_{v}(\tau )=\sum_{\gamma _{\phi}}\kappa (\gamma _{\phi},\tau ),  \label{eq: Qvkappa}
\end{equation}%
where the sum is taken over all independent $\phi -$graphs for $\Phi
_{\lambda }$. 
\end{theorem}
\bigskip

\begin{exercise}
\label{Example VFock}For $v=a[-1]^{3}a[-2]^{2}a[-5]\mathbf{1}$\ of weight $%
wt[v]=12$ with $\Phi _{\lambda }=\{1,1,1,2,2,5\}$ and $\left\vert F(\Phi
_{\lambda })\right\vert =5!!=15$ (cf. Exercise \ref{Exercise F(Phi)}) show
that%
\begin{eqnarray*}
Q_{v}(\tau ) &=&6C(1,1)C(1,2)C(2,5)+3C(1,1)C(2,2)C(1,5)+6C(1,2)^{2}C(1,5) \\
&=&0-90E_{2}(\tau )E_{4}(\tau )E_{6}(\tau )+0.
\end{eqnarray*}%
Thus only 3 elements of $F(\Phi _{\lambda })$ make a non-zero contribution
to $Q_{v}(\tau )$.
\end{exercise}

\begin{exercise}
\label{Exercise Fockgraphs} Find all the $\phi$-graphs for $v=a[-1]^{3}a[-2]^{2}a[-5]\mathbf{1}$.
\end{exercise}

\subsection{Sewing Two Tori}\label{Subsect_sew}

\label{subsect_epsilon} In this section we digress from VOA theory to
briefly review some aspects of Riemann surface theory and the construction
of a genus two surface by sewing together two punctured tori. A genus two
Riemann surface can also be constructed by sewing a handle to a torus but we
do not consider that situation here. For more details see  \cite{MT2},
\cite{MT4}.
\medskip

Let $\mathcal{S}^{(2)}$ denote a compact Riemann surface of genus $2$ and
let $a_{1},a_{2},b_{1},b_{2}$ be the canonical homology basis e.g. \cite{FK}%
. There exists two holomorphic 1-forms $\nu _{i}$, $i=1,2$ which we may
normalize by 
\begin{equation*}
\oint_{a_{i}}\nu _{j}=2\pi i\delta _{ij}.
\end{equation*}%
The genus $2$ period matrix $\Omega $ is defined by 
\begin{equation}
\Omega _{ij}=\frac{1}{2\pi i}\oint_{b_{i}}\nu _{j},  \label{eq: period}
\end{equation}%
for $i,j=1,2$. Using the Riemann bilinear relations, one finds that $\Omega $
is a complex symmetric matrix with positive-definite imaginary part, i.e. $\Omega \in \mathfrak{H}_{2}$, the genus $2$ Siegel complex upper half-space.
\medskip

The intersection form $\Xi $ is a natural non-degenerate symplectic bilinear
form on the first homology group $H_{1}(\mathcal{S}^{(2)},\mathbb{Z})\cong 
\mathbb{Z}^{4}$, satisfying 
\begin{equation*}
\Xi (a_{i},a_{j})=\Xi (b_{i},b_{j})=0,\quad \Xi (a_{i},b_{j})=\delta
_{ij},\quad i,j=1,2.
\end{equation*}%
The mapping class group is given by the symplectic group \ 
\begin{eqnarray*}
Sp(4,\mathbb{Z}) &=&\{\gamma =\left( 
\begin{array}{ll}
A & B \\ 
C & D%
\end{array}%
\right) \in SL(4,\mathbb{Z})| \\
AB^{T} &=&BA^{T},CD^{T}=D^{T}C,AD^{T}-BC^{T}=I_{2}\},
\end{eqnarray*}%
where $A^{T}$ denotes the transpose of $A$. $Sp(4,\mathbb{Z})$ acts on $%
\mathfrak{H}_{2}$ via 
\begin{equation}
\gamma .\Omega {=(A{\Omega }+B)(C\Omega +D)^{-1},}  \label{eq: Sp4zOmega}
\end{equation}%
and naturally on $H_{1}(\mathcal{S},\mathbb{Z})$, where it preserves $\Xi $.

\bigskip 

We now briefly review a general method originally due to Yamada \cite{Y} and
discussed at length in \cite{MT2} for calculating the period matrix (and
other structures) on a Riemann surface formed by sewing together two other
Riemann surfaces. In particular, we wish to describe $\Omega _{ij}$ on a
genus two Riemann surface formed by sewing together two tori $\mathcal{S}_{a}
$ for $a=1,2$. Consider an oriented torus $\mathcal{S}_{a}=\mathbb{C}%
/\Lambda _{\tau _{a}}$ with lattice $\Lambda _{\tau _{a}}$ with basis $(2\pi
i,2\pi i\tau _{a})$ for $\tau _{a}\in \mathfrak{H}$, the complex upper
half plane. For local coordinate $z_{a}\in \mathcal{S}_{a}$ consider the
closed disk $\left\vert z_{a}\right\vert \leq r_{a}$. This is contained in $%
\mathcal{S}_{a}$ provided\ $r_{a}<\frac{1}{2}D(q_{a})$ where 
\begin{equation}
D(q_{a})=\min_{\lambda \in \Lambda _{\tau _{a}},\lambda \neq 0}|\lambda |,
\label{eq: D(q)}
\end{equation}%
is the minimal lattice distance.
\medskip

Introduce a sewing parameter $\epsilon \in \mathbb{C}$ where $|\epsilon
|\leq r_{1}r_{2}<\frac{1}{4}D(q_{1})D(q_{2})$ and excise the disk $%
\{z_{a},\left\vert z_{a}\right\vert \leq |\epsilon |/r_{\bar{a}}\}$ centered
at $z_{a}=0$ to form a punctured torus 
\begin{equation*}
\hat{\mathcal{S}}_{a}=\mathcal{S}_{a}\backslash \{z_{a},\left\vert
z_{a}\right\vert \leq |\epsilon |/r_{\bar{a}}\},
\end{equation*}%
where we use the convention 
\begin{equation}
\overline{1}=2,\quad \overline{2}=1.  \label{eq: abar}
\end{equation}%
Define the annulus

\begin{equation}
\mathcal{A}_{a}=\{z_{a},|\epsilon |/r_{\bar{a}}\leq \left\vert
z_{a}\right\vert \leq r_{a}\}\subset \hat{\mathcal{S}}_{a},
\label{eq: Annulus}
\end{equation}%
We then identify $\mathcal{A}_{1}$ with $\mathcal{A}_{2}$ via the sewing
relation 
\begin{equation}
z_{1}z_{2}=\epsilon ,  \label{eq: pinch}
\end{equation}%
to obtain an explicit construction of a genus two Riemann surface%
\begin{equation*}
\mathcal{S}^{(2)}=\hat{\mathcal{S}}_{1}\cup \hat{\mathcal{S}}_{2}\cup (%
\mathcal{A}_{1}\simeq \mathcal{A}_{2}),
\end{equation*}%
which is parameterized by the domain 
\begin{equation}
\mathcal{D}^{\epsilon }=\{(\tau _{1},\tau _{2},\epsilon )\in \mathfrak{H}%
\times\mathfrak{ H}\times\mathbb{ C}\ |\ |\epsilon |<\frac{1}{4}%
D(q_{1})D(q_{2})\}.  \label{eq: Deps}
\end{equation}

\begin{center}
\bigskip

\begin{picture}(300,100)

\put(50,50){\qbezier(10,-20)(-30,0)(10,20)}
\put(50,52){\qbezier(10,18)(50,35)(90,18)}

\put(50,48){\qbezier(10,-18)(50,-35)(90,-18)}

\put(45,50){\qbezier(25,0)(45,17)(60,0)}
\put(45,50){\qbezier(20,2)(45,-17)(65,2)}


\put(175,52){\qbezier(10,18)(50,35)(90,18)}
\put(175,50){\qbezier(90,20)(130,0)(90,-20)}
\put(175,48){\qbezier(10,-18)(50,-35)(90,-18)}

\put(200,50){\qbezier(25,0)(45,17)(60,0)}
\put(200,50){\qbezier(20,2)(45,-17)(65,2)}

\put(140,50){\circle{16}}
\put(140,50){\circle{40}}

\put(140,50){\vector(-1,-2){0}}
\put(50,50){\qbezier(90,0)(100,15)(90,30)}%
\put(140,90){\makebox(0,0){$z_1=0$}}

\put(140,50){\line(-1,1){14.1}}
\put(127,55){\makebox(0,0){$r_1$}}

\put(140,50){\line(1,0){8}}
\put(145,50){\vector(1,4){0}}
\put(55,20){\qbezier(90,4)(85,17)(90,30)}%
\put(150,15){\makebox(0,0){$|\epsilon|/r_2$}}

\put(25,50){\makebox(0,0){$\mathcal{S}_1$}}


\put(185,50){\circle{16}}
\put(185,50){\circle{40}}

\put(185,50){\vector(1,-2){0}}
\put(95,50){\qbezier(90,0)(80,15)(90,30)}%
\put(185,90){\makebox(0,0){$z_2=0$}}

\put(185,50){\line(-1,-1){14.1}}
\put(171,45){\makebox(0,0){$r_2$}}

\put(185,50){\line(1,0){8}}
\put(190,50){\vector(1,4){0}}
\put(90,20){\qbezier(100,4)(95,17)(100,30)}%
\put(190,15){\makebox(0,0){$ |\epsilon|/r_1$}}

\put(300,50){\makebox(0,0){$\mathcal{S}_2$}}

\end{picture}

{\small Fig.~2 Sewing Two Tori}
\end{center}

\bigskip 
In \cite{Y}, Yamada describes a general method for computing
the period matrix on the sewn Riemann surface $\mathcal{S}^{(2)}$ in terms
of data obtained from the two tori. This is described in detail in \cite%
{MT2} where we obtain the explicit form for $\Omega $ in terms of the
infinite matrix $A_{a}(\tau _{a},\epsilon )=(A_{a}(k,l,\tau
_{a},\epsilon ))$ for $k,l\geq 1$ where 
\begin{equation}
A_{a}(k,l,\tau _{a},\epsilon )=\frac{\epsilon ^{(k+l)/2}}{\sqrt{kl}}%
C(k,l,\tau _{a}),  \label{eq: Akldef}
\end{equation}%
and where $C(k,l,\tau _{a})$ is given in (\ref{eq: Ckldef}). Thus, dropping
the subscript,
\begin{equation*}
A(\tau ,\epsilon )=\left( 
\begin{array}{ccccc}
\epsilon E_{2}(\tau ) & 0 & \sqrt{3}\epsilon ^{2}E_{4}(\tau ) & 0 & \cdots
\\ 
0 & -3\epsilon ^{2}E_{4}(\tau ) & 0 & -5\sqrt{2}\epsilon ^{3}E_{6}(\tau ) & 
\cdots \\ 
\sqrt{3}\epsilon ^{2}E_{4}(\tau ) & 0 & 10\epsilon ^{3}E_{6}(\tau ) & 0 & 
\cdots \\ 
0 & -5\sqrt{2}\epsilon ^{3}E_{6}(\tau ) & 0 & -35\epsilon ^{4}E_{8}(\tau ) & 
\cdots \\ 
\vdots & \vdots & \vdots & \vdots & \ddots%
\end{array}%
\right) .
\end{equation*}%
The matrices $A_{1},A_{2}$ not only play a central r\^{o}le here but also later
on in our discussion of the genus two partition for the Heisenberg VOA $M_0$.
In particular, the matrix $I-A_{1}A_{2}$ and $\det (I-A_{1}A_{2})$ (where $I$
is the infinite identity matrix here) are important, where $\det
(I-A_{1}A_{2})$ is defined by 
\begin{eqnarray}
\log \det (I-A_{1}A_{2}) &=&\mathrm{Tr}\log (I-A_{1}A_{2})  \notag \\
&=&-\sum_{n\geq 1}\frac{1}{n}\mathrm{Tr}((A_{1}A_{2})^{n}).
\label{eq: logdet}
\end{eqnarray}%
These expressions are power series in $\mathfrak{Q}[[\epsilon]]$. One finds \cite{MT2}

\begin{theorem}
\label{Theorem_A1A2} \ \ 

(a) The infinite matrix 
\begin{equation}
(I-A_{1}A_{2})^{-1}=\sum_{n\geq 0}(A_{1}A_{2})^{n},  \label{eq: I_minus_A1A2}
\end{equation}%
is convergent for $(\tau _{1},\tau _{2},\epsilon )\in \mathcal{D}^{\epsilon
} $.

(b) $\det (I-A_{1}A_{2})$ is non-vanishing and holomorphic  on $\mathcal{D%
}^{\epsilon }$.
\end{theorem}

Furthermore we may obtain an explicit formula for the genus two period
matrix on $\mathcal{S}^{(2)}$

\begin{theorem}
\label{Theorem_period_eps} The sewing procedure determines a holomorphic map 
\begin{eqnarray}
F^{\epsilon }:\mathcal{D}^{\epsilon } &\rightarrow &\mathfrak{H}_{2},  \notag
\\
(\tau _{1},\tau _{2},\epsilon ) &\mapsto &\Omega (\tau _{1},\tau
_{2},\epsilon ),  \label{eq: Fepsmap}
\end{eqnarray}%
where $\Omega =\Omega (\tau _{1},\tau _{2},\epsilon )$ is given by 
\begin{eqnarray*}
2\pi i\Omega _{11} &=&2\pi i\tau _{1}+\epsilon
(A_{2}(I-A_{1}A_{2})^{-1})(1,1), \\
2\pi i\Omega _{22} &=&2\pi i\tau _{2}+\epsilon
(A_{1}(I-A_{2}A_{1})^{-1})(1,1), \\
2\pi i\Omega _{12} &=&-\epsilon (I-A_{1}A_{2})^{-1}(1,1).
\end{eqnarray*}%
Here $(1,1)$ refers to the $(1,1)$-entry of a matrix. 
\end{theorem}

\bigskip

$\mathcal{D}^{\epsilon }$ is preserved under the action of $G\simeq (SL(2,%
\mathbb{Z})$ $\times SL(2,\mathbb{Z}))\rtimes \mathbb{Z}_{2}$, the direct
product of the left and right torus modular groups, which are interchanged
upon conjugation by an involution $\beta $ as follows%
\begin{eqnarray}
\gamma _{1}.(\tau _{1},\tau _{2},\epsilon ) &=&(\frac{a_{1}\tau _{1}+b_{1}}{%
c_{1}\tau _{1}+d_{1}},\tau _{2},\frac{\epsilon }{c_{1}\tau _{1}+d_{1}}), 
\notag \\
\gamma _{2}.(\tau _{1},\tau _{2},\epsilon ) &=&(\tau _{1},\frac{a_{2}\tau
_{2}+b_{2}}{c_{2}\tau _{2}+d_{2}},\frac{\epsilon }{c_{2}\tau _{2}+d_{2}}), 
\notag \\
\beta .(\tau _{1},\tau _{2},\epsilon ) &=&(\tau _{2},\tau _{1},\epsilon ),
\label{eq: GDeps}
\end{eqnarray}%
for $(\gamma _{1},\gamma _{2})\in SL(2,\mathbb{Z})\times SL(2,\mathbb{Z})$
with $\gamma _{i}=\left( 
\begin{array}{cc}
a_{i} & b_{i} \\ 
c_{i} & d_{i}%
\end{array}%
\right) $.
\medskip

There is a natural injection $G\rightarrow Sp(4,\mathbb{Z})$ in which the
two $SL(2,\mathbb{Z})$ subgroups are mapped to 
\begin{equation}
\Gamma _{1}=\left\{ \left[ 
\begin{array}{cccc}
a_{1} & 0 & b_{1} & 0 \\ 
0 & 1 & 0 & 0 \\ 
c_{1} & 0 & d_{1} & 0 \\ 
0 & 0 & 0 & 1%
\end{array}%
\right] \right\} ,\;\Gamma _{2}=\left\{ \left[ 
\begin{array}{cccc}
1 & 0 & 0 & 0 \\ 
0 & a_{2} & 0 & b_{2} \\ 
0 & 0 & 1 & 0 \\ 
0 & c_{2} & 0 & d_{2}%
\end{array}%
\right] \right\} ,  \label{eq: G1G2}
\end{equation}%
and the involution is mapped to 
\begin{equation}
\beta =\left[ 
\begin{array}{cccc}
0 & 1 & 0 & 0 \\ 
1 & 0 & 0 & 0 \\ 
0 & 0 & 0 & 1 \\ 
0 & 0 & 1 & 0%
\end{array}%
\right] .  \label{eq: beta}
\end{equation}%
Thus as a subgroup of $Sp(4,\mathbb{Z})$, $G$ also has a natural action on
the Siegel upper half plane $\mathfrak{H}_{2}$ as given in (\ref{eq: Sp4zOmega}%
). This action is compatible with respect to the map (\ref{eq: Fepsmap})
which is directly related to the observation that $A_{a}(k,l,\tau
_{a},\epsilon )$ of (\ref{eq: Akldef}) is a modular form of weight $k+l$ for 
$k+l>2$, whereas $A_{a}(1,1,\tau _{a},\epsilon )=\epsilon E_{2}(\tau _{a})$
is a quasi-modular form. The exceptional modular transformation
property of the latter term (\ref{gammaE2}) leads via Theorem \ref{Theorem_period_eps} to

\begin{theorem}
\label{TheoremGequiv} $F^{\epsilon }$ is equivariant with respect to the
action of $G$ i.e. there is a commutative diagram for $\gamma \in G$, 
\begin{equation*}
\begin{array}{ccc}
\mathcal{D}^{\epsilon } & \overset{F^{\epsilon }}{\rightarrow } & \mathfrak{H}%
_{2} \\ 
\gamma \downarrow &  & \downarrow \gamma \\ 
\mathcal{D}^{\epsilon } & \overset{F^{\epsilon }}{\rightarrow } & \mathfrak{H}%
_{2}%
\end{array}%
\end{equation*}
\end{theorem}

\bigskip 
\begin{exercise}
\label{Exercise Omega}Show that to $O(\epsilon ^{4})$%
\begin{align*}
2\pi i\Omega _{11}& =2\pi i\tau _{1}+E_{{2}}(\tau _{2}){\epsilon }^{2}+E_{{2}%
}(\tau _{1})E_{{2}}(\tau _{2})^{2}{\epsilon }^{4}, \\
2\pi i\Omega _{22}& =2\pi i\tau _{2}+E_{{2}}(\tau _{1}){\epsilon }^{2}+E_{{2}%
}(\tau _{1})^{2}E_{{2}}(\tau _{1})^{2}{\epsilon }^{4}, \\
2\pi i\Omega _{12}& =-\epsilon +E_{{2}}(\tau _{1})E_{{2}}(\tau _{2}){%
\epsilon }^{3}.
\end{align*}
\end{exercise}

\bigskip

\subsection{The Genus Two Partition Function for the Heisenberg VOA}
\label{Subsect_Genus Two}

In this section we define and compute the genus two partition function for
the Heisenberg VOA $M_0$ on the genus two Riemann surface $\mathcal{S}^{(2)}$
described in the last section. The partition function is defined in terms of
the genus one 1-point functions $Z_{M_0}^{(1)}(v,\tau _{a})$ on $\mathcal{S}%
_{a}=\mathbb{C}/\Lambda _{\tau _{a}}$ for all $v\in M$. The rationale behind
this definition, which is strongly influenced by ideas in CFT, can be
motivated by considering the following trivial sewing of a torus $\mathcal{S}%
_{1}=\mathbb{C}/\Lambda _{\tau _{1}}$ to a Riemann sphere $\mathbb{CP}^{1}$.
Let $z_{1}\in \mathcal{S}_{1}$ and $z_{2}\in \mathbb{CP}^{1}$ be local
coordinates and define the sewing by identifying the annuli $r_{a}\geq
|z_{a}|\geq |\epsilon |r_{\bar{a}}^{-1}$ via the sewing relation $%
z_{1}z_{2}=\epsilon $ (adopting the same notation as above). The resulting
surface is a torus\ described by the same modular parameter $\tau _{1}$.
\bigskip

Let $V$ be a VOA with LiZ metric $\langle \ ,\rangle $ and consider an 
$n$-point function\footnote{Here and below we include a superscript 
$(1)$ to indicate the genus of the
Riemann torus.} $F_{V}^{(1)}((v^{1},x_{1}),\ldots (v^{n},x_{n}),\tau _{1})$
for $x_{i}\in \mathcal{A}_{1}$, the torus annulus (\ref{eq: Annulus}). This
can be expressed in terms of a $1$-point function (\cite{MT1}, Lemma 3.1) by 
\begin{eqnarray}
F_{V}^{(1)}((v^{1},x_{1}),\ldots (v^{n},x_{n}),\tau _{1})
&=&Z_{V}^{(1)}(Y[v^{1},x_{1}]\ldots Y[v^{n},x_{n}]\mathbf{1},\tau _{1})
\label{eq: Z1Ysq2} \\
&=&Z_{V}^{(1)}(Y[v^{1},x_{1n}]\ldots Y[v^{n-1},x_{n-1n}]v^{n},\tau _{1}), 
\notag \\
&&  \label{eq: Z1Ysq21n}
\end{eqnarray}%
for $x_{in}=x_{i}-x_{n}$ (see (\ref{eq: Gn1pt})). Denote the
square bracket LiZ metric by $\langle \ ,\rangle _{\mathrm{sq}}$, and choose
a basis $\{u\}$ of $V_{[r]}$ with dual basis $\{\bar{u}\}$ with respect to $%
\langle \ ,\rangle _{\mathrm{sq}}$. Expanding in this basis we find that for
any $0\leq k\leq n-1$ 

\begin{equation*}
Y[v^{k+1},x_{k+1}]\ldots Y[v^{n},x_{n}]\mathbf{1=}\sum_{r\geq 0}\sum_{u\in
V_{[r]}}\langle \bar{u},Y[v^{k+1},x_{k+1}]\ldots Y[v^{n},x_{n}]\mathbf{1}%
\rangle _{\mathrm{sq}}u,
\end{equation*}%
so that 
\begin{eqnarray*}
F_{V}^{(1)}((v^{1},x_{1}),\ldots (v^{n},x_{n}),\tau _{1}) &=&\sum_{r\geq
0}\sum_{u\in V_{[r]}}Z_{V}^{(1)}(Y[v^{1},x_{1}]\ldots Y[v^{k},x_{k}]u,\tau
_{1}) \\
&&.\langle \bar{u},Y[v^{k+1},x_{k+1}]\ldots Y[v^{n},x_{n}]\mathbf{1}\rangle
_{\mathrm{sq}}.
\end{eqnarray*}%
Using (\ref{eq: Z1Ysq21n}) we have 
\begin{equation*}
Z_{V}^{(1)}(Y[v^{1},x_{1}]\ldots Y[v^{k},x_{k}]u,\tau _{1})=\mathrm{Res}%
_{z_{1}}z_{1}^{-1}F_{V}^{(1)}((v^{1},x_{1}),\ldots
(v^{k},x_{k}),(u,z_{1}),\tau _{1}).
\end{equation*}%
Let us now assume that each $v^{i}$ is quasi-primary of $L[0]$ weight $wt[v^{i}]$ and let $y_{i}=\epsilon/x_{i}\in \mathbb{CP}^{1}$. Then (\ref{eq: inv bil form}), (\ref{eq: LiZsym}%
), (\ref{eq: Y_D}) and (\ref{eq: adj op qp}) respectively imply 
\begin{eqnarray*}
&&\langle \bar{u},Y[v^{k+1},x_{k+1}]\ldots Y[v^{n},x_{n}]\mathbf{1}\rangle _{%
\mathrm{sq}} \\
&=&\langle \mathbf{1},Y^{\dagger }[v^{n},x_{n}]\ldots Y^{\dagger
}[v^{k+1},x_{k+1}]\bar{u}\rangle _{\mathrm{sq}} \\
&=&\langle \mathbf{1},\epsilon ^{L[0]}Y^{\dagger }[v^{n},x_{n}]\epsilon
^{-L[0]}\ldots \epsilon ^{L[0]}Y^{\dagger }[v^{k+1},x_{k+1}]\epsilon
^{-L[0]}\epsilon ^{L[0]}\bar{u}\rangle _{\mathrm{sq}} \\
&=&\epsilon ^{r}\langle \mathbf{1},Y[v^{n},y_{n}]\ldots Y[v^{k+1},y_{k+1}]%
\bar{u}\rangle _{\mathrm{sq}}\prod_{k+1\leq j\leq n}(-\frac{\epsilon }{%
x_{j}^{2}})^{wt[v^{j}]} \\
&=&\epsilon ^{r}\mathrm{Res}_{z_{2}}z_{2}^{-1}Z_{V}^{(0)}((v^{n},y_{n}),%
\ldots (v^{k+1},y_{k+1}),(\bar{u},z_{2}))\prod_{k+1\leq j\leq n}(\frac{dy_{j}%
}{dx_{j}})^{wt[v^{j}]}.
\end{eqnarray*}%
Note that we are also making use here of the isomorphism between the round
and square bracket formalisms in the identification of the genus zero
correlation function. The result of these calculations is that for any $%
0\leq k\leq n-1$ 
\begin{eqnarray}
&&\mathcal{F}_{V}^{(1)}(v^{1},\ldots v^{n};\tau _{1})\equiv
F_{V}^{(1)}((v^{1},x_{1}),\ldots (v^{n},x_{n}),\tau _{1})\prod_{1\leq i\leq
n}dx_{i}{}^{wt[v^{i}]} = \notag \\
&&\sum_{r\geq 0}\epsilon ^{r}\sum_{u\in V_{[r]}}\mathrm{Res}%
_{z_{1}}z_{1}^{-1}F_{V}^{(1)}((v^{1},x_{1}),\ldots
(v^{k},x_{k}),(u,z_{1}),\tau _{1}).  \notag \\
&&\mathrm{Res}_{z_{2}}z_{2}^{-1}F_{V}^{(0)}((v^{k+1},y_{k+1}),\ldots
(v^{1},y_{1}),(\bar{u},z_{2}))\prod_{1\leq i\leq
k}dx_{i}{}^{wt[v^{i}]}\prod_{k+1\leq j\leq n}dy_{j}{}^{wt[v^{j}]}.  \notag \\
&&  \label{eq: Fgenus1to0}
\end{eqnarray}%
Following Exercises \ref{Exercise SL2C} and \ref{Exercise F0form} the
(formal) form $\mathcal{F}_{V}^{(1)}(v^{1},\ldots v^{n};\tau _{1})$ is
invariant with respect to M\"{o}bius transformations. (Similarly to Remark %
\ref{Remark Fcal0form} we note that $\mathcal{F}_{V}^{(1)}(v^{1},\ldots
v^{n};\tau _{1})$ is a conformally invariant global form on $\mathcal{S}_{1}$
for primary $v^{1},\ldots v^{n}$). Geometrically, (\ref{eq: Fgenus1to0}) is
telling us that we express $\mathcal{F}_{V}^{(1)}(v^{1},\ldots v^{n};\tau
_{1})$ via the sewing procedure in terms of data arising from $\mathcal{F}%
_{V}^{(1)}(v^{1},\ldots v^{k},u;\tau _{1})$ and $\mathcal{F}%
_{V}^{(0)}(v^{k+1},\ldots v^{1},\bar{u})$ \ (cf. (\ref{eq: Fcal0form})).
Furthermore, we may choose to consider the contribution from a quasi-primary
vector $v^{i}$ as arising from either an "insertion" at $x_{i}\in \mathcal{S}%
_{1}$ or at the identified point $y_{i}=\epsilon /x_{i}\in \mathbb{CP}^{1}$.

\begin{center}
\begin{picture}(300,100)

\put(50,50){\qbezier(10,-20)(-30,0)(10,20)}
\put(50,52){\qbezier(10,18)(50,35)(90,18)}

\put(50,48){\qbezier(10,-18)(50,-35)(90,-18)}

\put(45,50){\qbezier(25,0)(45,17)(60,0)}
\put(45,50){\qbezier(20,2)(45,-17)(65,2)}


\put(175,52){\qbezier(10,18)(50,35)(90,18)}
\put(175,50){\qbezier(90,20)(130,0)(90,-20)}
\put(175,48){\qbezier(10,-18)(50,-35)(90,-18)}


\put(140,50){\circle{16}}
\put(140,50){\circle{40}}

\put(52,65){\qbezier(90,0)(100,8)(110,16)}%
\put(165,87){\makebox(0,0){$v_i$}}
\put(136,63){\makebox(0,0){$x_i$}}



\put(25,50){\makebox(0,0){$\mathcal{S}$}}


\put(185,50){\circle{16}}
\put(185,50){\circle{40}}

\put(95,65){\qbezier(90,0)(80,8)(70,16)}%
\put(191,63){\makebox(0,0){$y_i$}}



\put(300,50){\makebox(0,0){$\mathbb{CP}^1$}}

\end{picture}

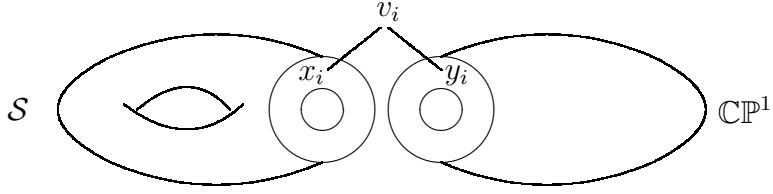
Fig.~3. Equivalent insertion of $v^{i}$ at $x_{i}$ or $y_{i}=\epsilon /x_{i}$%
.
\end{center}

A special case of (\ref{eq: Fgenus1to0}) is the partition (0-point) function
for which we find the trivial identity 
\begin{equation}
Z_{V}^{(1)}(\tau _{1})=\sum_{r\geq 0}\epsilon ^{r}\sum_{u\in
V_{[r]}}Z_{V}^{(1)}(u,\tau _{1})\mathrm{Res}_{z_{2}}z_{2}^{-1}F_{V}^{(0)}(%
\bar{u},z_{2})=Z_{V}^{(1)}(\tau _{1})+0,  \label{eq: Z1trivial}
\end{equation}%
since $F_{V}^{(0)}(\bar{u},z_{2})=0$ for $\bar{u}\notin V_{[0]}$.

\bigskip 

Motivated by this example, we define the genus two partition
function where we effectively replace the Riemann sphere on the rhs in Fig.~3 by a second torus $\mathcal{S}_{2}=\mathbb{C}/\Lambda _{\tau _{2}}$ as
described in the Section \ref{subsect_epsilon}. Thus replacing the genus
zero 1-point function $F_{V}^{(0)}(\bar{u},0)$ of (\ref{eq: Z1trivial}) by $%
Z_{V}^{(1)}(\bar{u},\tau _{2})$ we define the genus two partition function
for a VOA $V$ with a LiZ metric by 
\begin{equation}
Z_{V}^{(2)}(\tau _{1},\tau _{2},\epsilon )=\sum_{r\geq 0}\epsilon
^{r}\sum_{u\in V_{[r]}}Z_{V}^{(1)}(u,\tau _{1})Z_{V}^{(1)}(\bar{u},\tau
_{2}).  \label{eq: Z2defn}
\end{equation}%
The inner sum is taken over any basis $\{u\}$ for $V_{[r]}$ with dual basis $%
\{\bar{u}\}$ with respect to the square bracket LiZ metric. Although the
definition is associated with the specific genus two sewing scheme, it is
regarded at this stage as a purely formal expression which can be computed
to any given order in $\epsilon $. One can also define genus two correlation
functions by inserting appropriate genus one correlation functions in (\ref%
{eq: Z2defn}). We do not consider these here.

\bigskip

Let us now compute the genus two partition function for the rank one
Heisenberg VOA $M_0$ generated by $a$ of weight 1. We employ the square
bracket Fock basis of (\ref{eq: vFovkki}) which we alternatively notate here
(cf. (\ref{eq: Fockstate})) by%
\begin{equation}
v=v(\lambda )=a[-1]^{e_{1}}\ldots a[-p]^{e_{p}}\mathbf{1},
\label{eq: vlambda}
\end{equation}%
for non-negative integers $e_{i}$. We recall that $v(\lambda )$ is of square
bracket weight $wt[v]=\sum_{i}ie_{i}$ and is described by a label set $\Phi
_{\lambda }=\{1,\ldots ,p\}$ with $n=$ $\sum e_{i}$ elements corresponding
to an unrestricted partition $\lambda =\{1^{e_{1}}\ldots p^{e_{p}}\}$ of $%
wt[v]$. The Fock vectors (\ref{eq: vlambda}) form a diagonal basis for the
LiZ metric $\langle \ ,\rangle _{\mathrm{sq}}$ with 
\begin{equation}
\bar{v}=\frac{1}{\prod_{1\leq i\leq p}(-i)^{e_{i}}e_{i}!}v,  \label{eq: vsq2}
\end{equation}%
from (\ref{eq: inner prod}). Following (\ref{eq: Z2defn}), we find 
\begin{equation}
Z_{M_0}^{(2)}(\tau _{1},\tau _{2},\epsilon )=\sum_{v\in V}\frac{\epsilon
^{wt[v]}}{\prod_{i}(-i)^{e_{i}}e_{i}!}Z_{M_0}^{(1)}(v,\tau
_{1})Z_{M_0}^{(1)}(v,\tau _{2}),  \label{eq: Z2_def_eps}
\end{equation}%
where the sum is taken over the basis (\ref{eq: vlambda}). $Z_{M_0}^{(2)}(\tau _{1},\tau _{2},\epsilon )$ is given by the
following closed formula \cite{MT4}:

\begin{theorem}
\label{Theorem_Z2_boson} The genus two partition function for the rank one
Heisenberg VOA is%
\begin{equation}
Z_{M_0}^{(2)}(\tau _{1},\tau _{2},\epsilon )=\frac{1}{\eta (\tau _{1})\eta
(\tau _{2})}(\det (I-A_{1}A_{2}))^{-1/2},  \label{eq: Z2_1bos}
\end{equation}%
with $A_{a}$ of (\ref{eq: Akldef}).
\end{theorem}

\noindent\textbf{Proof.} The proof relies on an
interesting graph-theoretic interpretation of (\ref{eq: Z2_def_eps}). This
follows the technique introduced in Theorem \ref{theorem: HeisenGraph} for
graphically interpreting the genus one 1-point function $Z_{M_0}^{(1)}(v(%
\lambda ),\tau _{1})$ in terms the sum of weights for the $\phi $-graphs. We
sketch the main features of the proof leaving the interested reader to
explore the details in \cite{MT4}.
\medskip

Since $v(\lambda )$ is indexed by unrestricted partitions $\lambda
=\{1^{e_{1}}, 2^{e_{2}}, \ldots \}$ we may write (\ref{eq: Z2_def_eps}) as 
\begin{equation}
Z_{M_0}^{(2)}(\tau _{1},\tau _{2},\epsilon )=\sum_{\lambda =\{i^{e_{i}}\}}%
\frac{1}{\prod_{i}e_{i}!}.\prod_{i}\left( \frac{\epsilon ^{i}}{-i}\right)
^{e_{i}}Z_{M_0}^{(1)}(v(\lambda ),\tau _{1})Z_{M_0}^{(1)}(v(\lambda ),\tau _{2}).
\label{eq: part func diag}
\end{equation}%
Theorem \ref{theorem: HeisenGraph} implies $Z_{M_0}^{(1)}(v(\lambda ),\tau
_{1})=0$ for odd $n=\sum e_{i}$ whereas for $n$ even 
\begin{equation*}
Z_{M_0}^{(1)}(v(\lambda ),\tau _{1})Z_{M_0}^{(1)}(v(\lambda ),\tau _{2})=\frac{1%
}{\eta (\tau _{1})\eta (\tau _{2})}\sum_{\gamma _{\phi _{1}}}\sum_{\gamma
_{\phi _{2}}}\kappa (\gamma _{\phi _{1}},\tau _{1})\kappa (\gamma _{\phi
_{2}},\tau _{2}),
\end{equation*}%
where $\gamma _{\phi _{1}},\gamma _{\phi _{2}}$ independently range over the 
$\phi -$graphs for $\Phi _{\lambda }$. Any pair $\gamma _{\phi _{1}},\gamma
_{\phi _{2}}$ can be naturally combined to form a \emph{chequered diagram} $%
D $ consisting of $n$ vertices labelled by $\Phi _{\lambda }$ of valence 2
with $n$ unoriented edges $\overset{k}{\bullet }\overset{a}{-}\overset{l}{%
\bullet }$ consecutively labelled by $a=1,2$ as specified by $\phi
_{a}=\ldots (kl)\ldots $. Following Exercise \ref{Exercise F(Phi)} there are 
$(n!!)^{2}$ chequered diagrams for a given $v(\lambda )$. We illustrate an
example of such a diagram in Fig.~4 for $v=a[-1]^{3}a[-2]^{2}a[-5]\mathbf{1}$
with $\phi _{1}$ of Fig.~1 and a separate choice for $\phi _{2}$ with cycle
shape $(11)(22)(15)$

\begin{center}
\begin{picture}(250,80)


\put(100,50){\line(1,2){10}}
\put(90,50){\makebox(0,0){1}}
\put(100,50){\circle*{4}}
\put(102,63){\makebox(0,0){\scriptsize 1}}

\put(110,70){\line(1,-2){20}}
\put(105,78){\makebox(0,0){1}}
\put(110,70){\circle*{4}}
\put(117,47){\makebox(0,0){\scriptsize 2}}

\put(135,60){\qbezier(5,-10)(-8,0)(-5,10)}%

\put(135,60){\qbezier(5,-10)(8,0)(-5,10)}%

\put(139,78){\makebox(0,0){2}}
\put(130,70){\circle*{4}}
\put(141,65){\makebox(0,0){\scriptsize 1}}

\put(130,55){\makebox(0,0){\scriptsize 2}}

\put(150,50){\makebox(0,0){2}}
\put(140,50){\circle*{4}}

\put(110,30){\line(1,0){20}}
\put(135,20){\makebox(0,0){5}}
\put(130,30){\circle*{4}}
\put(120,25){\makebox(0,0){\scriptsize 1}}

\put(100,50){\line(1,-2){10}}
\put(105,20){\makebox(0,0){1}}
\put(110,30){\circle*{4}}
\put(102,36){\makebox(0,0){\scriptsize 2}}

\end{picture}

Fig.~4 A Chequered Diagram
\end{center}

For $\lambda =\{1^{e_{1}}\ldots p^{e_{p}}\}$ the symmetric group $\Sigma
(\Phi _{\lambda })$ acts on the chequered diagrams which have $\Phi
_{\lambda }$ as underlying set of labeled nodes. We define $\mathrm{Aut}(D)$%
, the \emph{automorphism group of }$D$, \ to be the subgroup of $\Sigma
(\Phi _{\lambda })$ which preserves node labels. $\mathrm{Aut}(D)$ is
isomorphic to $\Sigma _{e_{1}}\times \ldots \times \Sigma _{e_{p}}$ of order 
$|\mathrm{Aut}(D)|=\prod_{i}e_{i}!$. We may thus express (\ref{eq: part func
diag}) as a sum over the isomorphism classes of chequered diagrams $D$ with%
\begin{equation*}
Z_{M_0}^{(2)}(\tau _{1},\tau _{2},\epsilon )=\frac{1}{\eta (\tau _{1})\eta
(\tau _{2})}\sum_{D}\frac{\zeta (D)}{|\mathrm{Aut}(D)|},
\end{equation*}%
and 
\begin{equation}
\zeta (D)=\prod_{i}\left( \frac{\epsilon ^{i}}{i}\right) ^{e_{i}}\kappa
(\gamma _{\phi _{1}},\tau _{1})\kappa (\gamma _{\phi _{2}},\tau _{2}),
\label{eq: zetaD}
\end{equation}%
where $D$ is determined by $\gamma _{\phi _{1}},\gamma _{\phi _{2}}$ and
noting that $\prod_{i}\left( -1\right) ^{e_{i}}=1$ for $n$ even. From (\ref%
{eq: kappagam}) we recall that $\kappa (\gamma _{\phi _{a}},\tau _{1})$ is a
product of the weights of the $a$ labelled edges. Then $\zeta (D)$ can be
more naturally expressed in terms of a weight function on chequered diagrams
defined by 
\begin{equation}
\zeta (D)=\Pi _{E}\zeta (E),  \label{eq: zetaDprod}
\end{equation}%
where the product is taken over the edges $E$ of $D$ and where for an edge $%
E $ labeled \ $\overset{k}{\bullet }\overset{a}{-}\overset{l}{\bullet }$ \
we define 
\begin{equation*}
\zeta (E)=\frac{\epsilon ^{\frac{k+l}{2}}}{\sqrt{kl}}C(k,l,\tau
_{a})=A_{a}(k,l,\tau _{a},\epsilon ),
\end{equation*}%
for $A_{a}$ of (\ref{eq: Akldef}).

\bigskip

Every chequered diagram can be formally represented as a product 
\begin{equation*}
D=\prod_{i}L_{i}^{m_{i}},
\end{equation*}%
with $D$ a disjoint union of unoriented chequered cycles (connected
diagrams) $L_{i}$ with multiplicity $m_{i}$ (e.g. the chequered diagram of
Fig.~4 is the product of two disjoint cycles). Then $\mathrm{Aut}(D)$ is
isomorphic to the direct product of the groups $\mathrm{Aut}(L_{i}^{m_{i}})$
of order $\left\vert \mathrm{Aut}(L_{i}^{m_{i}})\right\vert =\left\vert 
\mathrm{Aut}(L_{i})\right\vert ^{m_{i}}m_{i}!$ so that 
\begin{equation*}
|\mathrm{Aut}(D)|=\prod_{i}\left\vert \mathrm{Aut}(L_{i})\right\vert
^{m_{i}}m_{i}!.
\end{equation*}%
But from (\ref{eq: zetaDprod}) it is clear that $\zeta (D)$ is
multiplicative over disjoint unions of diagrams, and we find

\begin{eqnarray*}
\sum_{D}\frac{\zeta (D)}{|\mathrm{Aut}(D)|} &=&\prod_{L}\sum_{m\geq 0}\frac{%
\zeta (L)^{m}}{|\mathrm{Aut}(L)|^{m}m!} \\
&=&\exp \left( \sum_{L}\frac{\zeta (L)}{|\mathrm{Aut}(L)|}\right) ,
\end{eqnarray*}%
where $L$ ranges over isomorphism classes of unoriented chequered cycles.
Further analysis shows that \cite{MT4} 
\begin{equation*}
\sum_{L}\frac{\zeta (L)}{|\mathrm{Aut}(L)|}=\frac{1}{2}\mathrm{Tr}%
(\sum_{n\geq 1}\frac{1}{n}(A_{1}A_{2})^{n})=-\frac{1}{2}\mathrm{Tr}(\log
(1-A_{1}A_{2})),
\end{equation*}%
so that we find 
\begin{equation*}
\sum_{D}\frac{\zeta (D)}{|\mathrm{Aut}(D)|}=(\det (1-A_{1}A_{2}))^{-1/2},
\end{equation*}%
following (\ref{eq: logdet}). Thus  Theorem \ref{Theorem_Z2_boson} holds.

\bigskip

The convergence and holomorphy of the determinant is the subject of Theorem %
\ref{Theorem_A1A2} (b) so that having computed the closed formula (\ref{eq:
Z2_1bos}) we may conclude that $Z_{M_0}^{(2)}(\tau _{1},\tau _{2},\epsilon )$
is not just a formal function but can be evaluated on $\mathcal{D}^{\epsilon
}$ to find

\begin{theorem}
\label{Theorem_Z2_boson_eps_hol} $Z_{M_0}^{(2)}(\tau _{1},\tau _{2},\epsilon )$
is holomorphic on\ the domain $\mathcal{D}^{\epsilon }$.
\end{theorem}

\bigskip

We next consider the automorphic properties of $Z_{M_0}^{(2)}(\tau _{1},\tau
_{2},\epsilon )$ with respect to the modular group $G\subset Sp(4,\mathbb{Z})
$ of (\ref{eq: GDeps}) which acts on $\mathcal{D}^{\epsilon }$. We first
recall a little from the classical theory of modular forms (cf. Section \ref{Section_ModForms}).
For a meromorphic function $f(\tau )$ on $\mathfrak{H}$, $k\in \mathbb{Z}$
and $\gamma =\left( 
\begin{array}{cc}
a & b \\ 
c & d%
\end{array}%
\right) \in SL(2,\mathbb{Z})$, we define the right action 
\begin{equation}
f(\tau )|_{k}\gamma =f(\gamma \tau )\ (c\tau +d)^{-k},
\label{eq: slashaction}
\end{equation}%
where, as usual 
\begin{equation*}
\gamma \tau =\frac{a\tau +b}{c\tau +d}.
\end{equation*}%
$f(\tau )$ is called a weak modular form for a subgroup $\Gamma \subseteq
SL(2,\mathbb{Z})$ of weight $k$ if $f(\tau )|_{k}\gamma =f(\tau )$ for all $%
\gamma \in \Gamma $.

\medskip

We have already discussed the (genus one) partition function for the rank $n$ Heisenberg VOA
$V=M_0^{\otimes n}$ in Subsection \ref{Subsect_gradeddim} (cf. (\ref{M0dgdim})). 
In particular, for $n=2$ we have
\begin{equation*}
Z_{M_0^{2}}^{(1)}(\tau )=Z_{M_0}^{(1)}(\tau )^{2}=\frac{1}{\eta (\tau )^{2}}.
\end{equation*}%
Then we find 
\begin{equation}
Z_{M_0^{2}}^{(1)}(\tau )|_{-1}\gamma =\chi (\gamma )Z_{M_0^{2}}^{(1)}(\tau ),
\label{eq: Z1modgam}
\end{equation}%
where $\chi$ is a character of $SL(2,\mathbb{Z})$ of order $12$ (cf. Exercise \ref{Exercise_8.1.4} and \cite{Se}), and
\begin{equation}
Z_{M_0^{24}}^{(1)}(\tau )^{-1}=\Delta (\tau ).  \label{eq: Z24genus1}
\end{equation}

\bigskip
Similarly, we consider the genus two partition function for the rank two
Heisenberg VOA given by 
\begin{equation}
Z_{M_0^{2}}^{(2)}(\tau _{1},\tau _{2},\epsilon )=Z_{M_0}^{(2)}(\tau _{1},\tau
_{2},\epsilon )^{2}=\frac{1}{\eta (\tau _{1})^{2}\eta (\tau _{2})^{2}\det
(I-A_{1}A_{2})}.  \label{eq: Z2mod_eps}
\end{equation}%
Analogously to (\ref{eq: slashaction}), we define for all $\gamma \in G$ 
\begin{equation}
f(\tau _{1},\tau _{2},\epsilon )|_{k}\gamma =f(\gamma (\tau _{1},\tau
_{2},\epsilon ))\det (C\Omega +D)^{-k},  \label{eq: Gaction}
\end{equation}%
where the action of $\gamma $ on the right-hand-side is as in (\ref{eq:
GDeps}) and $\Omega (\tau _{1},\tau _{2},\epsilon )$ is determined by
Theorem \ref{Theorem_period_eps}. Then (\ref{eq: Gaction}) defines a right
action of $G$ on functions $f(\tau _{1},\tau _{2},\epsilon )$. We next
obtain a natural genus two extension of (\ref{eq: Z1modgam}). Define the a
character $\chi ^{(2)}$ of $G$ by 
\begin{equation*}
\chi ^{(2)}(\gamma _{1}\gamma _{2}\beta ^{m})=(-1)^{m}\chi (\gamma _{1})\chi
(\gamma _{2}),\quad \ \gamma _{i}\in \Gamma _{i},\ i=1,2.
\end{equation*}%
with $\Gamma _{i},\beta $ of (\ref{eq: G1G2}) and (\ref{eq: beta}). $\chi
^{(2)}$ takes values which are twelfth roots of unity. Then, much as for
Theorem \ref{Theorem_period_eps}, the exceptional transformation law
of $A_{a}(1,1,\tau _{a},\epsilon )=E_{2}(\tau _{a})$ implies that

\begin{theorem}
\label{Theorem_Z2_G}If $\gamma \in G$ then 
\begin{equation*}
Z_{M_0^{2}}^{(2)}(\tau _{1},\tau _{2},\epsilon )|_{-1}\gamma =\chi
^{(2)}(\gamma )Z_{M_0^{2}}^{(2)}(\tau _{1},\tau _{2},\epsilon ).
\end{equation*}
\end{theorem}

\bigskip

The definition (\ref{eq: Gaction}) is analogous to that for a \emph{Siegel
modular form} for the symplectic group $Sp(4,\mathbb{Z})$ defined as follows
(e.g. \cite{Fr}). For a meromorphic function $F(\Omega )$ on $\mathfrak{H}_{2}$, 
$k\in \mathbb{Z}$ and $\gamma \in Sp(4,\mathbb{Z})$, we define the right
action 
\begin{equation}
F(\Omega )|_{k}\gamma =F(\gamma .\Omega )\ \det (C\Omega +D)^{-k},
\label{eq: FOmega}
\end{equation}%
with $\gamma .\Omega $ of (\ref{eq: Sp4zOmega}). $F(\Omega )$ is called a
modular form for $\Gamma \subseteq Sp(4,\mathbb{Z})$ of weight $k$ if $%
F(\Omega )|_{k}\gamma =F(\Omega )$ for all $\gamma \in \Gamma $.

\medskip

Theorem \ref{Theorem_Z2_G} implies that for the rank 24 Heisenberg VOA 
$M_0^{24}$ 
\begin{equation}
Z_{M_0^{24}}^{(2)}(\tau _{1},\tau _{2},\epsilon )|_{-12}\gamma
=Z_{M_0^{24}}^{(2)}(\tau _{1},\tau _{2},\epsilon ),  \label{eq: Z24}
\end{equation}%
for all $\gamma \in G$. This might lead one to speculate that, analogously
to the genus one case in (\ref{eq: Z24genus1}), $Z_{M_0^{24}}^{(2)}(\tau
_{1},\tau _{2},\epsilon )^{-1}$ is a holomorphic Siegel modular form of
weight 12. Indeed, there does exist a unique holomorphic Siegel 12 form, $%
\Delta _{^{12}}^{(2)}(\Omega )$, such that 
\begin{equation*}
\Delta _{^{12}}^{(2)}(\Omega )\rightarrow \Delta (\tau _{1})\Delta (\tau
_{2}),
\end{equation*}%
as $\epsilon \rightarrow 0$ but explicit calculations show that $%
Z_{M_0^{24}}^{(2)}(\tau _{1},\tau _{2},\epsilon )^{-1}\neq $ $\Delta
_{^{12}}^{(2)}(\Omega )$. In any case, we cannot naturally extend the action
of $G$ on $\mathcal{D}^{\epsilon }$ to $Sp(4,\mathbb{Z})$. These
observations are strongly expected to be related to the conformal anomaly 
\cite{BK} in string theory and to the non-existence of a\ global section for
the Hodge line bundle in algebraic geometry \cite{Mu2}.

\bigskip

Siegel modular forms do arise in the determination of the genus two
partition function for a lattice VOA $V_{L}$ for even lattice $L$ of rank $l$
(and conjecturally for all rational theories) as follows. We recall the
genus one partition function for $V_{L}$ is (cf. Subsection \ref{Subsect_LatticeVOA})  
\begin{equation}
Z_{V_{L}}^{(1)}(\tau )=Z_{M_0^{l}}^{(1)}(\tau )\theta _{L}^{(1)}(\tau ),
\label{eq: Zlattice1}
\end{equation}%
\newline
for $\theta _{L}^{(1)}(\tau )=\sum_{\alpha }q^{(\alpha ,\alpha )/2}$. In the
genus two case, we may define the Siegel lattice theta function by \cite{Fr}%
\begin{equation*}
\theta _{L}^{(2)}(\Omega )=\sum_{\alpha ,\beta \in L}\exp (\pi i((\alpha
,\alpha )\Omega _{11}+2(\alpha ,\beta )\Omega _{12}+(\beta ,\beta )\Omega
_{22})).
\end{equation*}%
$\theta _{L}^{(2)}(\Omega )$ is a Siegel modular form of weight $l/2$ for a
subgroup of $Sp(4,\mathbb{Z})$. The genus one result (\ref{eq: Zlattice1})
is naturally generalized to find \cite{MT4}:

\begin{theorem}
\label{Theorem_Z2_L_eps}For a lattice VOA $V_{L}$ we have 
\begin{equation*}
Z_{V_{L}}^{(2)}(\tau _{1},\tau _{2},\epsilon )=Z_{M_0^{l}}^{(2)}(\tau
_{1},\tau _{2},\epsilon )\theta _{L}^{(2)}(\Omega ).
\end{equation*}
\end{theorem}

\bigskip

\begin{exercise}
\label{Exercise ZB}Show that $Z_{M_0}^{(2)}(\tau _{1},\tau _{2},\epsilon )$ to 
$O({\epsilon }^{4})$ is given by%
\begin{equation*}
\frac{1}{\eta (\tau _{1})\eta (\tau _{2})}\left[ 1+\frac{1}{2}E_{2}(\tau
_{1})E_{2}(\tau _{2}){\epsilon }^{2}+\left( \frac{3}{8}\,E_{2}(\tau
_{1})^{2}E_{2}(\tau _{2})^{2}+\frac{15}{2}\,E_{4}(\tau _{1})E_{4}(\tau
_{2})\right) {\epsilon }^{4}\right] .
\end{equation*}
\end{exercise}

\begin{exercise}
\label{Exercise ZB det}Verify (\ref{eq: Z2_1bos}) to $O({\epsilon }^{4})$ by
showing that 
\begin{equation*}
\det (I-A_{1}A_{2})=1-E_{2}(\tau _{1})E_{2}(\tau _{2}){\epsilon }%
^{2}-15\,E_{4}(\tau _{1})E_{4}(\tau _{2}){\epsilon }^{4}+O({\epsilon }^{6}).
\end{equation*}
\end{exercise}

\bigskip

\section{Exceptional VOAs and the Virasoro Algebra}
\label{Section_Exceptional}

In this section we review some recent research concerning a r\^{o}le played by
the Virasoro algebra in certain exceptional VOAs \cite{T2}, \cite{T3}. We
will mainly concern ourselves here with simple VOAs $V$ of strong CFT-type
for which $\dim V_{1}>0$. We construct certain quadratic Casimir vectors
from the elements of $V_{1}$ and examine the constraints on $V$ arising from
the assumption that the Casimir vectors of low weight are Virasoro
descendants of the vacuum. This sort of assumption is similar to that of  `minimality'
in the holomorphic VOAs $V^{(k)}$ that we discussed in Subsection \ref{Subsect_Applications}.
 In particular we discuss how a special set of
simple Lie algebras: $A_{1},A_{2},G_{2},D_{4},F_{4},E_{6},E_{7},E_{8}$,
known as Deligne's exceptional series \cite{De}, arises in this context. We
also show that the genus one partition function is determined by the same
Virasoro condition. These constraints follow from an analysis of appropriate
genus zero matrix elements and genus one $2$-point functions. In particular,
we will make a relatively elementary use of rational matrix elements, the
LiZ metric, the Zhu reduction formula and modular differential equations. As
such, this example offers a useful and explicit application of many of the
concepts reviewed in these Notes.

\subsection{Quadratic Casimirs and Genus Zero Constraints}
\label{Subsect_Casimir}

Consider a simple VOA $V$ of strong CFT-type of central charge $c$ with $d=\dim V_{1}>0$. From Theorem \ref{theorem: LiZ}, $V$ possesses an LiZ metric 
$\langle \ ,\rangle$,  i.e. a unique (non-degenerate) normalized symmetric bilinear form.
For $a,b\in V_{1}$ define $[a,b]\equiv a_{0}b$. From Exercise \ref{Exercise_9.4.1} this defines a Lie algebra on $V_1$ with invariant bilinear form $\langle \ ,\rangle $. We denote this Lie algebra by $\mathfrak{g}$. The modes
of elements of $V_{1}$ satisfy the Kac-Moody algebra (cf. Exercise \ref{Exercise KacMoody}) 
\begin{equation}
\lbrack a_{m},b_{n}]=[a,b]_{m+n}-m\langle a,b\rangle \delta _{m+n,0}.
\label{eq: KacMoody}
\end{equation}
which we denote by $\hat{\mathfrak{g}}$.  

\medskip

Let $\{u^{\alpha
}|\alpha =1\ldots d\}$ and $\{\bar{u}^{\beta }|\beta =1\ldots d\}$ denote a $%
\mathfrak{g}$-basis and LiZ dual basis respectively. Define the \emph{quadratic Casimir vectors} by 
\begin{equation}
\lambda ^{(n)}=u_{1-n}^{\alpha }\bar{u}^{\alpha }\in V_{n},
\label{eq:lambdan}
\end{equation}%
where $\alpha $ is summed. Since $u^{\alpha }\in V_{1}$ is a primary vector it follows that $%
[L_{m},u_{n}^{\alpha }]=-nu_{m+n}^{\alpha }$ and hence 
\begin{equation}
L_{m}\lambda ^{(n)}=(n-1)\lambda ^{(n-m)}\text{\ for }m>0.
\label{eq:LMlambda}
\end{equation}

\bigskip 

Let $\mbox{Vir}_{c}$ denote the subVOA of $V$ generated by the Virasoro vector $%
\omega $. We then find

\begin{lemma}
\label{lemma: VirSimple}The LiZ metric is non-degenerate on $\mbox{Vir}_{c}$.
\end{lemma}

\noindent \textbf{Proof.} Let $v=L_{-n_{1}}L_{-n_{2}}\ldots L_{-n_{k}}\mathbf{1\in }%
\mbox{Vir}_{c}$. Then (\ref{eq: LLdagger}) gives%
\begin{equation*}
\langle v,a\rangle =\langle \mathbf{1},L_{n_{k}}\ldots
L_{n_{2}}L_{n_{1}}a\rangle =0,
\end{equation*}%
for $a\in V \setminus \mbox{Vir}_{c}$. Since $\langle \ ,\rangle $ is non-degenerate on 
$V$ it must be non-degenerate on $\mbox{Vir}_{c}$. 

\begin{remark}
\label{Remark Virc}We note that this implies from Theorem \ref{theorem: LiZ}
that $\mbox{Vir}_{c}$ is simple with $c\neq c_{p,q}$ of (\ref{eq: Cpq}).
\end{remark}

\medskip
We now consider the constraints on $\mathfrak{g}$ that follow from assuming
that $\lambda ^{(n)}\in \mbox{Vir}_{c}$ for small $n$.\footnote{%
The original motivation, due to Matsuo \cite{Mat}, for considering quadratic
Casimirs is that both they and $\mbox{Vir}_{c}$ are invariant under the
automorphism group of $V$. Matsuo considered VOAs for which the automorphism
invariants of $V_{n}$ consist \textbf{only} of Virasoro descendents for
small $n$. Hence for these VOAs it necessarily follows that $\lambda
^{(n)}\in \mbox{Vir}_{c}$.} Firstly let us note \cite{Mat}

\begin{lemma}
\label{lemma:lambda} If $\lambda ^{(n)}\in \mbox{Vir}_{c}$ then $\lambda ^{(m)}\in
\mbox{Vir}_{c}$ and is uniquely determined for all $m\leq n$.
\end{lemma}

\noindent \textbf{Proof.} If $\lambda ^{(n)}\in \mbox{Vir}_{c}$ then $\lambda
^{(n)}=\sum_{v\in (\mbox{Vir}_{c})_{n}}\langle \bar{v},\lambda ^{(n)}\rangle v$
summing over a basis for $(\mbox{Vir}_{c})_{n}$. But $\langle \bar{v},\lambda
^{(n)}\rangle $ is uniquely determined by repeated use of (\ref{eq:LMlambda}%
) and Exercise \ref{Exercise lambda01}. Furthermore, for $m\leq n$ we have $%
\lambda ^{(m)}=\frac{1}{n-1}L_{n-m}\lambda ^{(n)}\in (\mbox{Vir}_{c})_{m}$.
\medskip

If follows that $\lambda ^{(2)}\in \mbox{Vir}_{c}$ implies 
\begin{equation}
\lambda ^{(2)}=-\frac{2d}{c}\omega ,  \label{eq:lambda2}
\end{equation}%
where $c\neq 0$ following Remark \ref{Remark Virc}. Note that for $\mathfrak{%
g}$ simple, this is the standard Sugawara construction for $\omega $ of
(\ref{Sugawaradef}).
Similarly  $\lambda ^{(4)}\in \mbox{Vir}_{c}$ implies
\begin{equation}
\lambda ^{(4)}=-\frac{3d}{c\left( 5c+22\right) }\left( 4L_{-2}^{2}\mathbf{1}%
+\left( 2+c\right) L_{-4}\mathbf{1}\right) ,  \label{eq:lambda4}
\end{equation}%
with $c\neq 0,-\frac{22}{5}$ following Remark \ref{Remark Virc}.

\medskip
We next consider the constraints on $\mathfrak{g}$ if either (\ref%
{eq:lambda2}) or (\ref{eq:lambda4}) hold. We do this by analysing the
following genus zero matrix element 
\begin{equation}
F(a,b;x,y)=\langle a,Y(u^{\alpha },x)Y(\bar{u}^{\alpha },y)b\rangle ,
\label{eq: Fabdef}
\end{equation}%
where $\alpha $ is summed and $a,b\in V_{1}$. Using associativity 
and (\ref{eq:lambdan}) we find 
\begin{eqnarray}
F(a,b;x,y) &=&\langle a,Y(Y(u^{\alpha },x-y)\bar{u}^{\alpha },y)b\rangle 
\notag \\
&=&\frac{1}{(x-y)^{2}}\sum_{n\geq 0}\langle a,o(\lambda ^{(n)})b\rangle
\left( \frac{x-y}{y}\right) ^{n},  \label{eq:Fablambda}
\end{eqnarray}%
where $o(\lambda ^{(n)})= \lambda _{n-1}$ from (\ref{o(v)}). Thus Exercise \ref{Exercise lambda01} implies%
\begin{equation}
F(a,b;x,y)=\frac{1}{(x-y)^{2}}\left[ -d\langle a,b\rangle +0+\langle
a,o(\lambda ^{(2)})b\rangle \left( \frac{x-y}{y}\right) ^{2}+\ldots \right] .
\label{eq: Fab03}
\end{equation}%
Alternatively, we also have 
\begin{eqnarray*}
F(a,b;x,y) &=&\langle a,Y(u^{\alpha },x)e^{yL_{-1}}Y(b,-y)\bar{u}^{\alpha
}\rangle \\
&=&\langle a,e^{yL_{-1}}Y(u^{\alpha },x-y)Y(b,-y)\bar{u}^{\alpha }\rangle \\
&=&\langle e^{yL_{1}}a,Y(u^{\alpha },x-y)Y(b,-y)\bar{u}^{\alpha }\rangle \\
&=&\langle a,Y(u^{\alpha },x-y)Y(b,-y)\bar{u}^{\alpha }\rangle \\
&=&\frac{1}{y^{2}}\sum_{m\geq 0}\langle a,u_{m-1}^{\alpha }b_{1-m}\bar{u}%
^{\alpha }\rangle \left( \frac{-y}{x-y}\right) ^{m} \\
&=&\frac{1}{y^{2}}\left[ \langle a,u_{-1}^{\alpha }b_{1}\bar{u}^{\alpha
}\rangle -\langle a,u_{0}^{\alpha }b_{0}\bar{u}^{\alpha }\rangle \frac{y}{x-y%
}+\ldots \right] ,
\end{eqnarray*}%
using skew-symmetry and translation (cf. Exercises \ref{Exercise_skew} and \ref{Exercise_translation}), invariance of the LiZ metric
and that $a$ is primary. The leading term is 
\begin{eqnarray*}
\langle a,u_{-1}^{\alpha }b_{1}\bar{u}^{\alpha }\rangle =
-\langle a,u^{\alpha }\rangle \langle b,\bar{u}^{\alpha }\rangle 
=-\langle a,b\rangle.
\end{eqnarray*}
The next to leading term is 
\begin{eqnarray*}
-\langle a,u_{0}^{\alpha }b_{0}\bar{u}^{\alpha }\rangle  = 
\langle u^{\alpha},a_{0}b_{0}\bar{u}^{\alpha }\rangle 
=K(a,b),
\end{eqnarray*}%
the Lie algebra Killing form 
\begin{equation}
K(a,b)=Tr_{\mathfrak{g}}(a_{0}b_{0}).  \label{eq: KillingFormLie}
\end{equation}%
Thus we have 
\begin{equation}
F(a,b;x,y)=\frac{1}{y^{2}}\left[ -\langle a,b\rangle +K(a,b)\frac{y}{x-y}%
+\ldots \right] .  \label{eq: Fxy_expan2}
\end{equation}

\bigskip

From Theorem \ref{theorem: RationalYY} we know that $F(a,b;x,y)$ is given by
a rational function 
\begin{equation}
F(a,b;x,y)=\frac{f(a,b;x,y)}{x^{2}y^{2}(x-y)^{2}},  \label{eq: Fabfxy}
\end{equation}%
where $f(a,b;x,y)$ is a homogeneous polynomial of degree $4$. Furthermore $%
f(a,b;x,y)$ is clearly symmetric in $x,y$ so that it may parameterized 
\begin{equation}
f(a,b;x,y)=p(a,b)x^{2}y^{2}+q(a,b)xy(x-y)^{2}+r(a,b)(x-y)^{4},
\label{eq: fxy}
\end{equation}%
for some bilinears $p(a,b)$, $q(a,b)$ and $r(a,b)$. We find

\begin{proposition}
\label{Proposition_Lie} $p(a,b),q(a,b),r(a,b)$ are given by 
\begin{eqnarray}
p(a,b) &=&-d\langle a,b\rangle ,  \label{eq: P} \\
q(a,b) &=&K(a,b)-2\langle a,b\rangle ,  \label{eq: Q} \\
r(a,b) &=&-\langle a,b\rangle .  \label{eq: R}
\end{eqnarray}
\end{proposition}

\noindent\textbf{Proof.} Expanding (\ref{eq: Fabfxy}) in $\frac{x-y}{y}$ we have 
\begin{equation}
F(a,b;x,y)=\frac{1}{(x-y)^{2}}\left[ p(a,b)+q(a,b)\left( \frac{x-y}{y}%
\right) ^{2}+\ldots \right] ,  \label{eq: Fab_gexpan1}
\end{equation}%
whereas expanding (\ref{eq: Fabfxy}) in $\frac{y}{x-y}$ gives 
\begin{equation}
F(a,b;x,y)=\frac{1}{y^{2}}\left[ r(a,b)+(-2r(a,b)+q(a,b))\frac{y}{x-y}%
+\ldots \right] .  \label{eq: Fab_gexpan2}
\end{equation}%
Comparing to (\ref{eq: Fab03}) and (\ref{eq: Fxy_expan2}) gives the result. 

\bigskip

We next show that if $\lambda ^{(2)}\in \mbox{Vir}_{c}$ then the Killing form is
proportional to the LiZ metric:

\begin{proposition}
\label{Proposition_Lie_lambda2} If $\lambda ^{(2)}\in \mbox{Vir}_{c}$ then 
\begin{equation}
K(a,b)=-2\langle a,b\rangle \left( \frac{d}{c}-1\right) ,
\label{eq: KillingForm}
\end{equation}%
so that 
\begin{equation}
f(a,b;x,y)=-\langle a,b\rangle \left[ dx^{2}y^{2}+\frac{2d}{c}%
xy(x-y)^{2}+(x-y)^{4}\right] .  \label{eq: fxycd}
\end{equation}
\end{proposition}

\noindent\textbf{Proof.} (\ref{eq:lambda2}) implies $o(\lambda ^{(2)})=-\frac{2d}{c}%
L_{0}$. Comparing the next to leading terms in (\ref{eq: Fab03}) and (\ref%
{eq: Fab_gexpan1}) we find 
\begin{equation*}
q(a,b)=\langle a,o(\lambda ^{(2)})b\rangle =-\frac{2d}{c}\langle a,b\rangle ,
\end{equation*}%
which implies the result. 

\bigskip

Since the LiZ metric is non-degenerate, it follows from Cartan's criterion
in Lie theory that $\mathfrak{g}$ is solvable for $d=c$ and is semi-simple
for $d\neq c$ i.e. 
\begin{equation*}
\mathfrak{g}=\mathfrak{g}^{1}\oplus \mathfrak{g}^{2}\oplus \ldots \oplus 
\mathfrak{g}^{r},
\end{equation*}%
for simple components $\mathfrak{g}^{i}$ of dimension $d^{i}$. The
corresponding Kac-Moody algebra $\hat{\mathfrak{g}}^{i}$ has level $l^{i}=-%
\frac{1}{2}\langle \mathbf{\alpha }^{i},\mathbf{\alpha }^{i}\rangle $ where $%
\mathbf{\alpha }^{i}$ is a long root\footnote{%
Then $(a,b)_{i}\equiv -\langle a,b\rangle /l_{i}$ is the unique
non-degenerate form on $\hat{\mathfrak{g}}_{i}^{(l_{i})}$ with normalization 
$(\mathbf{\alpha }_{i},\mathbf{\alpha }_{i})_{i}=2$.} so that the dual
Coxeter number is 
\begin{equation}
h_{i}^{\vee }=l^{i}(\frac{d}{c}-1).  \label{eq: dualCoxeter}
\end{equation}%
Furthermore, (\ref{eq:lambda2}) implies that $\omega =\sum_{1\leq i\leq
r}\omega ^{i}$ with $\omega ^{i}$ the Sugawara Virasoro vector for central
charge $c^{i}=\frac{l^{i}d^{i}}{l^{i}+h_{i}^{\vee }}$ for the simple
component $\hat{\mathfrak{g}}^{i}$. It follows that for each component 
\begin{equation}
\frac{d^{i}}{c^{i}}=\frac{d}{c},  \label{eq: dcdici}
\end{equation}%
so that 
\begin{equation}
\lambda ^{i(2)}=-\frac{2d}{c}\omega ^{i},  \label{eq: lambdai2}
\end{equation}%
for the quadratic Casimir on $\mathfrak{g}^{i}$.

\bigskip

We next show that if $\lambda ^{(4)}\in \mbox{Vir}_{c}$ then $\mathfrak{g}$ must be
simple. Let $L_{n}^{i}$ denote the modes of $\omega ^{i}$ and $%
L_{n}=\sum_{i}L_{n}^{i}$ denote the modes of $\omega $ with $%
[L_{m}^{i},L_{n}^{j}]=0$ for $i\neq j$. Using $\lambda
^{(4)}=\sum_{i}\lambda ^{i(4)}$ (for quadratic Casimirs on $\mathfrak{g}^{i}$%
) it follows from (\ref{eq:LMlambda}) that 
\begin{equation}
L_{2}^{i}\lambda ^{(4)}=3\lambda ^{i(2)}.  \label{eq: L2lambda4}
\end{equation}%
Since $L_{n}^{i}$ satisfies the Virasoro algebra of central charge $c^{i}$
we find 
\begin{eqnarray*}
L_{2}^{i}L_{-2}^{2}\mathbf{1} &=&8\omega ^{i}+c^{i}\omega , \\
L_{2}^{i}L_{-4}\mathbf{1} &=&6\omega ^{i}.
\end{eqnarray*}%
If\ $\lambda ^{(4)}\in \mbox{Vir}_{c}$ then (\ref{eq:lambda4}) holds and hence 
\begin{equation*}
L_{2}^{i}\lambda ^{(4)}=-\frac{3d}{c\left( 5c+22\right) }\left(
(44+6c)\omega ^{i}+4c^{i}\omega \right) .
\end{equation*}%
Equating to (\ref{eq: L2lambda4}) and using (\ref{eq: lambdai2}) implies
that 
\begin{equation*}
\omega ^{i}=\frac{c^{i}}{c}\omega .
\end{equation*}%
But since the Virasoro vectors $\omega ^{1},\ldots \omega ^{r}$ are
independent it follows that $r=1$ i.e. $\mathfrak{g}$ is a simple Lie
algebra.

\medskip

If (\ref{eq:lambda4}) holds one also finds that 
\begin{equation*}
\langle a,o(\lambda ^{(4)})b\rangle =-\frac{9d(6+c)}{c\left( 5c+22\right) }%
\langle a,b\rangle .
\end{equation*}%
Comparing to the corresponding term in (\ref{eq:Fablambda}) this results in
a further constraint on the parameters $d,c$ in (\ref{eq: fxycd}) given by%
\begin{equation}
d=\frac{c\left( 5c+22\right) }{10-c}.  \label{eq: d(c)}
\end{equation}%
Notice that the numerator vanishes for $c=0,-22/5$, the zeros of the Kac
determinant $\det M_{4}(c)$ (\ref{eq: detM4}).

\bigskip

For integral $d>0$ there are only 42 rational values of $c$ satisfying (\ref%
{eq: d(c)}). This list is further restricted by the possible values of $d$
for $\mathfrak{g}$ simple. The level $l$ is necessarily rational from (\ref%
{eq: dualCoxeter}). Restricting $l$ to be integral (for example, if $V$ is
assumed to be $C_{2}$-cofinite \cite{DM1}) we find that $l=1$ and $\mathfrak{%
g}$ must be one of Deligne's exceptional Lie algebras:

\begin{theorem}
\label{theorem:Deligne} Suppose $\lambda ^{(4)}\in \mbox{Vir}_{c}$.

(a) Then $\mathfrak{g}$ is a simple Lie algebra.

(b) If $c$ is rational and the level $l$ of $\hat{\mathfrak{g}}$ is integral
then 
\begin{equation*}
\mathfrak{g}=A_{1},A_{2},G_{2},D_{4},F_{4},E_{6},E_{7}\text{ or }E_{8},
\end{equation*}%
with dual Coxeter number 
\begin{equation*}
h^{\vee }=\frac{d}{c}-1=\frac{12+6c}{10-c},
\end{equation*}%
for central charge $c=1,2,\frac{14}{5},4,\frac{26}{5},6,7,8$ respectively
and level $l=1$.
\end{theorem}

The simple Lie algebras appearing in Theorem \ref{theorem:Deligne} are known
as Deligne's exceptional Lie algebras \cite{De}. These algebras are of
particular interest because not only is the dimension $d$ of the adjoint
representation $\mathfrak{g}$ described by a rational function of $c$ in (%
\ref{eq: d(c)}) but also the dimensions of the irreducible representations
that arise in decomposition of up to four tensor products of $\mathfrak{g}$.
In Deligne's original calculations, these dimensions were expressed as
rational functions of a convenient parameter $\lambda $. In this VOA setting
we instead employ the canonical parameter $c$ where 
$\lambda =\frac{c-10}{2+c}.$
\bigskip

\begin{exercise}
\label{Exercise KacMoody} Verify (\ref{eq: KacMoody}).
\end{exercise}

\begin{exercise}
\label{Exercise lambda01} Show that $\lambda^{(0)}=-d\mathbf{1}$ and $\lambda ^{(1)}=0$.
\end{exercise}

\begin{exercise}
\label{Exercise lambda(2)} Verify (\ref{eq:lambda2}).
\end{exercise}

\begin{exercise}
\label{Exercise lambda(4)} Verify (\ref{eq:lambda4}) using (\ref{eq M4}).
\end{exercise}

\bigskip
\subsection{Genus One Constraints from Quadratic Casimirs}
\label{Subsect_Genus1 Casimir}

We next consider the constraints on the genus one partition function $%
Z_{V}(\tau )$ that follow if $\lambda ^{(4)}\in \mbox{Vir}_{c}$. We will show that
in this case, $Z_{V}(\tau )$ is the unique solution to a second order
Modular Linear Differential Equation (MLDE) (cf. Subsection \ref{Subsect_ExamplesVVMF}). As a
consequence, we prove that $V=L_{\mathfrak{g}}(1,0)$, the level 1 WZW VOA
where $\mathfrak{g}$ is an Deligne exceptional series. To prove this we
apply both versions of Zhu's recursion formulas of Theorems \ref{TheoremZhu} and \ref{Proposition_1pt}. In particular, we evaluate the $1$-point correlation function for a
Virasoro descendent of the vacuum from where an MLDE naturally arises. This is
similar in spirit to Zhu's \cite{Z} analysis of correlation functions for
the modules of $C_{2}$-cofinite VOAs but has the advantage of being
considerably less technical.

\bigskip

We recall the genus one partition function%
\begin{equation*}
Z_{V}(\tau )=\mathrm{Tr}_{V}(q^{L_{0}-c/24}),
\end{equation*}%
the 1-point correlation function for $a\in V$ 
\begin{equation}
Z_{V}(a,\tau )=\mathrm{Tr}_{V} o(a)q^{L_{0}-c/24},  \label{eq:1pt}
\end{equation}%
and the 2-point correlation function which can be expressed in terms of
1-point functions by%
\begin{eqnarray}
F_{V}((a,x),(b,y),\tau ) &=&Z_{V}(Y[a,x]Y[b,y]\mathbf{1},\tau )
\label{eq:2pt1} \\
&=&Z_{V}(Y[a,x-y]b,\tau ),  \label{eq:2pt2}
\end{eqnarray}%
for square bracket vertex operators $Y[a,z]=Y(q_{z}^{L_{0}}a,q_{z}-1).$

\bigskip

We define quadratic Casimir vectors in the square bracket VOA formalism%
\begin{equation*}
\lambda ^{\lbrack n]}=u^{\alpha }[1-n]\bar{u}^{\alpha }\in V_{[n]},
\end{equation*}%
(for $\alpha $ summed) for basis $\{u^{\alpha }\}$ and square bracket LiZ
dual basis $\{\bar{u}^{\alpha }\}$. Consider the genus one analogue of (\ref%
{eq: Fabdef}) given by the 2-point function 
\begin{equation*}
F_{V}((u^{\alpha },x),(\bar{u}^{\alpha },y),\tau )=Z_{V}(Y[u^{\alpha },x]Y[%
\bar{u}^{\alpha },y]\mathbf{1},\tau ),
\end{equation*}%
($\alpha $ summed). Associativity (\ref{eq:2pt2}) implies the genus one
analogue of (\ref{eq:Fablambda}) so that 
\begin{equation}
F_{V}((u^{\alpha },x),(\bar{u}^{\alpha },y),\tau )=\sum_{n\geq
0}Z_{V}(\lambda ^{\lbrack n]},\tau )(x-y)^{n-2}.  \label{eq:Fablambdatorus}
\end{equation}

From Zhu's recursion formula I of Theorem \ref{TheoremZhu} we may alternatively expand 
$F((u^{\alpha },x),(\bar{u}^{\alpha },y),\tau )$ in terms of Weierstrass
functions as follows%
\begin{eqnarray*}
F_{V}((u^{\alpha },x),(\bar{u}^{\alpha },y),\tau ) &=&\mathrm{Tr}_{V}\left(
o(u^{\alpha })o(\bar{u}^{\alpha })q^{L_{0}-c/24}\right) \\
&&+\sum\limits_{m\geq 1}\frac{(-1)^{m+1}}{m!}P_{1}^{(m)}(x-y,\tau
)Z_{V}(u^{\alpha }[m]\bar{u}^{\alpha },\tau ) \\
&=&\mathrm{Tr}_{V}\left( o(u^{\alpha })o(\bar{u}^{\alpha
})q^{L_{0}-c/24}\right) -dP_{2}(x-y,\tau )Z_{V}(\tau ).
\end{eqnarray*}%
Recalling Theorem \ref{thm5.1} we may compare the $(x-y)^{2}$ terms in this expression and
(\ref{eq:Fablambdatorus}) to obtain 
\begin{equation}
Z_{V}(\lambda ^{\lbrack 4]},\tau )=-3dE_{4}(\tau )Z_{V}(\tau ).
\label{eq:Zlambda4}
\end{equation}

\medskip

Since $(V,Y(\ ),\mathbf{1},\omega )$ is isomorphic to $(V,Y[\ ],\mathbf{1},\tilde{\omega})$ it follows that $\lambda ^{(n)}\in \mbox{Vir}_{c}$ iff 
$\lambda^{\lbrack n]}\in \mbox{Vir}_{c}$. Thus assuming $\lambda ^{(4)}\in \mbox{Vir}_{c}$ we have 
\begin{equation}
Z_{V}(\lambda ^{\lbrack 4]},\tau )=-\frac{3d}{c(5c+22)}\left[
4Z_{V}(L[-2]^{2}\mathbf{1},\tau )+(2+c)Z_{V}(L[-4]\mathbf{1},\tau )\right] ,
\label{eq: Zlambda41pt}
\end{equation}%
following (\ref{eq:lambda4}). The Virasoro 1-point functions $Z_{V}(L[-2]^{2}%
\mathbf{1},\tau )$, $Z_{V}(L[-4]\mathbf{1},\tau )$ can be evaluated via the
Zhu recursion formula II of Theorem \ref{Proposition_1pt}. In particular taking 
$u=\tilde{\omega}$ and $v$ of $L[0]$ weight $k$ in (\ref{1ptrecur}) we obtain the general
Virasoro recursion formula 
\begin{eqnarray}
Z_{V}(L[-n]v,\tau ) &=&\delta _{n,2}\mathrm{Tr}_{V}(o(\tilde{\omega}%
)o(v)q^{L_{0}-c/24})  \notag \\
&&+\sum\limits_{0\leq m\leq k}(-1)^{m}\binom{m+n-1}{m+1}E_{m+n}(\tau
)Z_{V}(L[m]v,\tau ).  \notag \\
&&  \label{eq: ZhuL}
\end{eqnarray}%
But $o(\tilde{\omega})=L_{0}-c/24$ and hence 
\begin{equation*}
\mathrm{Tr}_{V}(o(\tilde{\omega})o(v)q^{L_{0}-c/24})=\theta Z_{V}(v,\tau ),
\end{equation*}%
where $\theta =q\frac{d}{dq}$. It follows that 
\begin{equation}
Z_{V}(L[-2]v,\tau )=D_{k}Z_{V}(v,\tau )+\sum\limits_{2\leq m\leq
k}E_{2+m}(\tau )Z_{V}(L[m]v,\tau ),  \label{eq: ZL2}
\end{equation}%
where $D_{k}=\theta +kE_{2}(\tau )$ is the modular derivative (\ref{ModDer}).
(Zhu makes extensive use of the identities (\ref{eq: ZhuL}) and (\ref{eq: ZL2}) 
in his analysis of correlation functions for $C_{2}$-cofinite VOAs \cite{Z}. This is the origin of MLDEs as discussed in Section \ref{Section_VOAModInvariance}).

\medskip
We  immediately  find  from (\ref{eq: ZhuL}) that $Z_{V}(L[-4]\mathbf{1},\tau )=0$ and 
\begin{eqnarray*}
Z_{V}(L[-2]^{2}\mathbf{1},\tau ) &=&D_{2}Z_{V}(L[-2]\mathbf{1},\tau
)+E_{4}(\tau )Z_{V}(L[2]L[-2]\mathbf{1},\tau ) \\
&=&\left[ D^{2}+\frac{c}{2}E_{4}(\tau )\right] Z_{V}(\tau ),
\end{eqnarray*}%
where $D^{2}=D_{2}D_{0}=\left( q\frac{d}{dq}\right) ^{2}+2E_{2}(\tau )q\frac{%
d}{dq}$. Substituting into (\ref{eq: Zlambda41pt}) we find $Z_{V}(\tau )$
satisfies the following second order MLDE: 
\begin{equation}
\left[ D^{2}-\frac{5}{4}c(c+4)E_{4}(\tau )\right] Z_{V}(\tau )=0.
\label{eq:Zdeqn}
\end{equation}%
(\ref{eq:Zdeqn}) has a regular singular point at $q=0$ with indicial roots $%
-c/24$ and $(c+4)/24$. Applying (\ref{eq: d(c)}) it follows  that  there exists a unique solution with leading $q$
expansion $Z_{V}(\tau )=q^{-c/24}(1+O(q))$. Furthermore, since $E_{4}(\tau )$
is holomorphic then $Z_{V}(\tau )$ is also holomorphic for $0<|q|<1$. In
summary, we find

\begin{theorem}
\label{theorem:Zdeqn} If $\lambda ^{(4)}\in \mbox{Vir}_{c}$ then $Z_{V}(\tau )$ is
a uniquely determined holomorphic function in $\mathfrak{H}$.
\end{theorem}

\medskip
An immediate consequence of Theorems \ref{theorem:Deligne} and \ref{theorem:Zdeqn} is:

\begin{theorem}
\label{theorem:ZWZW} $V=L_{\mathfrak{g}}(1,0)$ the level one WZW model
generated by $\mathfrak{g}$.
\end{theorem}

\noindent \textbf{Proof.} Clearly $L_{\mathfrak{g}}(1,0)\subseteq V$ with $\omega
,\lambda ^{(2)},\lambda ^{(4)}\in L_{\mathfrak{g}}(1,0)$. Thus\ $L_{%
\mathfrak{g}}(1,0)$ satisfies the conditions of Theorem \ref{theorem:Zdeqn}
for the same central charge $c$. Hence $Z_{L_{\mathfrak{g}}(1,0)}(\tau
)=Z_{V}(\tau )$ and so $L_{\mathfrak{g}}(1,0)=V$. 

\bigskip

It is straightforward to substitute $Z(\tau )=q^{-c/24}\sum_{n}\dim
V_{n}q^{n}$ into (\ref{eq:Zdeqn}) and solve recursively for $\dim V_{n}$ as
a rational function in $c$. In this way we recover $\dim V_{1}=d$ of (\ref{eq: d(c)}). The next two terms are 
\begin{eqnarray*}
\dim V_{2} &=&\frac{c(804+508c+175c^{2}+25c^{3})}{2(22-c)(10-c)}, \\
\dim V_{3} &=&\frac{c(33344+148872c+68308c^{2}+10330c^{3}+975c^{4}+125c^{5})%
}{6(34-c)(22-c)(10-c)}.
\end{eqnarray*}

\medskip

These dimension formulas can be further refined as follows. Consider the
Virasoro decomposition of $V_{2}$: 
\begin{equation}
V_{2}=\mathbb{C}\omega \oplus L_{-1}\mathfrak{g}\oplus P_{2},  \label{eq:V_2}
\end{equation}%
where $P_{2}$ is the space of weight two primary vectors. Let $p_{2}=\dim
P_{2}$. Then $\dim V_{2}=1+d+p_{2}$ with 
\begin{equation}
p_{2}=\frac{5(5c+22)(c+2)^{2}(c-1)}{2(22-c)(10-c)}.  \label{eq: p2dim}
\end{equation}%
Comparing with Deligne's analysis of the irreducible decomposition of tensor
products of $\mathfrak{g}$ we find that 
\begin{equation*}
p_{2}=\dim Y_{2}^{\ast },
\end{equation*}%
where $Y_{2}^{\ast }$ denotes an irreducible representation of $\mathfrak{g}$
in Deligne's notation \cite{De}. This is explored further in \cite{T3}.
\bigskip

Similarly for $V_{3}$ we find 
\begin{equation*}
V_{3}=\mathbb{C[}L_{-1}\omega ]\oplus L_{-1}^{2}\mathfrak{g}\oplus L_{-2}%
\mathfrak{g}\oplus L_{-1}P_{2}\oplus P_{3},
\end{equation*}%
where $P_{3}$ is the space of weight three primary vectors. Let $p_{3}=\dim
P_{3}$. Then $\dim V_{3}=1+2d+p_{2}+p_{3}$ with 
\begin{eqnarray*}
p_{3} &=&\frac{5c(5c+22)(c-1)(c+5)(5c^{2}+268)}{6(34-c)(22-c)(10-c)} \\
&=&\dim X_{2}+\dim Y_{3}^{\ast },
\end{eqnarray*}%
where $X_{2},Y_{3}^{\ast }$ denote two other irreducible representations of $%
\mathfrak{g}$ in Deligne's notation of dimension 
\begin{eqnarray*}
\dim X_{2} &=&\frac{5c(5c+22)(c+6)(c-1)}{2(10-c)^{2}} \\
\dim Y_{3}^{\ast } &=&\frac{5c(5c+22)(c+2)^{2}(c-8)(5c-2)(c-1)}{%
6(10-c)^{2}(22-c)(34-c)}.
\end{eqnarray*}

\subsection{Higher Weight Constructions}
\label{Subsect_HigherWt}

We can generalize the arguments given above to consider a VOA $V$ with $\dim
V_{1}=0$. Here we construct Casimir vectors from the weight two primary
space $P_{2}$ (provided $\dim P_{2}>0$) and obtain constraints on $V$ that
follow from such Casimirs being Virasoro vacuum descendents. If $\dim
P_{2}=0 $ we consider primaries of weight 3 and so on. In general, let $V$
be a VOA with primary vector space $P_{K}$ of lowest weight $K$ i.e. $%
V_{n}=(\mbox{Vir}_{c})_{n}$ for all $n<K$ so that%
\begin{equation}
Z_{V}(\tau )=q^{-c/24}(\sum\limits_{n<K}\dim (\mbox{Vir}_{c})_{n}q^{n}+O(q^{K})).
\label{eq: ZVK}
\end{equation}%
(Recall (\ref{vircqdim}) that $\sum\limits_{n\geq 0}\dim
(\mbox{Vir}_{c})_{n}q^{n}=\prod\limits_{m\geq
2}(1-q^{m})^{-1}=1+q^{2}+q^{3}+2q^{3}+\ldots $.) We construct Casimir
vectors, as in (\ref{eq:lambdan}), from a $P_{k}$ basis $\{u^{\alpha }\}$
and LiZ dual basis $\{\bar{u}^{\alpha }\}$ 
\begin{equation*}
\lambda ^{(n)}=u_{2K-1-n}^{\alpha }\bar{u}^{\alpha }\in V_{n}.
\end{equation*}%
We find the following natural generalization of Theorems \ref{theorem:Zdeqn}
and \ref{theorem:ZWZW}:

\begin{theorem}
\label{theorem:Zdeqngen}Let $V$ be a VOA with primary vectors of lowest
weight $K=2$ or $3$. If $\lambda^{(2K+2)}\in \mbox{Vir}_{c}$ then

(a) $Z_{V}(\tau )$ of (\ref{eq: ZVK}) is a holomorphic function in $\mathfrak{H}$
and is the unique solution to a MLDE of order $K+1$.

(b) $V$ is generated by $P_{K}$. 
\end{theorem}

\begin{remark}
\label{Remark}We conjecture that Theorem \ref{theorem:Zdeqngen} holds for
all $K$. 
\end{remark}

For $K=2$ the elements of $P_{2}$ satisfy a commutative nonassociative
algebra with invariant (LiZ) form known as a Griess algebra (cf. Exercise \ref{Exercise_9.4.3}). Theorem \ref{theorem:Zdeqngen} implies the dimension of the Griess algebra is 
\begin{equation}
\dim P_{2}=\frac{1}{2}\frac{(5c+22)(2c-1)(7c+68)}{c^{2}-55c+748}.
\label{eq: Griess}
\end{equation}%
This result originally appeared in \cite{Mat} subject to stronger
assumptions. Following Remark \ref{Remark Virc} we note that the zeros of
the numerator are the zeros $c_{p,q}$ of the Kac determinant $\det M_{n}(c)$
for $n\leq 6$. There are 37 rational values of $c$ for which $\dim P_{2}$ is
a positive integer. Furthermore, we may solve $Z_{V}(\tau )$ iteratively for 
$\dim V_{n}$ as rational functions in $c$. There are 9 values of $c$ for
which $\dim V_{n}$ is a positive integer for $n\leq 400$ given by \cite{T3}: 

\begin{equation*}
\begin{tabular}{|c|c|c|c|c|c|c|c|c|c|}
\hline
$c$ & $-{\frac{44}{5}}$ & $8$ & $16$ & $\frac{{47}}{2}$ & $24$ & $32$ & ${%
\frac{164}{5}}$ & ${\frac{236}{7}}$ & $40$ \\ \hline
$\dim P_{2}$ & $1$ & $155$ & $2295$ & $96255$ & $196883$ & $139503$ & $90117$
& $63365$ & $20619$ \\ \hline
\end{tabular}%
\end{equation*}

The first five cases can all be realized by explicit constructions. Of
particular interest is the case $c=24$ realized by the FLM Moonshine Module $%
V^{\natural }$ with $Z_{V^{\natural }}(\tau )=j(\tau )-744$ for which $P_{2}$ is
the original Griess algebra of dimension 196883 and whose automorphism group is
the Monster group (cf. Subsection \ref{Subsect_Applications}). There are constructions for $c=32$ and $40$ with the
appropriate partition function but it is not known if $\lambda ^{(6)}\in
\mbox{Vir}_{c}$. There are no known constructions for $c={\frac{164}{5}}$ and $%
\frac{236}{7}$.

\bigskip

For $K=3$ we find 
\begin{equation*}
\dim P_{3}={\frac{(5c+22)(2c-1)(7c+68)(5c+3)(3c+46)}{-5{c}^{4}
+703{c}^{3}-32992{c}^{2}+517172c-3984},}
\end{equation*}%
where the zeros of the numerator are Kac determinant zeros $c_{p,q}$ for $(p-1)(q-1)=n\leq 8$. Iteratively solving the appropriate MLDE
for $Z_{V}(\tau )$ we find $\dim P_{3}$ and $\dim V_{n}$ are positive
integral for only 3 rational values of $c$

\begin{equation*}
\begin{tabular}{|c|c|c|c|}
\hline
$c$ & $-{\frac{114}{7}}$ & $\frac{4}{5}$ & $48$ \\ \hline
$\dim P_{3}$ & $1$ & $1$ & $42987519$ \\ \hline
\end{tabular}%
\end{equation*}%

The first two examples can be realized by known VOAs. For $c=48$ we find $%
Z_{V}(\tau )=J(\tau )^{2}-393767$ which, intruigingly, is the partition
function of the minimal holomorphic VOA $V^{(2)}$ briefly discussed  in Subsection \ref{Subsect_Applications}.

\section{Appendices}\label{Section_Appendices}

\subsection{Lie Algebras and Representations}\label{SubSec_LieAppendix}
An \emph{associative algebra} is a linear space $A$ equipped with a bilinear, associative  product 
$A \otimes A \rightarrow A$, denoted by juxtaposition. Thus $a \otimes b \mapsto ab$ and
\begin{eqnarray*}
(ab)c = a(bc).
\end{eqnarray*}

A \emph{Lie algebra} is a linear space $L$ equipped with a bilinear product
(usually called bracket) $[ \ \ ]: L \otimes L \rightarrow L$ such that 
\begin{eqnarray*}
&&[ab]  =  -[ba] \hspace{3.0cm} (\mbox{skew-commutativity}) \\
&&{[}a[bc{]]} + [b[ca]] + [c[ab]] = 0 \ \ (\mbox{Jacobi identity})
\end{eqnarray*}

An associative algebra $A$ gives rise to a Lie algebra $A^-$ on the \emph{same} linear space
by defining $[ab] = ab - ba$. A basic example is End$(V)$ for a linear space $V$, where the  
associative product is composition of endomorphisms. This situation can be exploited using another basic associative algebra, the \emph{tensor algebra} 
\begin{eqnarray*}
T(V) = \oplus_{n \geq 0} V^{\otimes n} = \mathbb{C} \oplus V \oplus V\otimes V \oplus \hdots
\end{eqnarray*}
over $V$. Let $\iota:V \rightarrow T(V)$ be canonical identification of $V$ with the 
degree $1$ piece of $T(V)$. The \emph{universal mapping property} (UMP) for tensor algebras says that
any linear map $f: V \rightarrow A$  into an associative algebra $A$ has a 
\emph{unique} extension to a morphism of associative algebras $\alpha: T(V) \rightarrow A$:  
\begin{eqnarray*}
&&V \  \stackrel{f}{\longrightarrow} A \\
&& \  \iota \searrow \ \  \ \  \uparrow \alpha\\
&& \hspace{1.00cm} T(V)
\end{eqnarray*}
with $f = \alpha \circ \iota$.

 \medskip
A \emph{representation} of a Lie algebra $L$ is a linear map
$\pi: L \rightarrow$ End$(V)$ for some $V$ such that
\begin{eqnarray*}
\pi([ab])  = \pi(a)\pi(b) - \pi(b)\pi(a).
\end{eqnarray*}
That is, $\pi: L \rightarrow$ End$(V)^-$ is a morphism of Lie algebras. 
We call $V$ an $L$-module. 

\medskip
UMP provides an extension of $\pi$ to a
morphism of associative algebras $\alpha: T(L) \rightarrow$ End$(V)$. Identifying $a \in L$
with its image in $T(L)$, we see that for $a, b \in L$
\begin{eqnarray*}
\alpha(a\otimes b - b \otimes a - [ab]) &=&\alpha(a)\alpha(b) - \alpha(b)\alpha(a) - \alpha([ab])\\
&=&\pi(a)\pi(b) - \pi(b)\pi(a) - \pi([ab])\\
&=& 0.
\end{eqnarray*}
Let $J \subseteq T(L)$ be the $2$-sided ideal generated by $a\otimes b - b \otimes a - [ab], \ a, b \in L$, and set
\begin{eqnarray*}
\mathcal{U}(L) = T(L)/J.
\end{eqnarray*}
This is the \emph{universal enveloping algebra} of $L$. Thus every representation
$\pi$ of $L$ extends to a representation of the universal enveloping algebra in a canonical way:
\begin{eqnarray*}
&&L \stackrel{\pi}{\longrightarrow}\mbox{End}(V) \\
&& \   \iota' \searrow  \ \ \ \ \  \uparrow \alpha\\
&& \hspace{1.40cm} \mathcal{U}(L)
\end{eqnarray*}
where $\iota'$ is the composition $L \stackrel{\iota}{\rightarrow} T(L) \rightarrow \mathcal{U}(L)$

\begin{theorem}(PBW Theorem). Fix an ordered basis $x_1, x_2, \hdots$ of $L$, with
$\bar{x_i}$ the image of $x_i$ in $\mathcal{U}(L)$. Then
\begin{eqnarray*}
\{ \bar{x}_{i_1} \bar{x}_{i_2} \hdots \bar{x}_{i_k} \ | \ i_1 \geq i_2 \geq \hdots \geq i_k\geq 1 \}
\end{eqnarray*}
is a basis for $\mathcal{U}(L)$.
\end{theorem}
From PBW we see that $\iota'$
is \emph{injective}. Then for a representation of $\mathcal{U}(L)$, restriction to the subspace
$L = \iota(L)$ furnishes a representation of $L$. In this way, representations of $L$ and $\mathcal{U}(L)$ determine each other in a canonical fashion - a statement that can be better stated using categories of modules.

\medskip
The Lie algebra $L$ has a \emph{triangular decomposition} if it decomposes as
\begin{eqnarray*}
L = L^+ \oplus L^0 \oplus L^-
\end{eqnarray*}
\begin{flushright}

\end{flushright}
such that $L^{\pm}, L^0$ are Lie subalgebras, and the bracket satisfies
\begin{eqnarray*}
[L^+ L^-] \subseteq L^0, \ \ [L^{\pm} L^0] \subseteq L^{\pm}.
\end{eqnarray*}
Use of PBW and an appropriate choice of (ordered) basis leads to an identification
\begin{eqnarray*}
\mathcal{U}(L) = \mathcal{U}(L^-) \otimes \mathcal{U}(L^0) \otimes \mathcal{U}(L^+).
\end{eqnarray*}

\medskip
Noting that $L^0 \oplus L^+ \subseteq L$ is a Lie subalgebra, let $\pi: L^0 \oplus L^+ \rightarrow$ End$(V)$ be a representation. The \emph{induced module}
is
\begin{eqnarray}\label{indmoddef}
\mbox{Ind}(V) = \mbox{Ind}_{\mathcal{U}(L^0\oplus L^+)}^{\mathcal{U}(L)} V := \mathcal{U}(L) \otimes_
{\mathcal{U}(L^0\oplus L^+)} V = \mathcal{U}(L^-) \otimes V.
\end{eqnarray}
It is a $\mathcal{U}(L)$-module, hence also an $L$-module upon restriction. A ubiquitous special case occurs when $V$ is an $L^0$-module, which then becomes an $L^0\oplus L^+$-module 
by letting $L^+$ \emph{annihilate} $V$.

\bigskip

\begin{exercise}\label{Exercise_10.1} Show that the following Lie algebras have natural triangular decompositions:\\
(a) Heisenberg algebra $\hat{A}$ with 
\begin{eqnarray*}
\hat{A}^+ = \oplus_{n > 0} \mathbb{C}a\otimes t^n, \ \hat{A}^- = \oplus_{n < 0} \mathbb{C}a \otimes t^n, \ \hat{A}^0 = \mathbb{C}a\otimes t^0 \oplus \mathbb{C}K.
\end{eqnarray*}
(b) Virasoro algebra $\mbox{Vir}$ with
\begin{eqnarray*}
\mbox{Vir}^+ = \oplus_{n > 0} \mathbb{C}L_n, \ \mbox{Vir}^- = \oplus_{n < 0} \mathbb{C}L_n, \ \mbox{Vir}^0 = \mathbb{C}L_0 \oplus \mathbb{C}K.
\end{eqnarray*}
(c) Finite-dimensional simple Lie algebra (equipped with a choice of Cartan subalgebra and root system) with
$L^+ = \{ \mbox{positive root spaces}\}, \ L^- = \{\mbox{negative root spaces}\}, \ L^0 = \{\mbox{Cartan subalgebra} \}$.
\end{exercise}

\subsection{The Square Bracket Formalism}\label{SubSec_SquareBracket}
We prove (\ref{square1})-(\ref{square4}) of  Subsection \ref{Subsect_SquareVOA}. The square bracket vertex operator (\ref{Ysquare}), 
(\ref{qopformalism})
is
\begin{equation*}
Y[v,z]=q_{z}^{wt(v)}Y(v,q_{z}-1).
\end{equation*}%
Thus the square bracket modes of $Y[v,z]=\sum_{m\in \mathbb{Z}}v[m]z^{-m-1}$
are given by%
\begin{eqnarray*}
v[m] &=&\mathrm{Res}_z Y(v,q_{z}-1)z^{m}q_{z}^{wt(v)} \\
&=&\mathrm{Res}_z Y(v,q_{z}-1)\frac{d}{dz}%
(q_{z}-1)z^{m}q_{z}^{wt(v)-1}.
\end{eqnarray*}%
We may rewrite this in terms of $w=q_{z}-1=z+O(z^{2})$ by means of a
(formal) chain rule \cite{FHL}, \cite{Z} so that 
\begin{eqnarray*}
v[m] &=&\mathrm{Res}_w Y(v,w)z(w)^{m}q_{z(w)}^{wt(v)-1} \\
&=& \mathrm{Res}_w  Y(v,w)\ln (1+w)^{m}(1+w)^{wt(v)-1}.
\end{eqnarray*}%
Defining $c(wt(v),i,m)$ for $i\geq m\geq 0$ by 
\begin{equation*}
\sum\limits_{i\geq m}c(wt(v),i,m)w^{i}=\frac{1}{m!}\ln
(1+w)^{m}(1+w)^{wt(v)-1},
\end{equation*}%
we obtain (\ref{square1}). 

\medskip
Next note that $\sum\limits_{m\geq 0}\frac{1}{%
m!}\ln (1+w)^{m}x^{m}=(1+w)^{x}$. Hence we find%
\begin{equation*}
\sum\limits_{i\geq
0}\sum\limits_{m=0}^{i}c(wt(v),i,m)w^{i}x^{m}=(1+w)^{wt(v)-1+x},
\end{equation*}%
from which (\ref{square2}) follows. Finally, 
\begin{eqnarray*}
\sum\limits_{m\geq 0}\frac{(k+1-wt(v))^{m}}{m!}v[m] &=&\sum\limits_{m\geq
0}\sum\limits_{i\geq m}c(wt(v),i,m)(k+1-wt(v))^{m}v_i \\
&=&\sum\limits_{i\geq 0}v_i\sum\limits_{m=0}^{i}c(wt(v),i,m)(k+1-wt(v))^{m},
\\
&=&\sum_{i\geq 0}\binom{k}{i}v_i.
\end{eqnarray*}%
giving (\ref{square4}).

 \bigskip
 \noindent
 Authors' addresses:\\
 Geoffrey Mason: Department of Mathematics,  University of California at Santa Cruz, 
CA 95064, gem@cats.ucsc.edu\\
Michael Tuite: Department of Applied Mathematics, National University of Ireland, Galway, University Road, Galway, Ireland, michael.tuite@nuigalway.ie

\end{document}